\input amstex
\documentstyle{amsppt}
\font\cmssqi=cmssqi8 at 7 pt

\TagsOnRight  
\def\vs{\vskip.3cm}
\def\p{\partial}

\def\no{\noindent}  
\def\io{{\infty}} 
\def\re{\operatorname{Re}}
\def\im{\operatorname{Im}} 
\def\Id{\operatorname{Id}} 
\def\moo{C^{\io}}
\def\mooc{C^{\io}_{\text c}}

\def\N{\Bbb N}

\def\R{\Bbb R}

\def\carno#1#2{\left\Vert#1\right\Vert^2_{#2}}
 
\def\poscal#1#2{\langle#1,#2\rangle}

\def\poi#1#2{\left\{#1,#2\right\}}
\def\norm#1{\left\Vert#1\right\Vert}
\def\nuorm#1{\Vert#1\Vert}
\def\Val#1{\left\vert#1\right\vert}
\def\val#1{\vert#1\vert}

\def\l2{L^2(\R^{n})}
\def\L2{L^2(\R^{2n})}
\def\hs{\hskip15pt}
\def\tab{\leaders\hbox to 3mm{\hfil.\hfil}\hfill}
\def\supp{\operatorname{supp}}

\def\op#1{{\text{Op}(#1)}}
\def\w#1{{#1^{\text {Wick}}}}

\def\RZ{\R^{2n}}

\def\sign{\operatorname{sign}}

\def\cps{condition $(\psi)$}
\def\tab{\leaders\hbox to 1.5mm{\hfil.\hfil}\hfill}
\let \dis=\displaystyle
\let\no=\noindent
\headline
{\ifodd\pageno\rightheadline 
\else\leftheadline\fi}
\def\rightheadline{{\cmssqi
}\hfill\rm\folio}
\def\leftheadline{\folio\hfill{
\cmssqi Condition $(\psi)$}}
\voffset -0.6 truecm
\hoffset 0truecm  
\pageheight{25 truecm}
\pagewidth{16.5 truecm}
\nopagenumbers
\topmatter
\title
Cutting the loss of derivatives
for solvability under condition $(\Psi)$
\endtitle
\date November 4, 2005\enddate
\abstract
For 
a
principal type pseudodifferential operator,
we  prove that 
  \cps\
implies local solvability with a loss of 3/2
derivatives.
We use many elements of Dencker's paper on the proof of the Nirenberg-Treves conjecture
and we provide some improvements of the key energy estimates
which allows us to cut
the loss of derivatives from $\epsilon +3/2$ for any $\epsilon>0$
(Dencker's most recent result)  to 3/2 (the present paper).
It is already known that \cps\ does {\it not} 
imply local solvability with a loss of 1 derivative, so we have
to content ourselves with a loss $>1$.\vs
\centerline{\smc Contents}\vs
{\baselineskip=0,75\normalbaselineskip
\hs { \smc 1. Introduction and statement of the result}\tab{1}\par
\hs \hs
1.1. Introduction\par
\hs \hs 1.2.  Statement of the result\par
\hs \hs 1.3.  Some notations\par
\hs \hs 1.4.  Partitions of unity\par
\par
\hs {\smc 2. The geometry of \cps}\tab 4\par
\hs \hs 2.1.  The basic structure\par
\hs \hs 2.2.  Some lemmas on $\scriptstyle {C^3}$ functions\par
\hs \hs 2.3.  Inequalities for symbols\par
\hs \hs 2.4.  Quasi-convexity\par
\par
\hs {\smc 3. Energy estimates}\tab{18}\par
\hs \hs 3.1.  Preliminaries\par
\par
\hs \hs 3.2.  Stationary estimates for the model cases \par
\hs \hs 3.3.  Stationary estimates \par
\hs \hs 3.4.  The multiplier method \par
\hs {\smc 4. From semi-classical to local estimates}\tab{25}\par
\hs \hs  4.1. From semi-classical to inhomogeneous estimates\par
\hs \hs 4.2. From  semi-classical to localized inhomogeneous estimates\par
\hs \hs 4.3. From inhomogeneous localization to homogeneous localization\par
\hs \hs 4.4. Proof of the solvability result stated in Theorem 1.2.2\par
\par
\hs {\smc A. Appendix}\tab{37}\par
\hs \hs A.1. The Wick quantization\par
\hs \hs A.2. Properties of some metrics\par
\hs \hs A.3. Proof of Lemma 2.1.5 on the proper class\par
\hs \hs A.4. Some a priori estimates and loss of derivatives\par
\hs \hs A.5. Some lemmas on symbolic calculus\par
\hs \hs A.6. The Beals-Fefferman reduction\par
\hs \hs A.7. On tensor products of homogeneous functions\par
\hs \hs A.8. Composition of symbols\par
}
\endabstract
\author
Nicolas Lerner
\endauthor
\affil University  of  Rennes\endaffil
\address
Universit\'e de Rennes 1, Irmar, 
Campus de Beaulieu, 35042 Rennes cedex, France
\endaddress
\email
nicolas.lerner\@univ-rennes1.fr
\vskip0pt
{\it Web-page:} 
http://www.perso.univ-rennes1.fr/nicolas.lerner/
\endemail
\endtopmatter
\document
\head
1. Introduction and statement of the results
\endhead
\subhead
1.1. Introduction
\endsubhead
In 1957, Hans Lewy \cite{Lw}
constructed  a counterexample showing that
very simple and natural differential 
equations can fail to have local solutions;
his example is the complex vector field
$L_{0}=\p_{x_{1}}+i\p_{x_{2}}+i(x_{1}+ix_{2})\p_{x_{3}}$
and one can show that there exists some $\moo$ function $f$
such that the equation $L_{0}u=f$ has no distribution solution,
even locally. A geometric interpretation and a generalization
of this counterexample were given in 1960
by L.H\"ormander in \cite{H2}  
and extended in
\cite{H3}
 to pseudodifferential operators. 
  In 1970,
L.Nirenberg and F.Treves 
(\cite{NT2-NT3-NT4}),
after a study of complex vector fields in \cite{NT1}
(see also \cite{Mi}),
refined this condition
on the principal symbol to the so-called
condition $(\psi)$, 
and  provided strong arguments
suggesting that it should be equivalent to 
local solvability. 
The necessity of condition 
$(\psi)$ for local solvability of  pseudodifferential
equations was proved 
in two dimensions
by
R.Moyer
in
\cite{Mo} and in general by L.H\"ormander (\cite{H5})
in 1981.
The sufficiency of condition $(\psi)$ for 
local solvability of differential equations 
was proved by R.Beals and C.Fefferman
(\cite{BF}) in 1973;
they created a new type of pseudodifferential 
calculus, based on a Calder\'on-Zygmund decomposition, and were able to 
remove the analyticity assumption
required by L.Nirenberg and F.Treves.
For differential equations in any dimension
(\cite{BF})
and for pseudodifferential equations
in two dimensions
(\cite{L1}, see also \cite{L2}), it was shown more precisely 
that $(\psi)$ implies local solvability with 
a loss of one derivative with respect to the
elliptic case:  for a differential
operator  $P$ of order $m$
(or a pseudodifferential operator in two dimensions),
 satisfying condition $(\psi)$, 
$f\in H^s_{\text{loc}}$, 
the equation $Pu=f$ has a solution
$u\in H^{s+m-1}_{\text{loc}}$.
In 1994, it was  proved by N.L. in \cite{L3} 
(see also \cite{H8}, \cite{L8})
that condition
$(\psi)$ does not imply local solvability with loss
of one derivative for pseudodifferential equations,
contradicting 
repeated claims by several authors. 
However in 1996, N.Dencker in \cite{D1}, proved that  these 
counterexamples were indeed locally solvable,  
but with a loss of two derivatives.
\par
In \cite{D2}, N.Dencker claimed 
that he can prove that \cps\
implies
local solvability with loss
of two derivatives;
this preprint 
contains several  breakthrough  ideas
on the control of the second derivatives subsequent
to 
\cps\ and on the choice of the multiplier.
The paper \cite{D3}
contains a
 proof
of local solvability with loss of two derivatives under
\cps,
providing the final step in the proof of the Nirenberg-Treves conjecture;
 the more recent paper
\cite{D4} is providing a proof
of local solvability with loss of
$\epsilon +\frac 3 2$ derivatives under
\cps, for any positive $\epsilon$.
In the present article,
we show that the loss can be limited to 3/2  derivatives, dropping the $\epsilon$ in the previous result. We follow the pattern
of Dencker's paper and give some improvements on the key energy estimates.
\remark{\bf Acknowledgement}{\it
For several months, I have  had the privilege
of exchanging several letters and files with Lars H\"ormander on the topic of solvability.
I am  most grateful for the  help generously provided. These personal communications
are referred to in the text as \cite{H9} and are important in all sections of the present paper.
}\endremark
\subhead
1.2. Statement of the result
\endsubhead
Let $P$ be a properly
supported
principal-type pseudodifferential operator
in a $C^\io$ 
manifold $\Cal M$, with 
principal (complex-valued)
\footnote{Naturally the local solvability of real principal type operators
is also a consequence of the next theorem, but 
much stronger results for real principal type equations
were already established in the 1955 paper \cite{H1}
(see also section 26.1 in \cite{H6}).} symbol $p$.
The symbol $p$
is assumed to be  a $\moo$
homogeneous
\footnote{Here and 
in the sequel, ``homogeneous" will always mean 
positively homogeneous.} 
function
of degree $m$
on
$\dot T^*(\Cal M)$,
the cotangent bundle minus the zero section.
The principal type assumption
 that we shall use here is that 
$$
 (x,\xi)\in
 \dot T^*(\Cal M),\quad
p(x,\xi)=0\Longrightarrow\p_\xi p(x,\xi)\not=0.
\tag 1.2.1$$
Also, the operator $P$
will be assumed of polyhomogeneous type,
which means
that its total symbol
is equivalent
to 
$
p+\sum_{j\ge 1}p_{m-j}
$,
where
$p_k$ is a smooth
homogeneous function of degree $k$
on
$\dot T^*(\Cal M)$.
\definition{Definition 1.2.1. Condition ($\psi$)}
Let $p$ be a $\moo$
homogeneous
function
on  
$\dot T^*(\Cal M)$.
The function $p$ is said to satisfy \cps\
if, for $z=1$ or $i$, 
$\im zp$
does not change sign from $-$ to $+$
along an oriented
bicharacteristic of
$\re zp$.
\enddefinition
It is a non-trivial fact that
\cps\
is
invariant by multiplication
by an complex-valued smooth
elliptic factor
(see section 26.4 in \cite{H6}).
\proclaim{Theorem 1.2.2}
Let $P$
be as above, such that
its principal symbol $p$ satisfies
\cps.
Let $s$ be a real number.
Then,
for all $x\in\Cal M$,
there exists a neighborhood $V$
such that
for all
$f\in H^s_{\text{loc}}$,
there exists 
$u\in {H}^{s+m-\frac{3}{2}}_{\text{loc}}$
such that
$$
Pu=f \quad\text{in $V$}.
$$
 \endproclaim
\demo{The proof of this theorem
will be given at the end of section 4}
\enddemo
Note that our loss of derivatives is equal to 
 3/2. The paper
\cite{L3}
proves that solvability with loss of one derivative
does {\it not} follow from
\cps,
so we have to content ourselves
with a loss strictly greater than one.
However, the number 3/2
is not likely to play
any significant r\^ole
and one should probably
expect
a loss of 1+$\epsilon$ derivatives under \cps.
In fact, for the counterexamples given
in \cite{L3}, it seems
(but it has not been proven) that there is only a 
``logarithmic"
loss, i.e. the solution $u$ should satisfy
$
 u\in \log{\langle D_{x}\rangle}\bigl({H}^{s+m-1}\bigr).
$
\par
Nevertheless,
the methods used in the present article
are strictly limited to providing a 3/2 loss.
We refer the reader to
our appendix A.4 for an  argument
involving a Hilbertian lemma on a simplified model.
This is of course
in sharp contrast
with
operators
satisfying
condition $(P)$
such as differential operators
satisfying
\cps. Let us recall that condition
$(P)$ is simply ruling out
any change of sign
of $\im(zp)$
along the oriented
Hamiltonian flow of
$\re (zp)$.
Under condition $(P)$
(\cite{BF})
or under \cps\
in two dimensions
(\cite{L1}),
local solvability occurs 
with a loss of one derivative,
the ``optimal" loss,
and in fact the same as for
$\p/\p x_1$.
One should also note that
the semi-global
existence theorems
of \cite{H4}
(see also theorem 26.11.2 in
\cite{H6})
involve a loss of 1+$\epsilon$ derivatives.
However in that case there is no known counterexample
which would ensure that this loss is unavoidable.
\remark
{Remark 1.2.3}
 Theorem 1.2.2 will be proved by a multiplier
 method, involving the computation of
 $\poscal{Pu}{Mu}
 $
 with a suitably chosen operator $M$.
It is interesting to notice
that, the greater is the loss of derivatives,
the more regular should be the multiplier in the energy
method.
As a matter of fact, the Nirenberg-Treves multiplier
of \cite{NT3} is not even a pseudodifferential operator
in the $S^0_{1/2,1/2}$ class,
since it could be as singular as the operator
$\sign{D_{x_1}}$; 
 this does not
create any difficulty,
since the loss of derivatives is only 1.
On the other hand,
in \cite{D1},
\cite{L6},
where estimates with loss of 2 derivatives are handled,
the regularity
of the multiplier
is much better than
$S^0_{1/2,1/2}$,
since we need to consider it
as an operator of order 0
in an asymptotic class
defined by an admissible metric on the phase
space.
\endremark
\remark{N.B} For microdifferential operators
acting on microfunctions, the sufficiency
of \cps\
was proven by J.-M.Tr\'epreau \cite{Tr}(see also \cite{H7}),
so the present paper is concerned only with the $\moo$ category.
\endremark
 \subhead 
 1.3. Some notations
 \endsubhead
 First of all, we recall the definition
of the Weyl quantization $a^w$
of a function $a\in\Cal S(\RZ)$:
for $u\in\Cal S(\R^n)$,
$$
(a^w u)(x)=\iint e^{2i\pi(x-y)\xi}a
(\frac{x+y}{2},\xi) u(y) dy.
\tag 1.3.1$$
Our definition of the Fourier transform
$\hat u$
of $u\in\Cal S(\R^n)$
is
$
\hat u (\xi)=\int 
e^{-2i\pi x\xi}u(x) dx
$
and the usual quantization 
$a(x,D_x)$
of $a\in \Cal S(\RZ)$ is 
$
(a(x,D_x) u)(x)=\int e^{2i\pi x\xi}a
({x},\xi) \hat u(\xi) d\xi.
$
The phase space
$\R^n_x\times \R^n_\xi$
is a symplectic vector space
with the standard
symplectic form
$$
\bigl[(x,\xi),(y,\eta)\bigr]=\poscal{\xi}{y}-
\poscal{\eta}{x}.
\tag 1.3.2$$ 
\definition{Definition 1.3.1}
Let $g$ be a metric on $\RZ$, i.e. a mapping $X\mapsto g_X$
from $\RZ$
to the cone of positive definite quadratic forms on $\RZ$.
Let $M$ be a positive function defined on $\RZ$.
\item {(1)} The metric $g$ is said to be slowly varying whenever
$\exists C>0, \exists r>0, \forall X,Y,T\in\RZ,$
$$
g_X(Y-X)\le r^2\Longrightarrow C^{-1}{g_Y(T)}\le {g_X(T)}\le C {g_Y(T)}.
$$
\item {(2)} The symplectic dual metric 
 $g^\sigma$ is defined as
 $
 g^\sigma_X(T)=\sup_{g_X(U)=1}[T,U]^2.
$
The parameter of $g$ is defined as
$
\lambda_g(X)=\inf_{T\not=0 }\bigl( g^\sigma_X(T)/ g_X(T)\bigr)^{1/2}
$
and we shall say that $g$ satisfies the uncertainty principle
if
$\inf_X\lambda_g(X)\ge 1$.
\item {(3)} The metric $g$ is said to be temperate when
$\exists C>0, \exists N\ge 0, \forall X,Y,T\in\RZ,$
$$
{g^\sigma_X(T)}\le C {g^\sigma_Y(T)}\bigl(
1+g_X^\sigma(X-Y)
\bigr)^N.
$$
When the three properties above are satisfied, we shall say that
$g$ is admissible.
The constants appearing in (1) and (3)
will be called the structure constants of the metric $g$.
\item{(4)}   The function $M$ is said to be $g$-slowly varying
if
$\exists C>0, \exists r>0, \forall X,Y\in\RZ,$
$$
g_X(Y-X)\le r^2\Longrightarrow C^{-1}\le \frac{M(X)}{M(Y)}\le C.
$$
\item{(5)} The function $M$ is said to be $g$-temperate if
$\exists C>0, \exists N\ge 0, \ \forall X,Y\in\RZ,$
$$
\frac{M(X)}{M(Y)}\le C
\bigl(
1+g_X^\sigma(X-Y)
\bigr)^N.
$$
When $M$ satisfies (4) and (5), we shall say that $M$ is a $g$-weight.
\enddefinition
\remark{Remark}
If $g$ is a slowly varying metric and $M$ is $g$-slowly varying, there exists
$M_*\in S(M,g)$
such that there exists $C>0$ depending only 
on the structure constants
of $g$ such that $\forall X\in \RZ,\
C^{-1}\le \frac{M_*(X)}{M(X)}\le C.
$
That remark is classical and its proof is sketched in the appendix A.2.
\endremark
\definition{Definition 1.3.2}
Let $g$ be a metric on $\RZ$ and $M$ be a positive function defined on $\RZ$.
The set $S(M,g)$ is defined as the set of functions
$a\in\moo(\RZ)$ such that, for all $l\in \N$,
$
\sup_X\nuorm{a^{(l)}(X)}_{g_X}M(X)^{-1}<\io,
$
where
$a^{(l)}$ is the $l$-th derivative. 
It means that
$\forall l\in \N, \exists C_{l}, \forall X\in\RZ,\forall T_{1},\dots,T_{l}\in\RZ$,
$$
\val{a^{(l)}(X) (T_{1},\dots,T_{l})}\le C_{l}M(X)\prod_{{1\le j\le l}}g_{X}(T_{j})^{1/2}.
$$
\enddefinition
\subhead 
1.4. Partitions of unity
\endsubhead
We refer the reader to the chapter 18 in \cite{H6}
for the basic properties of admissible  metrics as well as for
the following lemma.
\proclaim{Lemma 1.4.1}
{Let $g$ be an admissible metric on $\RZ$.
There exists a sequence  $(X_k)_{k \in \N}$
of points in the phase space $\RZ$ and  positive numbers 
$r_0,N_0,$
such that the following properties are satisfied.
We define
$U_k, U_{k}^{\ast}, U_{k}^{\ast \ast}$ 
as the  $g_{k}=g_{X_k}$ balls with center $X_{k}$  and
radius $r_0, 2r_0, 4r_0$.
There exist two families of non-negative smooth functions on 
$\R^{2n}$, $(\chi_{k})_{k\in \N}$,
$(\psi_{k})_{k \in \N}$ such that
$$
\sum_{k}\chi_{k}(X) = 1,\  
\supp{\chi_{k}} \subset U_{k},\hs
\psi_{k} \equiv 1\  {\text{ on}\ }U_{k}^{\ast},\
\supp{\psi_{k}} \subset  U_{k}^{\ast \ast}. 
$$
Moreover,
$\chi_{k}, \psi_k \in S(1, g_{k} )$
with semi-norms bounded independently of $k$.
The overlap of the balls $U_{k}^{\ast \ast}$ 
is bounded, i.e.
$
\bigcap_{k \in {\Cal N}} U_{k}^{\ast \ast} \not= \emptyset\Longrightarrow \#{\Cal N}\le N_0.
$
Moreover,
$g_{X}\sim g_{k}$ all over $U_{k}^{\ast \ast}$
(i.e. the ratios $g_X(T)/g_k(T)$ are bounded
above and below by a fixed constant, 
provided that $X \in U_{k}^{\ast \ast}$).
}
\endproclaim
The next lemma in proved in \cite{BC}(see also 
lemma 6.3 in \cite{L5}).
\proclaim{Lemma 1.4.2}
{Let $g$ be an admissible metric on $\RZ$ and 
$\sum_k\chi_k(x,\xi) =1$ be a partition of 
unity related to $g$ as in the previous lemma.
There exists a positive constant $C$ such that for all $u \in L^2(\R^n)$
$$
C^{-1}\carno{u}{L^2(\R^n)}\le\sum_{k}\carno{\chi_k^w u}{L^2(\R^n)}
\le C\carno{u}{L^2(\R^n)},
$$
where $a^w$ stands for the Weyl quantization
of the symbol $a$.
}
\endproclaim
The following lemma is proved in \cite{BL}.
\proclaim{Lemma 1.4.3}
{Let $g$ be an admissible metric on $\RZ$,
$m$ be a weight
for $g$,
$U_k$
and $g_k$
as in lemma 1.4.1.
Let $(a_k)$ be a sequence of bounded symbols in
$S\bigl(m(X_k),g_k\bigr)$ such that,
for all non-negative
integers
$l, N$
$$
\sup_{k\in \N, T\in \RZ}\val{
m(X_k)^{-1}a_k^{(l)}(X)T^l
\bigl(
1+g_k^\sigma(X-U_k)
\bigr)^{N}g_k(T)^{-l/2}
}<+\io.
$$
Then the symbol
$a=\sum_k a_k$ makes sense and belongs to 
$S(m,g)$.
The important point
here is that no support condition is required
for the
$a_k$, but instead some decay estimates with respect
to 
$g^\sigma$.
The sequence
$(a_k)$
will be called a confined sequence in $S(m,g)$.
}
\endproclaim
\head
2. The geometry of \cps
\endhead
In this section and also in section 3, we shall consider that the phase space  is
equipped  with a {\it symplectic
quadratic form $\Gamma$}
($\Gamma$ is a positive definite quadratic form such that
$\Gamma =\Gamma^\sigma$, see the definition 1.3.1(2) above).
It is possible to find some linear symplectic coordinates $(x,\xi)$ in $\RZ$
such that
$$
\Gamma(x,\xi)=\val{(x,\xi)}^2=\sum_{1\le j\le n}
x_j^2+\xi_j^2.
$$
The running point of our Euclidean symplectic
$\RZ$
will be usually
denoted by $X$
or by an upper-case letter such as $Y,Z$.
The open
$\Gamma$-ball
with center $X$ and radius
$r$
will be denoted by
$B(X,r)$.
\subhead
2.1. The basic structure\endsubhead
Let
$q(t,X, \Lambda)$ be a smooth real-valued
function
defined
on $\Xi=\R\times\RZ\times[1,+\io)$,
vanishing for 
$\val t\ge 1$ and
satisfying
$$\gather
\forall k\in \N,
\sup_{\Xi}\norm{\p_X^k q}_\Gamma
\Lambda^{-1+\frac{k}{2}}
= \gamma_k<+\io,\text{ i.e. $q(t,\cdot)\in
S(\Lambda,
\Lambda^{-1}\Gamma)$},
\tag 2.1.1
\\
s>t\quad \text{and $q(t,X,\Lambda)>0$}\Longrightarrow
q(s,X,\Lambda)\ge 0.
\tag 2.1.2
\endgather$$
\remark{Notation 
}{
In this section and in the next section, the Euclidean norm $\Gamma(X)^{1/2}$
is fixed and the norms of the vectors and of the multilinear forms
are taken with respect to that norm. We shall write everywhere
$\val{\cdot}$
instead of $\norm{\cdot}_\Gamma$.
Furthermore, we shall say that
$C$ is a ``fixed" constant
if it depends only on a finite number of $\gamma_{k}$ above and on the dimension $n$.}
\endremark
We shall always omit
the dependence of $q$ with respect to 
the large parameter
$\Lambda$
and  write
$q(t,X)$ instead of
$q(t,X,\Lambda)$.
The operator
$Q(t)=q(t)^w$
will stand for the operator with Weyl symbol
$q(t,X)$.
We introduce now
for $t\in \R$, following \cite{H9},
$$\gather
\Bbb X_+(t)
=
\cup_{s\le t}\{X\in  \RZ, q(s,X)>0\},
\
\Bbb X_-(t)
=\cup_{s\ge t}
\{X\in  \RZ, q(s,X)<0
\},
\tag 2.1.3
\\
\Bbb X_0(t)
=\Bbb X_-(t)^c\cap\Bbb X_+(t)^c,
\tag 2.1.4
\endgather
$$
Thanks to (2.1.2),
$\Bbb X_+(t),\Bbb X_-(t)$
are disjoint open subsets of $\RZ$; moreover
$\Bbb X_0(t),
\Bbb X_0(t)\cup\Bbb X_\pm(t)
$ are closed
since their complements
are open.
The three sets $\Bbb X_0(t),\Bbb X_\pm(t)$
are two by two disjoint
with union $\RZ$
(note also that
$\overline{\Bbb X_\pm(t)}\subset
{\Bbb X_0(t)\cup \Bbb X_\pm(t)}$
since
${\Bbb X_0(t)\cup \Bbb X_\pm(t)}$ are closed).
When $t$ increases,
$\Bbb X_+(t)$ increases and 
$\Bbb X_-(t)$ decreases. 
\proclaim{Lemma 2.1.1}
Let $(E,d)$ be a metric space, $A\subset E$ and $\kappa >0$ be given.
We define $\Psi_{A,\kappa}(x)=\kappa$ if $A=\emptyset$
and if $A\not=\emptyset$, we define
$
\Psi_{A,\kappa}(x)=
\min\bigl(d(x,A), \kappa\bigr).$
The function
$\Psi_{A,\kappa}$
is valued in $[0,\kappa]$,
Lipschitz continuous
with  a Lipschitz constant $\le 1$.
Moreover,
the following implication holds:
$
A_1\subset A_2\subset E\Longrightarrow
\Psi_{A_1,\kappa}\ge
\Psi_{A_2,\kappa}.
$
\endproclaim
\demo{Proof}
The Lipschitz continuity assertion is obvious since $x\mapsto d(x,A)$
is Lipschitz continuous with Lipschitz constant 1.
The monotonicity property is trivially inherited from the distance function.
\qed
\enddemo
\proclaim{Lemma 2.1.2}
For each $X\in \RZ$,
the function $t\mapsto\Psi_{\Bbb X_+(t),\kappa}(X) $
is decreasing and for each $t\in \R$, the function
$X\mapsto\Psi_{\Bbb X_+(t),\kappa}(X)$ is supported in 
$\Bbb X_+(t)^c=\Bbb X_-(t)\cup \Bbb X_0(t)$.
For each $X\in \RZ$,
the function $t\mapsto\Psi_{\Bbb X_-(t),\kappa}(X) $
is increasing and 
for each $t\in \R$, the function
$X\mapsto\Psi_{\Bbb X_-(t),\kappa}(X)$ is supported in 
$\Bbb X_-(t)^c=\Bbb X_+(t)\cup \Bbb X_0(t)$.
As a consequence the function
$X\mapsto
\Psi_{\Bbb X_+(t),\kappa}(X)
\Psi_{\Bbb X_-(t),\kappa}(X)$
is supported in $\Bbb X_0(t)$.
\endproclaim \demo{Proof}
The monotonicity in $t$
follows from the fact that
$\Bbb X_+(t)$(resp.
$\Bbb X_-(t)$)
is increasing (resp. decreasing)
with respect to $t$ and from Lemma 2.1.1.
Moreover, if $X$ belongs to the open set $\Bbb X_\pm(t)$, one has 
$\Psi_{\Bbb X_\pm(t),\kappa}(X)=0$, implying the support property.
\qed
\enddemo
\proclaim{Lemma 2.1.3} For $\kappa >0, t\in\R, X\in\RZ$, we define
\footnote{
When the distances of $X$ to both
$\Bbb X_\pm(t)$ are less than $\kappa$, we have
$\sigma(t,X,\kappa)=\val{X-\Bbb X_-(t)}-
\val{X-\Bbb X_+(t)}.$
}
$$
\sigma(t,X,\kappa)=\Psi_{\Bbb X_-(t),\kappa}(X)-
\Psi_{\Bbb X_+(t),\kappa}(X).
\tag 2.1.5$$
The function $t\mapsto \sigma(t,X,\kappa)$ is increasing and valued in $[-\kappa, \kappa]$,
the function $X\mapsto \sigma(t,X,\kappa)$ is Lipschitz continuous with Lipschitz constant less than $2$;
we have
$$
\sigma(t,X,\kappa) =
\cases
\min(\val{X-\Bbb X_-(t)},\kappa)&\text{if $X\in\Bbb X_+(t)$,}
\\
-\min(\val{X-\Bbb X_+(t)},\kappa)&\text{if $X\in\Bbb X_-(t)$.}
\endcases
$$
We have
$\{X\in \RZ, \sigma(t,X,\kappa)=0\}\subset\Bbb X_0(t)\subset
\{X\in \RZ, q(t,X)=0\},$
and
$$\multline
\{X\in\RZ, \pm q(t,X)> 0\} \subset
\Bbb X_\pm(t)\subset
\{X\in \RZ, \pm\sigma(t,X,\kappa)> 0\}
\\
\subset
\{X\in \RZ, \pm\sigma(t,X,\kappa)\ge 0\}
\subset
\{X\in\RZ, \pm q(t,X)\ge 0\} .
\endmultline
\tag 2.1.6$$
\endproclaim \demo{Proof}
Everything follows from the previous lemmas, except for the first, fourth and sixth inclusions.
Note that
if $X\in \Bbb X_+(t)$, $\sigma(t,X,\kappa)=\min(\val{X-\Bbb X_-(t)},\kappa)$
 is positive (otherwise it vanishes and  $X\in \Bbb X_+(t)\cap\overline{\Bbb X_-(t)} \subset
\Bbb X_+(t)\cap\bigl(\Bbb X_-(t)\cup\Bbb X_0(t)\bigr)=\emptyset$). As a consequence,
we get the penultimate  inclusions
$\Bbb X_+(t)\subset\{X\in\RZ, \sigma(t,X,\kappa)>0\}$
and similarly
$\Bbb X_-(t)\subset\{X\in\RZ, \sigma(t,X,\kappa)<0\},
$
so that
$$
\{X\in \RZ, \sigma(t,X,\kappa)=0\}\subset \Bbb X_+(t)^c\cap\Bbb X_-(t)^c=\Bbb X_0(t),$$
giving the first inclusion.
The last inclusion
follows from the already established
$$
\{X\in\RZ,  q(t,X))<0\} \subset \Bbb X_-(t)\subset
\{X\in\RZ,  \sigma(t,X,\kappa)<0\}.\qed$$
\enddemo
\definition
{Definition 2.1.4}
Let $q(t,X)$ be as above.
We define
$$
\delta_0(t,X)=\sigma(t,X,\Lambda^{1/2})
\tag 2.1.7$$
and we notice that from the previous lemmas,
$
t\mapsto\delta_0(t,X)\
\text{is increasing}$,
valued in
$[-\Lambda^{1/2},\Lambda^{1/2}]$,
satisfying
$$
\val{\delta_0(t,X)-\delta_0(t,Y)}\le 2\val{X-Y}
\tag 2.1.8$$
and such that
$$\gather
\{X\in\RZ,  \delta_0(t,X)=0\}\subset\{X\in\RZ,  q(t,X)=0\},
\tag 2.1.9\\
\{X\in\RZ,  \pm q(t,X)>0\}\subset
\{X, \pm \delta_0(t,X)>0\}
\subset\{X,  \pm q(t,X)\ge 0\}.\tag 2.1.10
\endgather$$
\enddefinition
\proclaim{Lemma 2.1.5}
Let $f$ be a symbol in $S(\Lambda^m,\Lambda^{-1}\Gamma)$
where $m$ is a positive real number.
We define
$$
\lambda(X)=1+\max_{0\le j<2m\atop j\in\N }\bigl(
\nuorm{f^{(j)}(X)}_\Gamma^{\frac{2}{2m-j}}\bigr).
\tag 2.1.11$$
Then $f\in S(\lambda^m, \lambda^{-1}\Gamma)$
and the mapping from $S(\Lambda^m, \Lambda^{-1}\Gamma)$ to $S(\lambda^m, \lambda^{-1}\Gamma)$
is continuous. Moreover, with $
\gamma=\max_{0\le j<2m\atop j\in\N }\gamma_j^{\frac{2}{2m-j}},
$
where the $\gamma_j$ are the semi-norms of
$f$, we have
for all $X\in\RZ$,
$$
1\le \lambda (X)\le 1+\gamma \Lambda.
\tag 2.1.12$$
The metric $\lambda^{-1}\Gamma$ is admissible(def.1.3.1),
with structure constants depending only on $\gamma$. It
will be called the $m$-proper metric of $f$.
The function $\lambda$ above is a weight for the metric
$\lambda^{-1}\Gamma$
and will be called the  $m$-proper weight of 
$f$.
\endproclaim
\demo{The proof of this lemma is given in the appendix A.3}
\enddemo
\proclaim{Lemma 2.1.6}
Let $q(t,X)$ and $\delta_0(t,X)$
be as above. We define, with $\langle s\rangle=(1+s^2)^{1/2}$,
$$
\mu(t,X)=\langle\delta_0(t,X)\rangle^2+\val{\Lambda^{1/2}q'_X(t,X)}
+\val{\Lambda^{1/2}q''_{XX}(t,X)}^2.
\tag 2.1.13$$
The metric $\mu^{-1}(t,\cdot)\Gamma$
is slowly varying
 with structure constants depending only on a finite number of semi-norms of $q$ in $S(\Lambda,\Lambda^{-1}\Gamma)$. 
Moreover, there exists $C>0$, depending only on a finite
number of semi-norms of $q$, such that
$$
\mu(t,X)\le C \Lambda,\quad\frac{\mu(t,X)}{\mu(t,Y)}\le C(1+\val{X-Y}^2),
\tag 2.1.14$$
and we have
$$
\Lambda^{1/2}q(t,X)\in S(\mu(t, X)^{3/2},\mu^{-1}(t,\cdot)\Gamma),
\tag 2.1.15$$
so that the semi-norms depend
only the  semi-norms of $q$
in $S(\Lambda,\Lambda^{-1}\Gamma)$.
\endproclaim
\demo{Proof}
We notice first that
$$1+\max\bigl(
\val{\Lambda^{1/2}q'_X(t,X)},
\val{\Lambda^{1/2}q''_{XX}(t,X)}^2\bigr)
$$
is the
$1$-proper weight
of the vector-valued symbol
$\Lambda^{1/2}q'_X(t,\cdot)$.
Using the lemma A.2.2,
we get that
$\mu^{-1}\Gamma$ is slowly varying,
and the lemma A.2.1 provides the second part of (2.1.14).
From  the definition 2.1.4 and (2.1.1),
we obtain that
$\mu(t,X)\le C \Lambda+\langle\delta_0(t,X)\rangle^2
\le   C' \Lambda$ and
$
\Lambda^{1/2}q'_X(t,\cdot)\in
S(\mu(t, X),\mu^{-1}(t,\cdot)\Gamma).
$\par
We are left with the proof of
$
\val{\Lambda^{1/2}q(t,X)}\le C \mu^{3/2}(t,X).
$
Let us consider
$\widetilde{\mu}(t,X)$ the 3/2-proper weight
of $\Lambda^{1/2}q(t,X)$:
$$
\widetilde{\mu}(t,X)=1+\max_{j=0,1,2}\val{\Lambda^{1/2}q^{(j)}(t,X)}^{\frac{2}{3-j}},
$$
where all  the derivatives are taken with respect to  $X$;
if the maximum is realized for $j\in\{1,2\}$,
we get from Lemma 2.1.5 and (2.1.15)
that
$$
\val{\Lambda^{1/2}q(t,X)}\le 
\widetilde{\mu}(t,X)^{3/2}=(1+
\max_{j=1,2}\val{\Lambda^{1/2}q^{(j)}(t,X)}^{\frac{2}{3-j}}
)^{\frac{3}{2}}
\le
(1+
\max_{j=1,2}{\mu^{(\frac{3}{2}-\frac{j}{2})
(\frac{2}{3-j})
}}
)^{\frac{3}{2}}
\le 2\mu(t,X)^{3/2},
$$
which is the result that we had to prove.
We have eventually to deal with the case where the maximum in the definition of
$\widetilde \mu$ is realized for $j=0$;
note that if $\widetilde \mu(t,X)\le C_0$,
we obtain
$$
\val{\Lambda^{1/2}q(t,X)}\le 
\widetilde{\mu}(t,X)^{3/2}\le C_0^{3/2}\le 
C_0^{3/2}{\mu}(t,X)^{3/2},
$$
so we may also assume
$\widetilde \mu(t,X)>C_0.$
If $C_0>1$, we have
$C_0<\widetilde{\mu}(t,X)=1+(\Lambda^{1/2}\val{q(t,X)})^{\frac{2}{3}}$
entailing
$$(1-C_0^{-1})\widetilde\mu(t,X)\le \val{\Lambda^{1/2}{q(t,X)}}^{\frac{2}{3}}\le \widetilde\mu(t,X).
$$
Now if $h\in\RZ$ is such that
$\val h\le r  \widetilde\mu(t,X)^{1/2}$,
we get from the slow variation of the metric
$\widetilde\mu^{-1}\Gamma$,
that the ratio
$ \widetilde\mu(t,X+h)/ \widetilde\mu(t,X)$
is bounded above and below,
provided $r$ is small enough.
Using now that $\Lambda^{1/2}q(t,\cdot)\in S(\widetilde\mu^{3/2}(t,\cdot),\widetilde\mu^{-1}(t,\cdot)\Gamma)$,
we get by Taylor's formula
$$
\Lambda^{1/2}{q(t,X+h)}
=
\Lambda^{1/2}{q(t,X)}+\Lambda^{1/2}{q'(t,X)h}+\frac{1}{2}\Lambda^{1/2}{q''(t,X)h^2}
+O(\gamma_3\val h^3/6),
$$
so that
$$\align
\Lambda^{1/2}\val{q(t,X+h)}&\ge 
\Lambda^{1/2}\val{q(t,X)}-\widetilde\mu(t,X) \val{h}-\frac{1}{2}\val h^2\widetilde\mu(t,X)^{1/2}
-\gamma_3\val h^3/6
\\
&\ge
\Lambda^{1/2}\val{q(t,X)}-\widetilde\mu(t,X)^{3/2}
\underbrace{\Bigl(r+\frac{r^2}{2}+\gamma_3\frac{r^3}{6}\Bigr)}_{=\epsilon(r)}.
\endalign
$$
This  gives
$\Lambda^{1/2}\val{q(t,X+h)}
\ge 
\Lambda^{1/2}\val{q(t,X)}-\epsilon(r)\widetilde\mu(t,X)^{3/2},\
\lim_{r\rightarrow 0}\epsilon(r)=0,
$
so that, for $r, C_0^{-1}$ small enough,
$$
\val{\Lambda^{1/2}{q(t,X+h)}}
\ge
\bigl((1-C_0^{-1})^{3/2}-\epsilon(r)\bigr)\widetilde\mu(t,X)^{3/2}\ge \frac{1}{2}
\widetilde\mu(t,X)^{3/2}.
$$
As a consequence,
 the $\Gamma$-ball $B\bigl(X, r\widetilde\mu(t,X)^{1/2}\bigr)$
 is included in
 $\Bbb X_+(t)$ or in $ \Bbb X_-(t)$
 and thus,
 in the first case (the second case is similar)
 $
 \val{X-\Bbb X_+(t)}=0,\
 \val{X-\Bbb X_-(t)}\ge r \widetilde\mu(t,X)^{1/2},
 $
 $\bigl($otherwise
 $\val{X-\Bbb X_-(t)}<r \widetilde\mu(t,X)^{1/2}$
 and
 $\emptyset\not=B\bigl(X, r\widetilde\mu(t,X)^{1/2}\bigr)\cap \Bbb X_-(t)
 \subset
 \Bbb X_+(t)\cap\Bbb X_-(t)=\emptyset
 \bigr)$,
 implying that, with a fixed $r_0>0$,
 $$\delta_0(t,X)\ge\min(\Lambda^{1/2}, r \widetilde\mu(t,X)^{1/2})\ge r_0
 \widetilde\mu(t,X)^{1/2}\ge r_0\val{\Lambda^{1/2}{q(t,X)}}^{1/3},
 $$
 so that, in both cases,
 $
 \val{\Lambda^{1/2}{q(t,X)}}\le r_0^{-3}\val{\delta_0(t,X)}^3\le r_0^{-3}
 \mu(t,X)^{3/2},
$
qed.
\qed
\enddemo
\proclaim{Lemma 2.1.7}
Let $q(t,X), \delta_0(t,X), \mu(t,X)$
be as above. We define, 
$$
\nu(t,X)=\langle\delta_0(t,X)\rangle^2+\val{\Lambda^{1/2}q'_X(t,X)\mu(t,X)^{-1/2}}^2.
\tag 2.1.16$$
The metric $\nu^{-1}(t,\cdot)\Gamma$
is slowly varying
 with structure constants depending only on a finite number of semi-norms of $q$ in $S(\Lambda,\Lambda^{-1}\Gamma)$. 
There exists $C>0$, depending only on a finite
number of semi-norms of $q$, such that
$$
\nu(t,X)\le 2\mu(t,X)\le C \Lambda,\quad\frac{\nu(t,X)}{\nu(t,Y)}\le C(1+\val{X-Y}^2),
\tag 2.1.17$$
and we have
$$
\Lambda^{1/2}q(t,X)\in S(\mu(t, X)^{1/2}\nu(t,X), \nu(t,\cdot)^{-1}\Gamma),
\tag 2.1.18$$
so that the semi-norms of this symbol  depend
only the  semi-norms of $q$
in $S(\Lambda,\Lambda^{-1}\Gamma)$.
Moreover the function $\mu(t,X)$
is a weight for the metric $\nu(t,\cdot)^{-1}\Gamma$.
\endproclaim
\demo{Proof}
Let us check the two first  inequalities in (2.1.17).
From
$\val{\Lambda^{1/2}q'}\le \mu(t,X)\le C\Lambda$,
established in the previous lemma, we get
$$
\nu(t,X)\le \langle\delta_0(t,X)\rangle^2+\mu(t,X)\le 2\mu(t,X)\le 2C\Lambda.
$$
We introduce now the weight $\mu_*(t,X)$ as in (1.3.3)
so that the ratios
$\mu_*(t,X)/\mu(t,X)$ are bounded above and below
by some constants depending only on a finite number
of semi-norms
of $q$.
That weight
$\mu_*(t,X)$ belongs to $S(\mu,\mu^{-1}\Gamma)
=S(\mu_*,\mu_*^{-1}\Gamma)$.
We notice first that
$$\multline
\val{\Lambda^{1/2}(q\mu_*^{-1/2})'}^2\le 
2\val{\Lambda^{1/2}q'\mu_*^{-1/2}}^2+C_1
\val{\Lambda^{1/2}q\mu^{-1}}^2
\le
C_2\val{\Lambda^{1/2}q'\mu^{-1/2}}^2+C_1
\val{\Lambda^{1/2}q\mu^{-1/2}}
\overbrace{\val{\Lambda^{1/2}q\mu^{-3/2}}}^{\lesssim 1}
\\
\le
C_2\val{\Lambda^{1/2}q'\mu^{-1/2}}^2+C_3
\val{\Lambda^{1/2}q\mu^{-1/2}}.
\endmultline$$
Since we have also
\footnote{Below, the inequality $a\lesssim b$ means that $a\le Cb$
where $C$ is a constant depending only on a finite number
of semi-norms of $q$. The equivalence $a\sim b$ stands for $a\lesssim b$ and 
$b\lesssim a$.}
$$\multline
\val{\Lambda^{1/2}q'\mu^{-1/2}}\sim \val{\Lambda^{1/2}q'\mu_{*}^{-1/2}}\lesssim
\val{\Lambda^{1/2}(q\mu_*^{-1/2})'}
+
\val{\Lambda^{1/2}q\mu_*^{-1}}
\\\lesssim
\val{\Lambda^{1/2}(q\mu_*^{-1/2})'}
+
\underbrace{\val{\Lambda^{1/2}q\mu_*^{-3/2}}^{1/2}}_{{\lesssim 1}}
\val{\Lambda^{1/2}q\mu_*^{-1/2}}^{1/2}
\endmultline
$$
we get that
$$\widetilde\nu(t,X)=1+\max\bigl(
\val{\Lambda^{1/2}q'_X(t,X)\mu(t,X)^{-1/2}}^2,
\val{\Lambda^{1/2}q(t,X)\mu(t,X)^{-1/2}}\bigr)
\tag 2.1.19$$
is equivalent to the 1-proper weight 
of the  symbol
$\Lambda^{1/2}q(t,X)\mu_*(t,X)^{-1/2}$ in $S(\mu,\mu^{-1}\Gamma)$.
As a consequence, from
the lemma A.2.2,
we get that
$(\widetilde{\nu}+\langle\delta_0\rangle^2)^{-1}\Gamma$ is slowly varying.\par\no
{\it  \underbar{We need only to prove that}}
$$
\val{\Lambda^{1/2}q(t,X))\mu(t,X)^{-1/2}}\le C
{\nu}(t,X).
\tag 2.1.20$$
In fact, from (2.1.20), we shall obtain
$
\nu(t,X)\le \widetilde{\nu}(t,X)+\langle\delta_0(t,X)\rangle^2\le
(C+1)\nu(t,X) 
$
so that the metrics
$(\widetilde{\nu}+\langle\delta_0\rangle^2)^{-1}\Gamma$
and
${\nu}^{-1}\Gamma$
are equivalent and thus both slowly varying
(that property will also give the last inequality in (2.1.17) from Lemma A.2.1).
Moreover,
from Lemma 2.1.5, we have
$\Lambda^{1/2}q(t,X)\mu_*(t,X)^{-1/2}\in S(\widetilde{\nu},
\widetilde{\nu}^{-1}\Gamma),$ so that 
$$
\Lambda^{1/2}(q\mu_*^{-1/2})^{(k)}
\lesssim
\cases 
\nu^{1-k/2}&\text{for $k\le 2$, since 
$\Lambda^{1/2}q\mu_*^{-1/2}\in S(\widetilde{\nu};,\widetilde{\nu}^{-1}\Gamma)$
and $\widetilde{\nu}\lesssim {\nu},$}
\\
\mu^{1-k/2}\lesssim 
 {\nu}^{1-k/2}&\text{for $k\ge 2$, since $\Lambda^{1/2}q\mu_*^{-1/2}\in S(\mu; \mu^{-1}\Gamma)$
 and $\nu\lesssim \mu$,}
\endcases
$$
which implies that
$
\Lambda^{1/2}q\mu_*^{-1/2}\in S(\nu,\nu^{-1}\Gamma);
$
moreover, we have
$\mu_*^{1/2}\in S(\mu_*^{1/2},\nu^{-1}\Gamma)$
since, using
$\nu\lesssim \mu$, we get
$$
\val{(\mu_*^{1/2})^{(k)}}\lesssim
\mu^{\frac{1-k}{2}}\lesssim\mu^{\frac{1}{2}}\nu^{-k/2},
$$
entailing
$
\Lambda^{1/2}q\in S(\mu^{1/2}\nu,\nu^{-1}\Gamma),
$
i.e. (2.1.18).
On the other hand,
$\mu$ is slowly varying for $\nu^{-1}\Gamma$,
since
$$\val{X-Y}\ll \nu(t,X)^{1/2}(\lesssim \mu(t,X)^{1/2})
\quad\text{implies
$\val{X-Y}\ll \mu(t,X)^{1/2}$
}$$
and thus $\mu(t,X)\sim \mu(t,Y)$,
which proves along with (2.1.14) that
$\mu$ is a weight for 
$\nu^{-1}\Gamma$.\par\no
{\underbar{\it Let us now check (2.1.20).}}
This inequality is obvious if
$\val{\Lambda^{1/2}q\mu^{-1/2}}\le$
$\val{\Lambda^{1/2}q'\mu^{-1/2}}^2.$
Note that if $\widetilde\nu(t,X)\le C_0$,
we obtain
$
\val{\Lambda^{1/2}q\mu^{-1/2}}\le C_0
\le C_0\nu$
so we may  also assume
$\widetilde \nu(t,X)>C_0.$
If $C_0>1$, we have
$C_0<\widetilde{\nu}(t,X)=1+(\Lambda^{1/2}\val{q}\mu^{-1/2})$
entailing
$$(1-C_0^{-1})\widetilde\nu(t,X)\le \val{\Lambda^{1/2}q\mu^{-1/2}}\le \widetilde\nu(t,X).
$$
Now if $h\in\RZ$ is such that
$\val h\le r  \widetilde\nu(t,X)^{1/2}$,
we get from the slow variation of the metric
$\widetilde\nu^{-1}\Gamma$,
that the ratio
$ \widetilde\nu(t,X+h)/ \widetilde\nu(t,X)$
is bounded above and below,
provided $r$ is small enough.
Using now that $\Lambda^{1/2}q\mu_*^{-1/2}\in S(\widetilde\nu,\widetilde\nu^{-1}\Gamma)$,
we get by Taylor's formula
$$
\Lambda^{1/2}{q(t,X+h)\mu_*^{-1/2}(t,X+h)}
=
\Lambda^{1/2}q(t,X)\mu_*^{-1/2}(t,X)+\epsilon(r)\widetilde\nu(t,X),\quad\lim_{r\rightarrow 0}\epsilon(r)=0,
$$
so that, for $r, C_0^{-1}$ small enough,
$$
\val{\Lambda^{1/2}{q(t,X+h)\mu_*^{-1/2}(t,X+h)}}
\ge
\bigl((1-C_0^{-1})-\epsilon(r)\bigr)\widetilde\nu(t,X)\ge \frac{1}{2}
\widetilde\nu(t,X).
$$
As a consequence,
 the $\Gamma$-ball $B(X, r\widetilde\nu(t,X)^{1/2})$
 is included in
 $\Bbb X_+(t)$ or in $ \Bbb X_-(t)$
 and thus,
 in the first case (the second case is similar)
 $
 \val{X-\Bbb X_+(t)}=0,\
 \val{X-\Bbb X_-(t)}\ge r \widetilde\nu(t,X)^{1/2},
 $
 implying that, with a fixed $r_0>0$,
 $$\delta_0(t,X)\ge\min(\Lambda^{1/2}, r \widetilde\nu(t,X)^{1/2})\ge r_0
 \widetilde\nu(t,X)^{1/2}\ge r_0\val{\Lambda^{1/2}{q(t,X)\mu(t,X)^{-1/2}}}^{1/2},
 $$
 so that, in both cases,
 $
 \val{\Lambda^{1/2}{q(t,X)\mu(t,X)^{-1/2}}}\le C_0 \val{\delta_0(t,X)}^2\le C_0
 \nu(t,X),
$
qed.
The proof of the lemma is complete.
\qed
\enddemo
We wish now to discuss the normal forms attached to the metric
$\nu^{-1}(t,\cdot)\Gamma$ for the symbol $q(t,\cdot)$. In the sequel of this section,
we consider that $t$ is fixed.
\definition{Definition 2.1.8}
Let $0<r_{1}\le 1/2$
be given.
With $\nu$ defined in (2.1.16),
we shall say that
\item{$(1)$} $Y$ is a nonnegative
(resp. nonpositive)
 point at level $t$ if
$\delta_0(t,Y)\ge r_{1} \nu(t,Y)^{1/2}$,
(resp. $\delta_0(t,Y)\le -r_{1} \nu(t,Y)^{1/2}$).
\item{$(2)$}$Y$ is a gradient point at level $t$
if
$
\val{\Lambda^{1/2}q_{Y}'(t,Y)\mu(t,Y)^{-1/2}}^2\ge \nu(t,Y)/4
\
\text{and }\
\delta_0(t,Y)^2<r_{1}^2 \nu(t,Y).
$
\item{(3)}
$Y$ is a negligible point in the remaining cases
$
\val{\Lambda^{1/2}q_{Y}'(t,Y)\mu(t,Y)^{-1/2}}^2< \nu(t,Y)/4\
\text{and }\
\delta_0(t,Y)^2<r_{1}^2 \nu(t,Y).
$
Note that this implies 
$
\nu(t,Y)\le 1+r_{1}^2 \nu(t,Y)+\nu(t,Y)/4\le 1+\nu(t,Y)/2
$
and thus
$\nu(t,Y)\le 2$.
\enddefinition 
Note that if $Y$ is a nonnegative point, from (2.1.8) we get, for $T\in \RZ$,
$\val T\le 1, 0\le r\le r_{1}/4$
$$
\delta_{0}\bigl(t,Y+ r\nu^{1/2}(t,Y)T\bigr)
\ge \delta_{0}(t,Y)-2r\nu^{1/2}(t,Y)\ge\frac{r_{1}}{2} \nu^{1/2}(t,Y)
$$
and from (2.1.10), this implies that 
$
q(t,X) \ge 0
$
on the ball
$B(Y,r\nu^{1/2}(t,Y))$.
Similarly  if $Y$ is a nonpositive point,
$
q(t,X) \le 0
$
on the ball
$B(Y,r\nu^{1/2}(t,Y))$.
Moreover if $Y$ is a gradient point, we have
$\val {\delta_{0}(t,Y)}<r_{1}\nu(t,Y)^{1/2}$
so that, if $Y\in \Bbb X_{+}(t)$, we have
$\min(\val{Y-\Bbb X_{-}(t)},\Lambda^{1/2})<r_{1}\nu(t,Y)^{1/2}$
and if  $r_{1}$ is small enough, since $\nu\lesssim \Lambda$, we get that
$
\val{Y-\Bbb X_{-}(t)}<r_{1}\nu(t,Y)^{1/2}
$
which implies that there exists $Z_{1}\in\Bbb X_{-}(t)$ such that
$\val{Y-Z_{1}}<r_{1}\nu(t,Y)^{1/2}$.
On the segment $[Y,Z_{1}]$, the Lipschitz continuous function
is such that
$\delta_{0}(t,Y)>0$
($Y\in \Bbb X_{+}(t)$ cf. Lemma 2.1.3)
and 
$\delta_{0}(t,Z_{1})<0$
($Z_{1}\in \Bbb X_{-}(t)$); 
as a result, there exists a point $Z$ (on that segment) such that
$\delta_{0}(t,Z)=0$ and thus  $q(t,Z)=0$.
Naturally the discussion
for a gradient point $Y$ in $\Bbb X_{-}(t)$,
is analogous.
If the gradient point $Y$ belongs to $\Bbb X_{0}(t)$,
we get right away $q(t,Y)=0$, also from the lemma 2.1.3.
The function
$$
f(T)=\Lambda^{1/2}q\Bigl(t,Y+ r_{1}\nu^{1/2}(t,Y)T\Bigr)\mu(t,Y)^{-1/2} \nu(t,Y)^{-1}
\tag 2.1.21$$
satisfies for $r_{1}$ small enough with respect to the semi-norms of $q$
and $c_{0}, C_{0}, C_{1}, C_{2}$
fixed positive constants, $\val T\le 1$, from (2.1.18),
$
\val{f(T)}\le\val{S-T}C_{0}r_{1}\le C_{1}r_{1}^2,\quad\val{f'(T)}\ge r_{1}c_{0}, \quad
\val{f''(T)}\le C_{2} r_{1}^2.
$
The standard analysis
(see our appendix A.6)
of the Beals-Fefferman metric \cite{BF}
shows that, on 
$B(Y,r_{1}\nu^{1/2}(t,Y))$
$$\gather
q(t,X)= \Lambda^{-1/2}\mu^{1/2}(t,Y)
\nu^{1/2}(t,Y)
e(t,X)\beta(t,X), \tag 2.1.22\\
1\le e\in S(1,\nu(t,Y)^{-1}\Gamma),\ \beta\in S(\nu(t,Y)^{1/2},\nu(t,Y)^{-1}\Gamma),
\tag 2.1.23\\
\beta(t,X)=\nu(t,Y)^{1/2}
(X_{1}+\alpha(t,X')),\alpha\in S(\nu(t,Y)^{1/2},\nu(t,Y)^{-1}\Gamma).
\tag 2.1.24
\endgather$$
\proclaim{Lemma 2.1.9}
Let $q(t,X)$ be a smooth function satisfying $(2.1.1$-$2)$
and 
let $t \in[-1,1]$ be given. The metric $g_{t}$ on $\RZ$
 is defined as $\nu{(t,X)^{-1}\Gamma}$
where $\nu$ is defined in $(2.1.16)$.
There exists $r_{0}>0$, depending only on a finite number of semi-norms of $q$ in $(2.1.1)$
such that, for any $r\in]0,r_{0}]$, there exists a sequence of points
$(X_{k})$
in $\RZ$,
and  sequences of functions $(\chi_{k}), (\psi_{k})$ satisfying the properties in the  lemma 1.4.1
such that there exists a partition of $\N$,
$$
\N=E_{+}\cup E_{-}\cup E_{0} \cup E_{00}
$$
so that, according to the definition 2.1.8,
$k\in E_{+}$ means that $X_{k}$ is a nonnegative point,
($k\in E_{-}$:$X_{k}$ nonpositive point;
$k\in E_{0}$:$X_{k}$  gradient point,
$k\in E_{00}$:$X_{k}$  negligible  point).
\endproclaim
\demo{Proof}
This lemma is an immediate consequence
of the definition 2.1.8, of lemma 1.4.1
and of lemma 2.1.7, 
asserting that the metric $g_{t}$ is admissible.
\qed
\enddemo
\subhead
2.2. Some lemmas on $ {C^3}$ functions
\endsubhead
We prove in this section a key result
on the second derivative $f''_{XX}$ of a real-valued smooth
function $f(t,X)$
such that
$\tau-i f(t,x,\xi)$ satisfies \cps. 
The following claim gives a good qualitative
version of what is needed for our estimates;
we shall not use this result, so the reader may skip the proof
and proceed directly to the more technical Lemma 2.2.2.
\proclaim{Claim 2.2.1}
Let $f_1,f_2$ be two real-valued twice differentiable functions defined on
an open set $\Omega$ of 
$\R^N$ and such that
$
f_1^{-1}(\R_+^*)\subset f_2^{-1}(\R_+)
$
(i.e. $f_1(x)>0\Longrightarrow f_2(x)\ge 0$).
If  for some $\omega\in\Omega$, the conditions 
$
f_1(\omega)=f_2(\omega)=0,\quad df_1(\omega)\not=0,df_2(\omega)=0
$
are satisfied,
we have
$
f_2''(\omega)\ge 0
$
(as a quadratic form).
\endproclaim
\demo{Proof}
Using the obvious invariance by change of coordinates
of the statement,
we may assume $f_1(x)\equiv x_1$ and $\omega=0$.
The assumption is then for $x=(x_1,x')\in\R\times \R^{N-1}$
in a neighborhood
of the origin
$$
f_2(0)=0, df_2(0)=0,\quad x_1>0\Longrightarrow f_2(x_1,x')\ge 0.
$$
Using the second-order Taylor-Young formula for
 $f_2$,
we get
$
f_2(x)= \frac{1}{2}\poscal{f''_2(0)x}{x}+\epsilon(x)\val{x}^2,\
\lim_{x\rightarrow 0}\epsilon(x)=0,
$
and thus for
$T=(T_1,T'),\val T=1$, $\rho\not=0$ small enough,
the implication
$
T_1>0\Longrightarrow
\poscal{f''_2(0)T}{T}+2\epsilon(\rho T)\ge 0.
$
Consequently we have
$
\{S, \poscal{f''_2(0)S}{S}\ge 0\}\supset\{S, S_1> 0 \}
$
and since the larger set
is closed and stable by the symmetry with respect to the origin,
we get
that it contains also
$\{S, S_1\le 0 \}$,
which is the result
$f''_2(0)\ge 0.$
\enddemo
\remark{Remark}
This claim has the following consequence: take three functions
$f_1, f_2, f_3$,
twice differentiable on $\Omega$, such that, for $1\le j\le k\le 3$,
$f_j(x)>0\Longrightarrow f_k(x)\ge 0$.
Assume that,
at some point $\omega$
we have
$
f_1(\omega)=f_2(\omega)=f_3(\omega)=0,\quad
df_1(\omega)\not=0,
df_3(\omega)\not=0,
df_2(\omega)=0.
$
Then one has $f_2''(\omega)=0$. 
The claim 2.2.1 gives 
$f_2''(\omega)\ge 0$
and it can be applied to the couple $(-f_3, -f_2)$ to get
$-f_2''(\omega)\ge 0$.
\endremark
\remark{Notation}
The open Euclidean ball of $\R^N$
with center 0 and radius $r$
will be denoted by
$B_r$.
For a $k$-multilinear symmetric form
$A$ on $\R^N$,
we shall note
$
\norm A=\max_{\val T=1}{\val{AT^k}}
$
which is easily seen to be equivalent to the norm
$
\max_{\val {T_1}=\dots=\val{T_k}=1}{\val{A(T_1,\dots,T_k)}}
$
since the symmetrized
$T_1\otimes\dots\otimes T_k$
can be written a sum of $k^{\text{th}}$ powers.\endremark
\proclaim{Lemma 2.2.2}
Let $R_{0}>0$ and $f_1, f_2$ be real-valued functions
defined in $\bar B_{R_0}$.
We assume that
$f_1$ is $C^2$,
$f_2$ is $C^3$ 
 and for
 $ x\in \bar  B_{R_0},
$
$$
f_1(x)>0\Longrightarrow f_2(x)\ge 0.
\tag 2.2.1$$
We define the non-negative numbers $\rho_1, \rho_2,$ by
$$\rho_1=\max\bigl(\val{f_1(0)}^{\frac{1}{2}},\val{f'_1(0)}\bigr),\quad
\rho_2=\max\bigl(\val{f_2(0)}^{\frac{1}{3}},\val{f'_2(0)}^{\frac{1}{2}},\val{f''_2(0)}\bigr),
\tag 2.2.2
$$
and we assume that, with a positive $C_0$,
$$0<\rho_1,\quad
\rho_2\le C_0 \rho_1\le R_0.
\tag 2.2.3
$$
We define the non-negative numbers $C_1, C_2,C_3,$ by
$$
C_1=1+C_0\norm{f''_1}_{L^\io(\bar B_{R_0})},\quad
C_2=4+\frac{1}{3}\norm{f'''_2}_{L^\io(\bar B_{R_0})},\quad
C_3=C_2+4\pi C_1.
\tag 2.2.4$$
Assume that
 for some $\kappa_2\in[0,1]$, with $\kappa_2C_1\le 1/4$,
$$\gather
 \rho_{1}=\val{f'_1(0)}>0,\tag 2.2.5
 \\
\max\bigl(\val{f_2(0)}^{1/3},\val{f'_2(0)}^{1/2}\bigr)\le \kappa_2 \val{f''_2(0)},\tag 2.2.6
\\
B(0,\kappa_2^2\rho_2)
\cap\{x\in \bar B_{R_0}, f_1(x)\ge  0 \}\not=\emptyset.\tag 2.2.7
\endgather$$
Then 
we have
$$
\val{f_2''(0)_-}\le C_3\kappa_2\rho_2,\tag 2.2.8
$$
where
$f_2''(0)_-$
stands for the negative part of the quadratic form
$f_2''(0)$.
Note that, whenever $(2.2.7)$ is violated,
we get
$B(0,\kappa_2^2\rho_2)\subset\{x\in \bar B_{R_0}, f_1(x)< 0 \}$
(note that $\kappa_2^2\rho_2\le \rho_2\le R_0$)
and thus
$$
\text{distance}\bigl(0,\{x\in \bar B_{R_0}, f_1(x)\ge  0 \}\bigr)\ge \kappa_2^2\rho_2.\tag 2.2.9
$$
\endproclaim
\demo{Proof}
We may assume that for $x=(x_1,x')\in\R\times \R^{N-1}$,
$
\rho_1=\val{f'_1(0)}=\frac{\p f_1}{\p x_1}(0,0)$,\quad
$\frac{\p f_1}{\p x'}(0,0)=0,
$
so that
$$
f_1(x)\ge f_1(0)+\rho_1 x_1-\frac{1}{2}\norm{f''_1}_\io\val x^2.
\tag 2.2.10$$
Moreover, from  (2.2.7), we know that there exists
$z\in B(0,\kappa_2^2\rho_2)$ such that $f_1(z)\ge 0$.
As a consequence, we have
$
0\le f_1(z)\le f_1(0)+\rho_1z_1
+
\frac{1}{2}\norm{f''_1}_\io\kappa_2^4\rho_2^2
$
and thus
$$
f_1(x)\ge \rho_1 x_1-\rho_1\kappa_2^2 \rho_2-\frac{1}{2}\norm{f''_1}_\io(\val x^2+\kappa_2^4\rho_2^2).
\tag 2.2.11$$
On the other hand, we have
$$
 f_2(x)\le 
f_2(0)+f'_2(0)x+\frac{1}{2}f''_2(0)x^2+\frac{1}{6}\norm{f'''_2}_\io\val x^3
\le \kappa_2^3\rho_2^3
+\kappa_2^2\rho_2^2\val x+
\frac{1}{6}\norm{f'''_2}_\io\val x^3+
\frac{1}{2}f''_2(0)x^2
$$
and the implications, for $\val x\le R_0$,
$$\multline
 \rho_1 x_1>\rho_1\kappa_2^2\rho_2+\frac{1}{2}\norm{f''_1}_\io(\val x^2+\kappa_2^4\rho_2^2)
 \Longrightarrow f_1(x)>0\Longrightarrow f_2(x)\ge 0\Longrightarrow
\\ -\frac{1}{2}f''_2(0)x^2\le 
 \kappa_2^3\rho_2^3
+\kappa_2^2\rho_2^2\val x+
\frac{1}{6}\norm{f'''_2}_\io\val x^3.
\endmultline
\tag 2.2.12
$$
Let us take $x= \kappa_2\rho_2 y$ with $\val y=1$
(note that $\val x=\kappa_2 \rho_2\le R_0 $); 
the property (2.2.12) gives, using
$\rho_2/\rho_1\le C_0$,
$$
y_1>\kappa_2(1
+\norm{f''_1}_\io C_0)
 \Longrightarrow
  -f''_2(0)y^2\le\kappa_2\rho_2 
   (4
+
\frac{1}{3}\norm{f'''_2}_\io ),
$$
so that
$
\{ y\in\Bbb S^{N-1},
 -f''_2(0)y^2\le\kappa_2\rho_2 
   (4
+
\frac{1}{3}\norm{f'''_2}_\io )\}\supset
\{y\in\Bbb S^{N-1},y_1>\kappa_2(1
+\norm{f''_1}_\io C_0)\}
$
and since the larger set is closed
and stable by symmetry with respect to the origin,
we get, with 
$$C_1=1
+\norm{f''_1}_\io C_0, C_2=4
+
\frac{1}{3}\norm{f'''_2}_\io,
$$
the implication
$$
 y\in\Bbb S^{N-1}, \val{y_1}\ge \kappa_2C_1\Longrightarrow
  -f''_2(0)y^2\le\kappa_2\rho_2 
   C_2.
\tag 2.2.13$$
Let us now take
$y\in\Bbb S^{N-1}$,  such that $\val{y_1}< \kappa_2C_1(\le 1/4)$.
We may assume
$y=y_1\vec{e_1}\oplus y_2 \vec{e_2},$
with 
$\vec{e_1},\vec{e_2},$ orthogonal unit vectors
and $y_2=(1-y_1^2)^{1/2}$.
We consider the following rotation in the $(\vec{e_1},\vec{e_2})$ plane
with $\epsilon_0=\kappa_2C_1\le 1/4,$
$$
R=\pmatrix
\format \r&\quad\r
\\
\cos(2\pi \epsilon_0)&\sin(2\pi \epsilon_0)
\\
-\sin(2\pi \epsilon_0)&\cos(2\pi \epsilon_0)
\endpmatrix,
\quad\text{so that}\quad
\val{(Ry)_1}=\val{y_1\cos(2\pi \epsilon_0)+
y_2\sin(2\pi \epsilon_0)},
$$
and since $\epsilon_0\le 1/4$,
$$
\val{(Ry)_1}\ge-\val{y_1}+(1-y_1^2)^{1/2}4 \epsilon_0\ge
\epsilon_0(\sqrt{15} -1)>\epsilon_0=\kappa_2C_1.
$$
Moreover the rotation $R$ satisfies
$
\norm{R-\Id}\le 2\pi \epsilon_0=2\pi \kappa_2C_1.
$
We have, using (2.2.13) and
$\val{(Ry)_1}\ge \kappa_2 C_1$,  $\val y=1$,
$$\multline
-f''_2(0)y^2=
-f''_2(0)(Ry)^2-\poscal{f_2''(0)(y-Ry)}{y+Ry}
\le -f''_2(0)(Ry)^2+\val{f_2''(0)}\val{y-Ry}\val{y+Ry}
\\
\le
\kappa_2\rho_2 
   C_2+2\rho_2\val{y-Ry}\le 
   \kappa_2\rho_2 
   C_2+2\rho_22\pi \kappa_2C_1.
\endmultline$$
Eventually, for all $y\in\Bbb S^{N-1}$,
we have
$$
-f''_2(0)y^2\le 
  \kappa_2\rho_2(C_2+4\pi C_1)=C_3\kappa_2\rho_2.
\tag 2.2.14$$
Considering now the quadratic form
$Q=f''_2(0)$
and its canonical decomposition $Q=Q_+-Q_-$, 
we have, for all $y\in \R^N$,
$
\poscal{Q_-y}{y}\le \kappa_2\rho_2C_3\val y^2+\poscal{Q_+y}{y}.
$
Using now the canonical orthogonal projections $E_\pm$
on the positive
(resp. negative) eigenspaces, we write
$y=E_+y\oplus E_-y$ and 
we get that
$$
\poscal{Q_-y}{y}=\poscal{Q_-E_-y}{E_-y}
\le C_3\kappa_2\rho_2\val {E_-y}^2+\poscal{Q_+E_-y}{E_-y}
=C_3\kappa_2\rho_2\val {E_-y}^2\le
C_3\kappa_2\rho_2\val {y}^2,
$$
yielding (2.2.8).
The proof of Lemma 2.2.2 is complete.\qed
\enddemo
\proclaim{Lemma 2.2.3}
Let $f_1, f_2, f_3$ be real-valued functions
defined in $\bar B_{R_0}$.
We assume that
$f_1, f_3$ are $C^2$,
$f_2$ is $C^3$ 
 and for
 $ x\in \bar  B_{R_0}, 1\le j\le k\le 3,
$
$$
f_j(x)>0\Longrightarrow f_k(x)\ge 0.
\tag 2.2.15$$
We define the non-negative numbers $\rho_1, \rho_2,\rho_3$ by
$$\left.\matrix
\rho_1=\max\bigl(\val{f_1(0)}^{\frac{1}{2}},\val{f'_1(0)}\bigr)
\\
\rho_3=\max\bigl(\val{f_3(0)}^{\frac{1}{2}},\val{f'_3(0)}\bigr)
\endmatrix\right.
\quad
\rho_2=\max\bigl(\val{f_2(0)}^{\frac{1}{3}},\val{f'_2(0)}^{\frac{1}{2}},\val{f''_2(0)}\bigr),
\tag 2.2.16
$$
and we assume that, with a positive $C_0$,
$$0<\rho_1, \rho_3\quad\text{and}\quad
\rho_2\le C_0\min( \rho_1,\rho_3)
\le C_0\max( \rho_1,\rho_3)\le R_0.
\tag 2.2.17
$$
We define the non-negative numbers $C_1, C_2,C_3,$ by
$$\aligned
C_1&=1+C_0\max(\norm{f''_1}_{L^\io(\bar B_{R_0})},
\norm{f''_3}_{L^\io(\bar B_{R_0})}),
\\
C_2&=4+\frac{1}{3}\norm{f'''_2}_{L^\io(\bar B_{R_0})},\quad
C_3=C_2+4\pi C_1.
\endaligned
\tag 2.2.18
$$
Assume that
 for some $\kappa_1,\kappa_3\in[0,1]$, and $0<\kappa_2C_3\le 1/2$,
$$\gather
 \val{f_1(0)}^{1/2}\le \kappa_1\val{f'_1(0)},\quad
  \val{f_3(0)}^{1/2}\le \kappa_3\val{f'_3(0)},\tag 2.2.19
 \\
B(0,\kappa_2^2\rho_2)
\cap\{x\in \bar B_{R_0}, f_1(x)\ge  0 \}\not=\emptyset,\tag 2.2.20
\\
B(0,\kappa_2^2\rho_2)
\cap\{x\in \bar B_{R_0}, f_3(x)\le  0 \}\not=\emptyset.\tag 2.2.21
\endgather$$
Then 
we have
$$
\max\bigl(\val{f_2(0)}^{1/3},\val{f'_2(0)}^{1/2}\bigr)\le 
\rho_2\le \kappa_2^{-1}\max\bigl(\val{f_2(0)}^{1/3},\val{f'_2(0)}^{1/2}\bigr)
.\tag 2.2.22
$$
Note that, whenever $(2.2.20)$ or 
$(2.2.21)$ is violated,
we get
$$\text{
$B(0,\kappa_2^2\rho_2)\subset\{x\in \bar B_{R_0}, f_1(x)< 0 \}$
or
$B(0,\kappa_2^2\rho_2)\subset\{x\in \bar B_{R_0}, f_3(x)>0 \}$
}$$
and thus
$$
\text{dist}\bigl(0,\{x\in \bar B_{R_0}, f_1(x)\ge  0 \}\bigr)\ge \kappa_2^2\rho_2
\text{\ or\  }
\text{dist}\bigl(0,\{x\in \bar B_{R_0}, f_3(x)\le  0 \}\bigr)\ge \kappa_2^2\rho_2.
\tag 2.2.23$$
\endproclaim
\demo{Proof}  
This follows almost immediately from the previous lemma and it is analogous to the remark
following
the claim 2.2.1:
assuming that we have
$$
\max\bigl(\val{f_2(0)}^{1/3},\val{f'_2(0)}^{1/2}\bigr)\le  \kappa_2\val{f''_2(0)}
\tag 2.2.24$$
will yield
$
\val{f''_2(0)}\le C_3\kappa_2\rho_2
$
by applying lemma 2.2.2 (note that $\kappa_2 C_1\le \kappa_2 \frac{C_3}{4\pi}\le \frac{1}{8\pi}<1/4)$
to the couples
$(f_1, f_2)$ and $(-f_3, -f_2)$;
consequently, if (2.2.24) is satisfied, we get
$$
\max\bigl(\val{f_2(0)}^{1/3},\val{f'_2(0)}^{1/2}\bigr)\le
\rho_2\le
\max\bigl(\val{f_2(0)}^{1/3},\val{f'_2(0)}^{1/2}, C_3\kappa_2 \rho_2\bigr)
$$
and since $C_3\kappa_2<1$,
it yields
$$
\max\bigl(\val{f_2(0)}^{1/3},\val{f'_2(0)}^{1/2}\bigr)=
\rho_2,
\tag 2.2.25$$
which implies (2.2.22).
Let us now suppose that (2.2.24) does not hold,
and that we have 
$\kappa_2\val{f''_2(0)}
<
\max\bigl(\val{f_2(0)}^{1/3},\val{f'_2(0)}^{1/2}\bigr).
$
This implies (2.2.22):
$$
\max\bigl(\val{f_2(0)}^{1/3},\val{f'_2(0)}^{1/2}\bigr)\le 
\rho_2\le \kappa_2^{-1}\max\bigl(\val{f_2(0)}^{1/3},\val{f'_2(0)}^{1/2}\bigr).
$$
The proof of the lemma is complete.\qed
\enddemo
\remark{Remark}
We shall apply this lemma to a ``fixed" $\kappa_2$, depending only on the constant $C_3$ such as
$\kappa_2=1/(2C_3)$.
\endremark
\subhead
2.3. Inequalities for symbols\endsubhead
In this section, we apply the results of the previous section to obtain various inequalities
on symbols
linked to our symbol $q$ introduced in (2.1.1). Our main result is the following theorem.
\proclaim{Theorem 2.3.1}
Let $q$ be a symbol satisfying $(2.1.1$-$2)$ and $\delta_0, \mu,\nu$
as defined
above in $(2.1.7)$, $(2.1.13)$ and $(2.1.16)$.
For the real numbers $t',t,t''$,  and $X\in \RZ$,
we define
$$\align
N(t',t'',X)&=\frac{\langle\delta_0(t',X)\rangle}{\nu(t',X)^{1/2}}
+\frac{\langle\delta_0(t'',X)\rangle}{\nu(t'',X)^{1/2}},\tag 2.3.1
\\
R(t,X)&=\Lambda^{-1/2}\mu(t,X)^{1/2} \nu(t,X)^{-1/2}\langle\delta_0(t,X).\tag 2.3.2
\endalign$$
Then there exists a constant $C_0\ge 1$,
depending only
on a finite number of semi-norms
of $q$ in $(2.1.1)$, such that, for $t'\le t\le t''$,
we have
$$
C_0^{-1}R(t,X)\le N(t',t'',X)+\frac{\delta_0(t'',X)-\delta_0(t,X)}{ \nu(t'',X)^{1/2}}+
\frac{\delta_0(t,X)-\delta_0(t',X)}{ \nu(t',X)^{1/2}}.
\tag 2.3.3
$$
\endproclaim
\demo{Proof}
We are given $X\in \RZ$ and $t'\le t\le t''$ real numbers.
\par\no
{\it \underbar{First reductions.}}
First of all, we may assume 
that, for some positive (small) $\kappa$ to be chosen later, we have
$$
\langle\delta_0(t',X)\rangle\le \kappa\nu(t',X)^{1/2}\quad\text{and}\quad
\langle\delta_0(t'',X)\rangle\le \kappa\nu(t'',X)^{1/2}.
\tag 2.3.4$$
In fact, otherwise, we have
$
N(t',t'',X)> \kappa
$
and since from (2.1.14), we have $\mu(t,X)\le C \Lambda$
where $C$ depends only
on a finite number of semi-norms
of $q$, we get from (2.3.2), (2.1.16)
$$
R(t,X)\le C^{1/2}
\nu(t,X)^{-1/2}\langle\delta_0(t,X)\rangle\le C^{1/2}
\le C^{1/2}\kappa^{-1}N(t',t'',X),
$$
so that we shall only need
$$
\boxed{C_0\ge C^{1/2}\kappa^{-1}}
\tag 2.3.5$$
to obtain (2.3.3).
Also, we may assume
that, with the same positive (small) $\kappa$,
$$
\nu(t,X)\le \kappa^2 \nu(t',X)
\quad\text{and}\quad
\nu(t,X)\le \kappa^2 \nu(t'',X).
\tag 2.3.6$$
Otherwise, we would have
for instance
$\nu(t,X)>\kappa^2 \nu(t',X)$
and since $t\ge t'$,
$$\multline
R(t,X)\le \Lambda^{-1/2}\mu(t,X)^{1/2}\kappa^{-1}\frac{\langle\delta_0(t,X)\rangle}{\nu(t',X)^{1/2}}
\\
\le C^{1/2}\kappa^{-1}
\biggl(\frac{\langle\delta_0(t',X)\rangle+\val{\overbrace{\delta_0(t,X)-\delta_0(t',X)}^{\ge 0}}}{\nu(t',X)^{1/2}}\biggr)
\\
\le 
C^{1/2}\kappa^{-1} N(t',t'',X)+C^{1/2}\kappa^{-1} 
\frac{{\delta_0(t,X)-\delta_0(t',X)}}{\nu(t',X)^{1/2}},
\endmultline$$
which implies (2.3.3) provided that (2.3.5) holds.
Finally, we may also  assume that
$$
\nu(t,X)\le \kappa^2\mu(t,X),
\tag 2.3.7$$
otherwise we would have,
using that
$\delta_0(t',X)\le \delta_0(t,X)\le \delta_0(t'',X)\ $
and the convexity of $s\mapsto\sqrt{1+s^2}=\langle s\rangle$,
$$
R(t,X)
\le\kappa^{-1}\frac{\langle\delta_0(t,X)\rangle}{\Lambda^{1/2}}
\le
\kappa^{-1}\frac{\langle\delta_0(t',X)\rangle}{\Lambda^{1/2}}
+\kappa^{-1}\frac{\langle\delta_0(t'',X)\rangle}{\Lambda^{1/2}}
$$
and this implies, using $\nu(t',X),\nu(t'',X)\le C \Lambda$ (see (2.1.17)),
$$
R(t,X)\le C^{1/2}\kappa^{-1}\frac{\langle\delta_0(t',X)\rangle}{\nu(t',X)^{1/2}}
+C^{1/2}\kappa^{-1}\frac{\langle\delta_0(t'',X)\rangle}{\nu(t'',X)^{1/2}},
$$
which gives (2.3.3) provided that (2.3.5) holds.
On the other hand,
we may assume that
$$
\max\bigl(\langle\delta_0(t,X)\rangle,\kappa^{1/2}\val{\Lambda^{1/2}q'(t,X)}^{1/2}\bigr)
\le 2\kappa\mu(t,X)^{1/2}.
\tag 2.3.8$$
Otherwise, we would have either
$$
\mu(t,X)^{1/2}\le \frac{1}{2}\kappa^{-1}\langle\delta_0(t,X)\rangle\le
\frac{1}{2}\kappa^{-1}\nu(t,X)^{1/2}\underbrace{\le}_{\text{from (2.3.7)}} \frac{1}{2}\mu(t,X)^{1/2}
$$
which is impossible,
or we would have
$$\multline
\mu(t,X)^{1/2}\le \frac{1}{2}\kappa^{-1/2}\val{\Lambda^{1/2}q'(t,X)}^{1/2}
\overbrace{\le}^{\text{from (2.1.16) }}
\frac{1}{2}\kappa^{-1/2}
\nu(t,X)^{1/4}\mu(t,X)^{1/4}
\\
\underbrace{\le}_{\text{from (2.3.7)}} 
\frac{1}{2}\mu(t,X)^{1/2},
\quad\text{(which is also impossible).}
\endmultline$$
The estimate (2.3.8) implies that, for $\boxed{\kappa<1/16},$
$$\multline
\Lambda\val{q''(t,X)}^2
\underbrace{\le}_{(2.1.13)} 
 \mu(t,X)
\underbrace{\le}_{(2.1.13)} 
\langle\delta_0(t,X)\rangle^2
+
\val{\Lambda^{1/2}q'(t,X)}
+
\Lambda\val{q''(t,X)}^2
\\
\underbrace{\le}_{(2.3.8)} 
(4\kappa^2+4\kappa)\mu(t,X)
+
\Lambda\val{q''(t,X)}^2,
\endmultline$$
and thus
$$
\Lambda\val{q''(t,X)}^2\le \mu(t,X)\le\frac{1}{1-8\kappa}\Lambda\val{q''(t,X)}^2\le 2
\Lambda\val{q''(t,X)}^2.
\tag 2.3.9$$
This implies that
$$
\multline
R(t,X)\le \Lambda^{-1/2}2^{1/2 }
\Lambda^{1/2}\val{q''(t,X)}\frac{\langle\delta_0(t,X)\rangle}{
\bigl(
\langle\delta_0(t,X)\rangle^2+\Lambda\val{q'(t,X)}^2
\mu(t,X)^{-1}
\bigr)
^{1/2}}
\\
\le
2^{1/2 }\val{q''(t,X)}.
\endmultline
\tag 2.3.10$$
{\it \underbar{Rescaling the symbols}}. We sum-up our situation,
changing the notations so that $X=0,t'=t_1, t=t_2, t''=t_3,
\nu_1=\nu(t',0),
\nu_2=\nu(t,0),
\nu_3=\nu(t'',0),
\delta_j=\delta_0(t_j,0), \mu_j=\mu(t_j,0)$.
The following conditions are satisfied:
$$
\left.\aligned\langle\delta_1\rangle\le\kappa\nu_1^{1/2},\quad&\text{}\quad
\langle\delta_3\rangle\le \kappa\nu_3^{1/2},
\\
\nu_2\le \kappa^2 \nu_1,
\quad&\text{}\quad
\nu_2\le \kappa^2 \nu_3,\quad \nu_2\le \kappa^2\mu_2
\\
R(t_2,0)\le
2\val{q''(t_2,0)}
\frac{\langle\delta_2\rangle}{\langle\delta_2\rangle+\val{q'(t_2,0)}/\val{q''(t_2,0)}}
&\le 2\val{q''(t_2,0)},
\\
\Lambda\val{q''(t_2,0)}^2\le \mu_2&\le 2\Lambda\val{q''(t_2,0)}^2,
\\
\kappa<1/16,\quad &\quad C_0\ge \kappa^{-1}C^{1/2},
\endaligned\right\}
\tag 2.3.11
$$
where $\kappa>0$ is to be chosen later 
and
$C$ depends only
on a finite number of semi-norms
of $q$.
We define now the smooth functions
$f_1,f_2$ defined on $\RZ$ by
$$
f_1(Y)=q(t_1,Y)\Lambda^{1/2}\mu_1^{-1/2}
,\quad
f_2(Y)=\nu_1^{1/2}q(t_2,Y),
\tag 2.3.12$$
and we note (see (2.1.1)-(2.1.15)) that $\norm{f_1''}_{L^\io}$
and
$\norm{f_2'''}_{L^\io}$ are bounded above by semi-norms of $q$;
moreover the assumption (2.2.1) holds for that couple of functions, from (2.1.2).
\proclaim{Lemma 2.3.2}
We define
$$
\mu_{12}^{1/2}=\max\bigl(\langle\delta_2\rangle,\val{\nu_1^{1/2}q'(t_2,0)}^{1/2},
\val{\nu_1^{1/2}q''(t_2,0)}\bigr).
\tag 2.3.13$$
If
$
\max\bigl(\langle\delta_2\rangle,\kappa^{1/2}\val{\nu_1^{1/2}q'(t_2,0)}^{1/2}\bigr)
> 2\kappa\mu_{12}^{1/2},
$
then $(2.3.3)$ is satisfied provided $C_0\ge 3/\kappa$.
\endproclaim
\demo{Proof}
We  have either
$
\val{\nu_1^{1/2}q''(t_2,0)}\le \mu_{12}^{1/2}\le \frac{1}{2}\kappa^{-1}\langle\delta_2\rangle
$
implying
$$
\val{q''(t_2,0)}\le\frac{1}{2\kappa}\frac{\langle\delta_2\rangle}{\nu_1^{1/2}}
\le
\frac{1}{2\kappa}\frac{\langle\delta_1\rangle}{\nu_1^{1/2}}+\frac{1}{2\kappa}
\frac{\delta_2-\delta_1}{\nu_1^{1/2}}
$$
which gives  (2.3.3)
(using $R(t_2,0)\le 2\val{q''(t_2,0)}$ in (2.3.11)), provided
$
{C_0\ge 1/\kappa},
$
or we  have
$$
\val{\nu_1^{1/2}q''(t_2,0)}\le \mu_{12}^{1/2}< \frac{1}{2}\kappa^{-1/2}
\val{\nu_1^{1/2}q'(t_2,0)}^{1/2},
$$
implying
$\dis
\frac{\val{q''(t_2,0)}^2}{\val{q'(t_2,0)}}\le\frac{1}{4\kappa\nu_1^{1/2}}
$ so that
(using  $R(t_2,0)\le 2\val{q''(t_2,0)}^2\langle\delta_{2}\rangle/\val{q'(t_2,0)}$ in (2.3.11)),
we get
$
R(t_2,0)\le
\frac{1}{2\kappa}\frac{\langle\delta_2\rangle}{\nu_1^{1/2}},
$
which gives similarly (2.3.3), provided
$
{C_0\ge 1/(2\kappa)}.
$
\qed
\enddemo
A consequence of this lemma is that we may assume
$$
\max\bigl(\langle\delta_2\rangle,\kappa^{1/2}\val{\nu_1^{1/2}q'(t_2,0)}^{1/2}\bigr)
\le 2\kappa\mu_{12}^{1/2}
= 2\kappa
\max\bigl(\langle\delta_2\rangle,\val{\nu_1^{1/2}q'(t_2,0)}^{1/2},
\val{\nu_1^{1/2}q''(t_2,0)}\bigr),
$$
and since $\kappa <1/4$, we get
$$
\mu_{12}^{1/2}=\val{\nu_1^{1/2}q''(t_2,0)},
\quad
\max\bigl(\langle\delta_2\rangle,\kappa^{1/2}\val{\nu_1^{1/2}q'(t_2,0)}^{1/2}\bigr)\le 
2\kappa
\val{\nu_1^{1/2}q''(t_2,0)}.
\tag 2.3.14$$
\proclaim{Lemma 2.3.3}
The functions $f_1, f_2$ defined in $(2.3.12)$ satisfy the assumptions $(2.2.1$-$2$-$3$-$4$-$5$-$6)$
in  the lemma
2.2.2.
\endproclaim
\demo{Proof}
We have already checked
(2.2.1).
We know from Lemma 2.1.7
that,
with a constant $C$ depending only on a finite number
of semi-norms of $q$ (see (2.1.18)),
$$
\val{f_1(0)=q(t_1,0)\Lambda^{1/2}\mu_1^{-1/2}}^{1/2}\le C \nu_1^{1/2},
$$
but we may assume here that $C\le 1/2$:
if we had 
$\val{f_1(0)}>\nu_1^{1/2}/2$,
the function $f_{1}$ would be positive (resp.negative) on $B(0, r_{0}\nu_{1}^{1/2})$,
with some fixed $r_{0}>0$
and consequently we would have
$\val {\delta_{1}}\ge r_{0}\nu_1^{1/2}$. But we know that
$\langle\delta_{1}\rangle\le \kappa\nu_{1}^{1/2}$,
so we can choose a priori $\kappa$ small enough
so that $\val {\delta_{1}}\ge r_{0}\nu_1^{1/2}$
does not occur.
From (2.3.11), we have 
$\langle\delta_1\rangle\le \kappa \nu_1^{1/2}$, the latter implying $f'_1(0)\not=0$ from (2.1.16)
since $\kappa^2 <3/4$
and more precisely
$$
\rho_1=\val{f'_1(0)}\ge (1-\kappa^2)^{1/2}\nu_1^{1/2}\ge \nu_1^{1/2}/2.
\tag 2.3.15$$
Moreover we have, from (2.1.18) and
$\nu_2\le \kappa^2 \nu_1$ in (2.3.11),
$$
\max(\val{\nu_1^{1/2}q'(t_2,0)}^{1/2},
\val{\nu_1^{1/2}q''(t_2,0)})\le \mu_{12}^{1/2}\le C_1 \nu_1^{1/2},
$$
with a constant $C_1$ depending only on a finite number
of semi-norms of $q$ and thus 
$$
\max(\val{f'_2(0)}^{1/2},\val{f''_2(0)})\le 2C_1 \rho_1.
\tag 2.3.16$$
Moreover, we have from Lemma 2.1.7,
$
\Lambda^{1/2}\val{q(t_2,0)}\mu_2^{-1/2}\le C_2 \nu_2,
$
so that with constants $C_2, C_3$
depending only on a finite number
of semi-norms of $q$, using (2.3.8), we get
$$
\val{f_2(0)}\le\nu_1^{1/2} C_2\nu_2\Lambda^{-1/2}\mu_2^{1/2}\le
\nu_1^{1/2} C_3 \nu_2\le C_3\kappa^2\nu_1^{3/2}.
$$
That  property and (2.3.16-15) give
(2.2.3) with $R_0=C \rho_1$,
where
$C$
depends only on a finite number
of semi-norms of $q$.
We have already seen that the constants occurring in (2.2.4)
are bounded
above by semi-norms of $q$ and that (2.2.5) holds.
Let us now check (2.2.6).
We already know that, from (2.3.14),
$$
\val{f'_2(0)}^{1/2}=\val{\nu_1^{1/2}q'(t_2,0)}^{1/2}
\le 2\kappa^{1/2}\val{\nu_1^{1/2}q''(t_2,0)}=2\kappa^{1/2}\val{f''_2(0)}.
\tag 2.3.17
$$
If we have
$
\val{\nu_1^{1/2}q(t_2,0)}\ge \kappa^{1/2}\mu_{12}^{3/2}
$
then for 
$\val h\le \kappa^{1/3}\mu_{12}^{1/2}$,
we get, using
$\nu_1\lesssim \Lambda$
and Taylor's formula along with (2.3.13-14),
$$\multline
\val{\nu_1^{1/2}q(t_2,h)}\ge 
\kappa^{1/2}\mu_{12}^{3/2}-4\kappa^{4/3}\mu_{12}^{3/2}-\frac{1}{2}\kappa^{2/3}\mu_{12}^{3/2}
-C\nu_1^{1/2}\Lambda^{-1/2}\kappa \mu_{12}^{3/2}
\\
=\mu_{12}^{3/2}(\kappa^{1/2}-4\kappa^{4/3}-\frac{\kappa^{2/3}}{2}-C'\kappa)
\ge \mu_{12}^{3/2}\kappa^{1/2}/2>0,
\endmultline$$
provided
$\kappa$ is small enough with respect
to a constant depending only on a finite number
of semi-norms of $q$;
that inequality
implies that the ball $B(0,\kappa^{1/3}\mu_{12}^{1/2})$ is 
included in $\Bbb X_+(t_2)$
or in
$\Bbb X_-(t_2)$
implying that
$\val{\delta_0(t_2,0)=\delta_2}\ge \kappa^{1/3}\mu_{12}^{1/2}$
which is incompatible
with (2.3.14), provided 
$\kappa<2^{-3/2}$,  since (2.3.14) implies
$\val {\delta_2}\le 2\kappa\mu_{12}^{1/2}$.
Eventually, we get
$$
\val{f_2(0)}^{1/3}=\val{\nu_1^{1/2}q(t_2,0)}^{1/3}\le \kappa^{1/6}\mu_{12}^{1/2}=
\kappa^{1/6} \val{f''_2(0)}
\tag 2.3.18$$
and with (2.3.18) we obtain (2.2.6) with
$$
\kappa_2=\kappa^{1/6}.
\tag 2.3.19$$
The proof of Lemma 2.3.3 is complete.
\qed
\enddemo\no
{\it \underbar{End of the proof of Theorem 2.3.1.}}
To apply Lemma 2.2.2, we have to suppose (2.2.7).
In that case we get
$
\nu_1^{1/2}\val{q''(t_2, 0)_-}=\val{f''_2(0)_-}\le C\kappa_2\rho_2= C\kappa^{1/6}\nu_1^{1/2}\val{q''(t_2, 0)}
$
i.e.
$$
\val{q''(t_2, 0)_-}\le C\kappa^{1/6}\val{q''(t_2, 0)}.
\tag 2.3.20$$
If (2.2.7) is not satisfied,
we obtain, according to (2.2.9),
(2.3.19)
and 
$\mu_{12}=\nu_1^{1/2}\val{q''(t_2,0)},$
$$
\delta_0(t_1,0)=\delta_1\le -\kappa^{1/3}\nu_1^{1/2}\val{q''(t_2, 0)},
$$
which gives
$
\frac{1}{3}R(t_2,0)\le \val{q''(t_2, 0)}\le \kappa^{-1/3}\frac{\val {\delta_1}}{\nu_1^{1/2}}
$
and (2.3.3) provided
$C_0\ge 3\kappa^{-1/3}$.
If we introduce now the smooth functions
$F_1,F_2$
defined on $\RZ$ by
$$
F_1(Y)= -q(t_3, Y)\Lambda^{1/2}\mu_3^{-1/2}
,\quad
F_2(Y)=-\nu_3^{1/2}q(t_2,Y),
\tag 2.3.21$$
starting over our discussion,
we see that (2.3.3) is satisfied,
provided 
$$\kappa\le \kappa_0\quad\text{and}\quad C_0\ge \gamma_0\kappa^{-1},
\tag 2.3.22$$
where
$\kappa_0,\gamma_0$
are positive constants depending only on the
semi-norms of $q$, except
in the case where we have (2.3.20) \underbar{and}
$$
\val{q''(t_2, 0)_+}\le C\kappa^{1/6}\val{q''(t_2, 0)}.
\tag 2.3.23$$ 
Naturally,
since
$\val{q''(t_2, 0)}=\val{q''(t_2, 0)_+}+\val{q''(t_2, 0)_-}$,
the estimates
(2.3.20-23) cannot be both true
for a  $\kappa$
small
enough
with respect to a constant depending on  a finite number of semi-norms of $q$
and a non-vanishing
$q''(t_2, 0)$
(that vanishing is prevented by the penultimate line in (2.3.11)).
The proof of Theorem 2.3.1 is complete.
\enddemo
\remark{Remark 2.3.4}
The reader may find our proof quite tedious, but referring him to the simpler remark following claim
2.2.1,
we hope that he can find there some motivation
to read the details of our argument, which is the rather natural quantitative statement
following from that remark.
On the other hand, Theorem 2.3.1 is analogous to one of the key argument provided by N.Dencker
in \cite{D3}
 in which he proves, using our notations in the theorem,
$$
R(t,X)\lesssim N(t',t'',X)+\delta_0(t'',X)-\delta_0(t',X)
\tag 2.3.24$$
which is weaker than our (2.3.3). In particular, $R$ (and $N$) looks like a symbol of order 0
(weight 1)
whereas the right-hand-side of  (2.3.24) contains the difference
$\delta_0(t'',X)-\delta_0(t',X)$, which looks like a symbol of order 1/2.
Our theorem gives a stronger and in some sense
more homogeneous version of N.Dencker's result, which will lead
to improvements in the remainder's estimates.
Also, we note the (inhomogeneous)
estimate
$$
\Lambda^{-1/2}\mu(t,X)^{1/2}\nu(t,X)^{-1/2}\lesssim N(t',t'',X),
$$
which is in fact a consequence of our proof, but is not enough
to handle the remainder's estimate below in our proof, and which will not be used:
in fact (2.3.3) implies
$$\multline
\Lambda^{-1/2}\mu^{1/2}\nu^{-1/2}=R\langle\delta_{0}\rangle^{-1}\\\lesssim
\frac{N(t',t'',X)}{\langle\delta_{0}(t,X)\rangle}+\frac{\delta_{0}(t'',X)-\delta_{0}(t,X)}{
\nu(t'',X)^{1/2}\langle\delta_{0}(t,X)\rangle}
+\frac{\delta_{0}(t,X)-\delta_{0}(t',X)}{
\nu(t',X)^{1/2}\langle\delta_{0}(t,X)\rangle}
\\\lesssim
\frac{N(t',t'',X)}{\langle\delta_{0}(t,X)\rangle}+\frac{1}{\nu(t'',X)^{1/2}}+\frac{1}{\nu(t',X)^{1/2}}
\lesssim N(t',t'',X).
\endmultline$$
\endremark
\subhead
2.4. Quasi-convexity
\endsubhead
A differentiable function $\psi$
of one variable
is said to be quasi-convex
on $\R$
if
$\dot \psi (t)$
does not change sign from $+$ to $-$ 
for increasing $t$
(see \cite{H7}). In particular, a differentiable
convex
function
is such that
$\dot \psi (t)$ is increasing
and is thus quasi-convex.
\definition{Definition 2.4.1 }
Let 
$\sigma_1:\R\rightarrow \R$
be an increasing function, $C_1>0$ and let 
$\rho_1:\R\rightarrow \R_+$.
We shall say that
$\rho_1$ is quasi-convex
with respect to $(C_1,\sigma_1)$
if for $t_1,t_2,t_3\in\R$,
$$
t_1\le t_2\le t_3\Longrightarrow
\rho_1(t_2)\le C_1\max\bigl(\rho_1(t_1),\rho_1(t_3)
\bigr)
+\sigma_1(t_3)-\sigma_1(t_1).
\tag 2.4.1$$
When $\sigma_1$
is a constant function and $C_1=1$,
this is the definition of quasi-convexity.
\enddefinition
\proclaim{Lemma 2.4.2}
Let 
$\sigma_1:\R\rightarrow \R$
be an increasing function and
let $\omega:\R\rightarrow \R_+$.
We define
$$
\rho_1(t)=\inf_{t'\le t\le t''}\Bigl(
\omega(t')+\omega(t'')+
\sigma_1(t'')-\sigma_1(t')\Bigr).
\tag 2.4.2$$
Then the function $\rho_1$
is quasi-convex with respect to $(2,\sigma_1)$.
\endproclaim
\demo{Proof}
We consider
$t_1\le t_2\le t_3$ three real numbers.
We have
$$\align
\rho_1(t_2)&=\inf_{t'\le t_2\le t''}
\Bigl(
\omega(t')+\omega(t'')+
\sigma_1(t'')-\sigma_1(t')
\Bigr)
\\
&\le
\inf_{t'\le t_1, t_3\le t''}\Bigl(
\omega(t')+\omega(t'')+
\sigma_1(t'')-\sigma_1(t_3)+\sigma_1(t_1)
-\sigma_1(t')\Bigr)
+\sigma_1(t_3)-\sigma_1(t_1)
\\
&\le
\inf_{t'\le t_1\le t''_1,\atop t'_3\le
t_3\le t''
}\Bigl(
\omega(t')+\omega(t''_1)+
\omega(t'_3)+
\omega(t'')+
\sigma_1(t'')-\sigma_1(t'_3)+\sigma_1(t''_1)
-\sigma_1(t')\Bigr)
+\sigma_1(t_3)-\sigma_1(t_1)
\\
&=\rho_1(t_1)+\rho_1(t_3)
+\sigma_1(t_3)-\sigma_1(t_1)
\le
2\max(\rho_1(t_1),\rho_1(t_3))
+\sigma_1(t_3)-\sigma_1(t_1).\qed
\endalign$$
\enddemo
The following lemma is  due to L.H\"ormander \cite{H9}.
\proclaim{Lemma 2.4.3}
Let 
$\sigma_1:\R\rightarrow \R$
be an increasing function and
let $\omega:\R\rightarrow \R_+$.
Let $T>0$ be given.
We consider the function
$\rho_1$
as defined in Lemma 2.4.2
and we define
$$
\Theta_T(t)=
\sup_{-T\le s\le t}
\left\{
\sigma_1(s)-\sigma_1(t)+\frac{1}{2T}
\int_s^t\rho_1(r) dr-\rho_1(s)
\right\}.
\tag 2.4.3$$
Then we have
$$
2T\p_t(\Theta_T+\sigma_1)\ge \rho_1,\quad
\text{and for $\val t\le T$,}
\quad
\val{\Theta_T(t)}\le \rho_1(t).
\tag 2.4.4$$
\endproclaim
\demo{Proof}
We have
$\Theta_T(t)\ge-\rho_1(t),$
and
$$
\Theta_T(t)+\sigma_1(t)=
\underbrace{
\sup_{-T\le s\le t}
\left\{
\sigma_1(s)+\frac{1}{2T}
\int_s^0\rho_1(r) dr-\rho_1(s)
\right\}}_{\text{increasing with $t$}}+\frac{1}{2T}\int_0^t\rho_1(r) dr,
$$
so that
$
\p_t(\Theta_T+\sigma_1)\ge
\frac{1}{2T} \rho_1.
$
Moreover, from the proof of
Lemma 2.4.2, we obtain
for
$s\le r\le t$ that
$
\rho_1(r)\le \rho_1(s)+\rho_1(t)+
\sigma_1(t)-\sigma_1(s)
$
and thus
$$\frac{1}{2T}\int_s^t \rho_1(r) dr\le
\frac{1}{t-s}\int_s^t \rho_1(r) dr
\le
\rho_1(s)+\rho_1(t)+
\sigma_1(t)-\sigma_1(s)
$$
which gives
$
\Theta_T(t)\le\rho_1(t),
$
ending the proof of the lemma.\qed
\enddemo
\definition{Definition 2.4.4} For $T>0, X\in \RZ, \val t\le T$, we define
$$
\omega(t,X)=\frac{\langle\delta_0(t,X)\rangle}{\nu(t,X)^{1/2}},\
\sigma_1(t,X)=\delta_0(t,X),\
\eta(t,X)=\int_{-T}^t \delta_{0}(s,X)\Lambda^{-1/2}ds+2T,
\tag 2.4.5
$$
where $\delta_0, \nu$ are defined in (2.1.7),(2.1.16).
For $T>0$, $(t,X)\in \R\times\RZ$,
we define
$\Theta(t,X)$
by the formula
(2.4.3)
$$
\Theta(t,X)=
\sup_{-T\le s\le t}
\left\{
\sigma_1(s,X)-\sigma_1(t,X)+\frac{1}{2T}
\int_s^t\rho_1(r,X) dr-\rho_1(s,X)
\right\},
\tag 2.4.6$$
where $\rho_{1}$ is defined by (2.4.2). We define also
$$
m(t,X)=\delta_{0}(t,X)+\Theta(t,X)+T^{-1}\delta_{0}(t,X)\eta(t,X).
\tag 2.4.7$$
\enddefinition
\proclaim{Theorem 2.4.5}
With the notations above for $\Theta,\rho_1, m$, 
with $R$ and $C_0$ defined in Theorem 2.3.1, we have for $T>0$, $\val t\le T$,
$X\in \RZ, \Lambda \ge 1$,
$$\gather
\val{\Theta(t,X)}\le \rho_1(t,X)\le 2\frac{\langle\delta_0(t,X)\rangle}{\nu(t,X)^{1/2}},
\quad
\val{\sigma_1(t,X)}= \val{\delta_{0}(t,X)},
\tag 2.4.8
\\
C_0^{-1}R(t,X)\le \rho_1(t,X)
\le
2T\frac{\p}{\p t}\Bigl(\Theta(t,X)+\sigma_1(t,X)\Bigr),
\tag 2.4.9
\\
0\le \eta(t,X)\le 4T,\quad
\frac{d}{dt}\bigl(\delta_0\eta\bigr)\ge
\delta_{0}^2
\Lambda^{-1/2},\quad\val{\eta'_{X}(t,X)}\le 4T\Lambda^{-1/2},
\tag 2.4.10
\\
T\frac{d}{dt} m\ge
\frac{1}{2}\rho_{1}+\delta_{0}^2\Lambda^{-1/2}\ge\frac{1}{2C_{0}} R+\delta_{0}^2\Lambda^{-1/2}\ge
\frac{1}{2^{3/2}C_{0}}
\langle\delta_{0}\rangle^2\Lambda^{-1/2}.
\tag 2.4.11
 \endgather$$
\endproclaim
\demo{Proof}
It follows immediately from the previous results:
the first estimate in (2.4.8) is (2.4.4),
whereas the second is due to 
$\rho_{1}\le 2\omega
$
which follows from (2.4.2).
The equality in (2.4.8) follows from  Definition 2.4.4.
The first inequality in (2.4.9)
is a consequence of (2.4.2) and (2.3.3) and the second is (2.4.4).
The first two inequalities in (2.4.10)
are a consequence of 
$
\val{\delta_{0}(t,X)}\le \Lambda^{1/2}
$
which follows from definition 2.1.4.
The third inequality
reads
$$
\frac{d}{dt}\bigl(\delta_0\eta\bigr)=\dot\delta_{0}\eta+\delta_{0}\dot\eta\ge
\delta_{0}\dot\eta=\delta_{0}^2\Lambda^{-1/2},
$$
and the fourth inequality in (2.4.10)
follows from
(2.1.8).
Let us check finally (2.4.11): 
since
$
m=\delta_{0}+\Theta+T^{-1}\delta_{0}\eta,
$
(2.4.4) and the already proven (2.4.10) imply
$
T\frac{d}{dt} m\ge\frac{1}{2}\rho_{1}+\delta_{0}^2\Lambda^{-1/2}
$
and (2.4.9)(proven) gives 
$$\multline
\frac{1}{2}\rho_{1}+\delta_{0}^2\Lambda^{-1/2}\ge\frac{1}{2C_{0}} R+\delta_{0}^2\Lambda^{-1/2}=
\frac{1}{2C_{0}} \Lambda^{-1/2}\mu^{1/2}\nu^{-1/2}\langle\delta_{0}\rangle
+\delta_{0}^2\Lambda^{-1/2}\\
\underbrace{\ge}_{\text{from (2.1.17)}}
\frac{1}{2C_{0}}\Lambda^{-1/2}\bigl(2^{-1/2}\langle\delta_{0}\rangle+\delta_{0}^2\bigr)
\ge \frac{1}{2^{3/2}C_{0}}\Lambda^{-1/2}\langle\delta_{0}\rangle^2,
\endmultline$$
completing the proof of Theorem 2.4.5.
\qed\enddemo
\head
3. Energy estimates
\endhead
\subhead
3.1. Preliminaries
\endsubhead
\definition{Definition 3.1.1}
Let $T>0$ be given. 
With
$m$ defined in (2.4.7),
we define for $\val t \le T$,
$$
M(t)= \w{m(t,X)},\tag 3.1.1
$$
where the Wick quantization is given by the definition A.1.1.
\enddefinition
\proclaim{Lemma 3.1.2}
With $T>0$ and $M$ given above, we have with $\rho_1$ given in $(2.4.2)$,
for $\val t\le T$, $\Lambda\ge 1,$
$$\multline
\frac{d}{dt}M(t)\ge \frac{1}{2T}\w{\rho_1(t,X)}+
\frac{1}{T}
\w{{(\delta_{0}^2 )}}
\Lambda^{-1/2}\ge 
\frac{1}{2C_{0}T} \w{R}+T^{-1}\w{{(\delta_{0}^2 )}}\Lambda^{-1/2}
\\
\ge 
\frac{1}{2^{3/2}C_{0}T}
\w{({\langle\delta_{0}\rangle^2})}\Lambda^{-1/2}.
\endmultline\tag 3.1.2$$
$$\gather
\val{\Theta(t,X)}\le \rho_1(t,X)\le 2\frac{\langle\delta_0(t,X)\rangle}{\nu(t,X)^{1/2}},
\tag 3.1.3
\\
T^{-1}\val{\delta_{0}(t,X)\eta(t,X)}\le 4\val{\delta_0(t,X)},\tag 3.1.4\\
T^{-1}\val{
\delta'_{0X}(t,X)\eta(t,X)}
+
T^{-1}\val{
\delta_{0X}(t,X)\eta'_{X}(t,X)}\le 12.\tag 3.1.5
\endgather$$
\endproclaim
\demo{Proof}
The derivative in (3.1.2)
is taken in the distribution sense, i.e. the first inequality in (3.1.2)
means
that (A.1.5) is satisfied with
$$
a(t,X)=m(t,X)-\frac{1}{2T}\int_{-T}^t
\rho_{1}(s,X)ds-\frac{1}{T}\Lambda^{-1/2}\int_{-T}^t\delta_{0}(s,X)^2ds.
$$
It follows in fact from (2.4.11).
The other inequalities in (3.1.2)
follow directly from
(2.4.11)
and the fact that the Wick quantization is positive
(see (A.1.3)).
The inequality
(3.1.3) is (2.4.8) and (3.1.4) follows from (2.4.10)
whereas (3.1.5) is a consequence of
(2.1.8), (2.4.10) and Definition 2.1.4.
\qed
\enddemo
\proclaim{Lemma 3.1.3}
Using the definitions above and the notation $(A.1.4)$,
we have
$$\align
\Theta(t,\cdot)\ast\exp-2\pi\Gamma&\in S\bigl(
{\langle\delta_0(t,\cdot)\rangle}{\nu(t,\cdot)^{-1/2}}
,\Gamma\bigr),\tag 3.1.6
\\
\delta_{0}(t,\cdot)\ast\exp-2\pi\Gamma&\in S\bigl(
{\langle\delta_0(t,\cdot)\rangle}
,\Gamma\bigr),\tag 3.1.7
\\
\delta'_{0X}(t,\cdot)\ast\exp-2\pi\Gamma&\in S\bigl(
1,\Gamma\bigr),\tag 3.1.8
\\
T^{-1}\eta(t,\cdot)\ast\exp-2\pi\Gamma&\in S\bigl(
1,\Gamma\bigr),\tag 3.1.9
\\
T^{-1}\eta(t,\cdot)'_{X}\ast\exp-2\pi\Gamma&\in S\bigl(
\Lambda^{-1/2},\Gamma\bigr),\tag 3.1.10
\endalign
$$
with semi-norms independent of $T\le 1$ and of $t$
for $\val t\le T$.
According to the definition 1.3.1, the function $X\mapsto \langle\delta_{0}(t,X)\rangle$ is a $\Gamma$-weight.
\endproclaim
\demo{Proof} The last statement follows from (2.1.8).
The inequalities  ensuring (3.1.6---10)
are then immediate consequences
of the lemmas 3.1.2 and A.1.3.
\qed
\enddemo
\subhead
3.2.  Stationary estimates for the model cases
\endsubhead
Let $T>0$ be given and
$Q(t)= q(t)^w$  given by (2.1.1-2). We define $M(t)$ according to (3.1.1).
We consider
$$
\re\bigl(Q(t)M(t)\bigr)
=\frac{1}{2}
Q(t)M(t)+\frac{1}{2}
M(t) Q(t)=P(t).
\tag 3.2.1$$
We have, omitting now the variable $t$ fixed throughout all this section 3.2, 
$$
P= \re{\Bigl[q^w\w{\bigl(\delta_{0}(1+T^{-1}\eta)\bigr)}
+q^w\w{\Theta
}\Bigr]}.
\tag 3.2.2$$
\subsubhead {{\bf [1].}}
Let us assume first that $q=\Lambda^{-1/2}
\mu^{1/2}\nu^{1/2}\beta e_{0}$ with $\beta\in S(\nu^{1/2}, \nu^{-1}\Gamma), 
1\le e_{0}\in S(1, \nu^{-1}\Gamma)$
and $\delta_{0}=\beta$. 
Moreover, we assume
$0\le T^{-1}\eta\le 4,  T^{-1}\val{\eta'}\le 4\Lambda^{-1/2}$, \quad
$\val{\Theta}\le C \langle\delta_{0}\rangle\nu^{-1/2}$.
Here $\Lambda, \mu, \nu$
are assumed to be positive constants such that
$\Lambda\ge \mu\ge \nu\ge 1$
\endsubsubhead
Then using the lemma A.1.5
with $$a_{1}=\beta e_{0},\quad m_{1}=\langle \beta\rangle,\quad a_{2}=(1+T^{-1}\eta)e_{0}^{-1},m_{2}=\nu^{-1/2},$$
we get, with obvious notations,
$$
\w{(\delta_{0}e_{0})}\w{\bigl(e_{0}^{-1}(1+T^{-1}\eta)\bigr)}=\w{\bigl(\delta_{0}(1+T^{-1}\eta)\bigr)}+{S(\langle \delta_{0}\rangle\nu^{-1/2},\Gamma)}^w
$$
and as a consequence from the proposition A.1.2(2), 
we obtain, with $$\beta_{0}=\beta e_{0},\quad \eta_{0}= e_{0}^{-1}(1+T^{-1}\eta),
\tag 3.2.3$$
the identity
$
\bigl({\beta_{0}}^w+{S(\nu^{-1/2},\nu^{-1}\Gamma)}^w\bigr)
\w{\eta_{0}}
=\w{\bigl(\delta_{0}(1+T^{-1}\eta)\bigr)}+{S(\langle \delta_{0}\rangle\nu^{-1/2},\Gamma)}^w,
$
entailing
$$
\w{\bigl(\delta_{0}(1+T^{-1}\eta)\bigr)}={\beta_{0}}^w\w{\eta_{0}}+{S(\langle \delta_{0}\rangle\nu^{-1/2},\Gamma)}^w.
$$
As a result, we have
$$\multline
QM=
\Lambda^{-1/2}
\mu^{1/2}\nu^{1/2}\beta_{0}^w\beta_{0}^w\w{\eta_{0}}
+\beta_{0}^w S(\overbrace{\Lambda^{-1/2}
\mu^{1/2}\nu^{1/2}\langle \delta_{0}\rangle\nu^{-1/2}}^{=
\Lambda^{-1/2}
\mu^{1/2}\langle \delta_{0}\rangle
},\Gamma)^w
\\+\beta_{0}^w
 S(\underbrace{\Lambda^{-1/2}
\mu^{1/2}
\nu^{1/2}\langle \delta_{0}\rangle \nu^{-1/2}}_{=\Lambda^{-1/2}
\mu^{1/2}\langle \delta_{0}\rangle},\Gamma)^w.
\endmultline
$$
This implies that, with $\gamma_{0}=1/\sup e_{0}>0,$ (so that $1\le e_{0}\le \gamma_{0}^{-1}$)
$$\align
2\re QM&=2
\Lambda^{-1/2}
\mu^{1/2}\nu^{1/2}\beta_{0}^w\w{\eta_{0}}\beta_{0}^w
+2\re 
\beta_{0}^w\Lambda^{-1/2}
\mu^{1/2}\nu^{1/2}\overbrace{\bigl[\beta_{0}^w,\w{\eta_{0}}\bigr]}^{\in
S({
\nu^{-1/2}},\Gamma)^w.
}
\\&\hfill\hskip195pt
+\re \beta_{0}^w
 S({\Lambda^{-1/2}
\mu^{1/2}
\langle \delta_{0}\rangle},\Gamma)^w
\\&=
2\Lambda^{-1/2}
\mu^{1/2}\nu^{1/2}\beta_{0}^w\w{\eta_{0}}\beta_{0}^w
+\re \beta_{0}^w
 S({\Lambda^{-1/2}
\mu^{1/2}
\langle \delta_{0}\rangle},\Gamma)^w
\\&
\underbrace{\ge}_{\text{
since $\eta_{0}\ge e_{0}^{-1}$ }
\atop
\text{
from $\eta \ge 0$ in (2.4.10)}
} 
2\Lambda^{-1/2}
\mu^{1/2}\nu^{1/2}\beta_{0}^w\gamma_{0}\beta_{0}^w
+\beta_{0}^w b_{0}^w+\bar b_{0}^w\beta_{0}^w, 
\tag 3.2.4
\endalign
$$
with 
$b_{0}\in  S(\Lambda^{-1/2}
\mu^{1/2}
\langle \delta_{0}\rangle,\Gamma).$
With the notation $\lambda=
\Lambda^{-1/2}
\mu^{1/2}\nu^{1/2}\gamma_{0}$,
we use  the identity,
$$
\Lambda^{1/2}
\mu^{1/2}\nu^{1/2}\beta_{0}^w\gamma_{0}\beta_{0}^w
+\beta_{0}^w b_{0}^w+\bar b_{0}^w\beta_{0}^w
=
\bigr(\lambda^{1/2}\beta_{0}^w+\lambda^{-1/2}\bar b_{0}^w\bigl)
\bigr(\lambda^{1/2}\beta_{0}^w+\lambda^{-1/2} b_{0}^w\bigl)
-\lambda^{-1}\bar b_{0}^wb_{0}^w,
$$
so that from (3.2.4), we obtain
with $b_{1}$ real valued in
$S(
\underbrace{\Lambda^{1/2}
\mu^{-1/2}\nu^{-1/2}\Lambda^{-1}\mu
\langle \delta_{0}\rangle^2}_{
\Lambda^{-1/2}
\mu^{1/2}\nu^{-1/2}
\langle \delta_{0}\rangle^2
},\Gamma)$, the inequality
$$
2\re QM+b_{1}^w\ge \Lambda^{-1/2}\mu^{1/2}\nu^{1/2}\gamma_{0}\beta_{0}^w\beta_{0}^w.
\tag 3.2.5$$
Using now
(A.1.11), 
we get, with a ``fixed" constant $C$, that
$$\multline
b_{1}^w
\le C
\Lambda^{-1/2}\mu^{1/2}\nu^{-1/2} \w{(1+\beta^2)}
=C \Lambda^{-1/2}\mu^{1/2}\nu^{-1/2} \Id
+C \Lambda^{-1/2}\mu^{1/2}\nu^{-1/2} \w{(\beta_{0}^2e_{0}^{-2})}
\\
\le C
\Lambda^{-1/2}\mu^{1/2}\nu^{-1/2} \Id
+C\Lambda^{-1/2}\mu^{1/2}\nu^{-1/2} \w{(\beta_{0}^2)},
\endmultline$$
and since, from the proposition A.1.2(2), we have
$$
\w{(\beta_{0}^2)}={(\beta_{0}^2)}^w+S(1,\nu^{-1}\Gamma)^w=\beta_{0}^w\beta_{0}^w
+S(1,\nu^{-1}\Gamma)^w,
$$
the inequality (3.2.4) implies
$$\multline
2\re QM+
 C
\Lambda^{-1/2}\mu^{1/2}\nu^{-1/2} \Id
+C\Lambda^{-1/2}\mu^{1/2}\nu^{-1/2} {\beta_{0}}^w{\beta_{0}}^w
+S(\Lambda^{-1/2}\mu^{1/2}\nu^{-1/2},\nu^{-1}\Gamma)^w
\\
\ge 2\re QM+
b_{1}^w
\ge
\Lambda^{-1/2}\mu^{1/2}\nu^{1/2}\gamma_{0}\beta_{0}^w\beta_{0}^w,
\endmultline$$
so that 
$$
\re QM + S(
\Lambda^{-1/2}\mu^{1/2}\nu^{-1/2} 
,\Gamma)^w\ge 
\beta_{0}^w\beta_{0}^w
(\Lambda^{-1/2}\mu^{1/2}\nu^{1/2}\gamma_{0}- C'
\Lambda^{-1/2}\mu^{1/2}\nu^{-1/2}).
\tag 3.2.6
$$
The rhs of (3.2.6) is nonnegative provided
$
\nu\ge C'\gamma_{0}^{-1}
$
and since
$C'\gamma_{0}^{-1}$ is a fixed constant, we may first suppose that this condition is satisfied; 
if it is not the case, we would have that
$\nu$ is bounded above by a fixed constant and since $\nu\ge 1$,
that would imply
$q\in S(\Lambda^{-1/2}\mu^{1/2},\Gamma)
$
and $P\in S(\Lambda^{-1/2}\mu^{1/2},\Gamma)^w$.
In both cases, we get
$$
\re QM + S(
\Lambda^{-1/2}\mu^{1/2}\nu^{-1/2} 
,\Gamma)^w\ge 0.\tag 3.2.7
$$
\subsubhead {{\bf [2].}}
Let us assume now that $q\ge 0$,
$q\in S(\Lambda^{-1/2}\mu^{1/2}\nu,\nu^{-1}\Gamma)$,
$\gamma_{0 }\nu^{1/2}\le \delta_{0}\le \gamma_{0 }^{-1}\nu^{1/2}$
with a positive fixed constant $\gamma_{0}$.
Moreover, we assume
$0\le T^{-1}\eta\le 4,  T^{-1}\val{\eta'}\le 4\Lambda^{-1/2}$, 
$\val{\Theta(X)}\le C$, $\Theta$ real-valued.
Here $\Lambda, \mu, \nu$
are assumed to be positive constants such that
$\Lambda\ge \mu\ge \nu\ge 1$
\endsubsubhead
We start over our discussion from the identity (3.2.2):
$$
P= \re{\Bigl[q^w\w{\Bigl(\delta_{0}(1+T^{-1}\eta)+\Theta
\Bigr)}\Bigr]}.
\tag 3.2.8$$
We define
$$
a_{0}=\delta_{0}(1+T^{-1}\eta)
\tag 3.2.9$$
and
we note that
$
\gamma_{0} \nu^{1/2}\le a_{0}\le 5 \gamma_{0}^{-1}\nu^{1/2}.
$
\remark{Remark 3.2.1}{ 
We may assume that
$
\nu^{1/2}\ge 2C/\gamma_{0} 
$
which implies
$C\le \frac{1}{2}\gamma_{0} \nu^{1/2}$
so that
$$
\frac{1}{2}\gamma_{0} \nu^{1/2}\le a_{0}+\Theta\le (5 \gamma_{0}^{-1}
+C\gamma_{0}/2)\nu^{1/2}.
\tag 3.2.10$$
In fact if 
$
\nu^{1/2}<2C/\gamma_{0} 
$
we have
$\w{\bigl(\delta_{0}(1+T^{-1}\eta)+\Theta\bigr)}\in S(1,\Gamma)^w$,
$\Lambda^{1/2}\mu^{-1/2}q\in S(1,\Gamma)$
and 
$P\in S(\Lambda^{-1/2}\mu^{1/2},\Gamma)^w$
so that (3.2.7) holds also in that case.}
\endremark
We have 
the identity
$$q^w\w{\bigl(\delta_{0}(1+T^{-1}\eta)\bigr)}
=q^w\w{a_{0}}\quad\text{with}\
{\cases
\gamma_{0}\nu^{1/2} \le a_{0} \le 5 \gamma_{0}^{-1}\nu^{1/2},
\\
\val{a'_{0}}\le 10+\gamma_{0}^{-1}4\frac{\val{\delta_{0}}}{\Lambda^{1/2}}\le 10 +4\gamma_{0}^{-1}.
\endcases}
\tag 3.2.11$$
The Weyl symbol of
$\w{(a_{0}+\Theta)}$,
which is
$$a=(a_{0}+\Theta)\ast2^n\exp-2\pi\Gamma,
\tag 3.2.12$$
belongs to 
$S_{1}(\nu^{1/2},\nu^{-1}\Gamma)$(see definition A.5.1):
this follows from  the lemma A.5.3 and  (3.2.11) for $a_{0}\ast\exp-2\pi\Gamma$
and is obvious for
$\Theta\ast2^n\exp-2\pi\Gamma$
which belongs to $S(1,\Gamma)$.
Moreover the estimates (3.2.10) imply that   the symbol $a$ satisfies
$$
\frac{1}{2}\gamma_{0}\nu^{1/2}
\le a(X)=\int (a_{0}+\Theta)(X+Y)2^ne^{-2\pi {\val Y}^2}dY\le 
(5 \gamma_{0}^{-1}+C\gamma_{0}/2)\nu^{1/2}.
\tag 3.2.13$$
As a result,
the symbol $b=a^{1/2}$ belongs to $S_{1}(\nu^{1/4},\nu^{-1}\Gamma)$
and $1/b\in S_{1}(\nu^{-1/4},\nu^{-1}\Gamma)$:
we have 
$${2}^{-1/2}\gamma_{0}^{1/2}\nu^{1/4}\le \val b\le (5 \gamma_{0}^{-1}+C\gamma_{0}/2)^{1/2}\nu^{1/4}$$
and moreover
$
a'=a_{0}'\ast2^n\exp-2\pi\Gamma+\Theta\ast2^n(\exp-2\pi\Gamma)',
$
so that, using
$$
\val{a'}\le 10+4\gamma_{0}^{-1}+C\norm{2^n(\exp-2\pi\Gamma)'}_{L^1(\RZ)}=C_{1},
$$
we get
$
2\val{b'}=\val{a'(X)}{a(X)}^{-1/2}\le 2^{1/2}\gamma_{0}^{-1/2}\nu^{-1/4}C_{1},
$
and the derivatives of $a^{1/2}$ of order $k\ge 2$ are a sum of terms of type
$$
a^{\frac{1}{2}-m}a^{(k_{1})}\dots a^{(k_{m})},\quad\text{with $k_{1}+\dots+k_{m}=k,$ all $ k_{j}\ge 1$,}
$$
which can be estimated by
$
C\nu^{\frac{1}{4}-\frac{m}{2}}\le C\nu^{-\frac{1}{4}}
$
since $m\ge 1$.
Similarly we obtain that $b^{-1}\in S_{1}(\nu^{-1/4}, \nu^{-1}\Gamma).$
From the lemma A.5.2, we have
$
b^w b^w= a^w+S(\nu^{-1/2},\Gamma)^w=\w{(a_{0}+\Theta)}+S(\nu^{-1/2},\Gamma)^w,
$
which means
$
\w{(a_{0}+\Theta)}=b^wb^w+r_{0}^w,\quad r_{0}\in S(\nu^{-1/2},\Gamma),
\
\text{real-valued}.
$
Using  that $1/b$ belongs to $S_{1}(\nu^{-1/4},\nu^{-1}\Gamma)$, we write,
using again the lemma A.5.2, 
$$
\bigl(b+\frac{1}{2}b^{-1}r_{0}\bigr)^w\bigl(b+\frac{1}{2}b^{-1}r_{0}\bigr)^w
=
b^wb^w+r_{0}^w+S(\nu^{-1/4}\nu^{-1/2} \nu^{-1/4},\Gamma)^w,
$$
which gives, $$
\w{(a_{0}+\Theta)}=
\bigl(b+\frac{1}{2}b^{-1}r_{0}\bigr)^w\bigl(b+\frac{1}{2}b^{-1}r_{0}\bigr)^w
+S(\nu^{-1},\Gamma)^w.
\tag 3.2.14$$
Note that $b_{0}=b+\frac{1}{2}b^{-1}r_{0}$
belongs to $S_{1}(\nu^{1/4},\nu^{-1}\Gamma)$
since it is true for $b$
and $b^{-1}r_{0}\in S(\nu^{-3/4},\Gamma)$: we get then
$$
2\re\bigl(q^w\w{(a_{0}+\Theta)}\bigr)= 2 b_{0}^w q^wb_{0}^w+
\overbrace{[
\underbrace{[q^w,b_{0}^w]}_{S(\Lambda^{-1/2}\mu^{1/2}\nu^{1/2}\nu^{-1/4},\Gamma)}, b_{0}^w]}^{
S(\Lambda^{-1/2}\mu^{1/2}\nu^{1/4}\nu^{-1/4},\Gamma)}+\re(q^w S(\nu^{-1},\Gamma)^w)
$$
so that
$$
P=b_{0}^w q^wb_{0}^w
+
S(\Lambda^{-1/2}\mu^{1/2},\Gamma)^w.
\tag 3.2.15$$
Using now the Fefferman-Phong inequality (\cite{FP}, Theorem 18.6.8 in \cite{H6})  
for the nonnegative symbol $q$, we get
$
b_{0}^w q^wb_{0}^w=b_{0}^w (q^w+C\Lambda^{-1/2}\mu^{1/2}\nu^{-1})b_{0}^w+
S(\Lambda^{-1/2}\mu^{1/2}\nu^{-1/2},\Gamma)^w
\ge
S(\Lambda^{-1/2}\mu^{1/2}\nu^{-1/2},\Gamma)^w,
$
so that, from (3.2.15) we get eventually
$$
\re (QM)+S(\Lambda^{-1/2}\mu^{1/2},\Gamma)^w\ge 0.
\tag 3.2.16
$$
\subhead
3.3.  Stationary estimates 
\endsubhead
Let $T>0$ be given and
$Q(t)= q(t)^w$  given by (2.1.1-2). We define $M(t)$ according to (3.1.1).
We consider
$$
\re\bigl(Q(t)M(t)\bigr)
=\frac{1}{2}
Q(t)M(t)+\frac{1}{2}
M(t) Q(t)=P(t).
\tag 3.3.1$$
We have, omitting now the variable $t$ fixed throughout all this section 3.3, 
$$
P= \re{\Bigl[q^w\w{\bigl(\delta_{0}(1+T^{-1}\eta)\bigr)}
+q^w\w{\Theta
}\Bigr]}.
\tag 3.3.2$$
\proclaim{Lemma 3.3.1}
Let $p$ be the Weyl symbol of $P$ defined in $(3.3.2)$
and
$\widetilde\Theta=\Theta\ast2^n \exp-2\pi\Gamma$,
where $\Theta$ is defined in $(2.4.6)$
(and satisfies $(2.4.8)$).
Then we have
$$
p(t,X)\equiv p_{0}(t,X)=
q(t,X)\Bigl(\delta_{0}(1+T^{-1}\eta)\ast2^n \exp-2\pi\Gamma\Bigr)+q(t,X)\widetilde\Theta(t,X),
\tag 3.3.3$$
modulo
$S(
\Lambda^{-1/2}     \mu^{1/2}  \nu^{-1/2}\langle\delta_{0}\rangle,\Gamma
)$.
\endproclaim
\demo{Proof}
Using the results of section 2.1, we know that the symbol
$X\mapsto q(t,X)$
belongs to the class
$$S(\Lambda^{-1/2}\mu(t,X)^{1/2}\nu(t,X), \nu(t,X)^{-1}\Gamma)$$
as shown in lemma 2.1.7.
In fact
from (2.1.18) we know that 
$
q\in S(\Lambda^{-1/2}\mu^{1/2}\nu,\nu^{-1}\Gamma),$
and from (3.1.3) and Lemma A.1.3, we obtain,
using Theorem 18.5.5
in
\cite{H6}, 
$$
q\sharp \widetilde\Theta= q\widetilde\Theta+\frac{1}{4i\pi}\poi{q}{\widetilde\Theta}
+
S(
\Lambda^{-1/2}     \mu^{1/2}  \nu^{-1/2}\langle\delta_{0}\rangle,\Gamma
).
$$
This implies that 
$
\re{(q\sharp \widetilde\Theta)}\in q\widetilde\Theta+S(
\Lambda^{-1/2}     \mu^{1/2}  \nu^{-1/2}\langle\delta_{0}\rangle,\Gamma
).
$
On the other hand, we know that
$$
\re\Bigl(
q\sharp \overbrace{\bigl[\delta_{0}(1+T^{-1}\eta)\ast \exp-2\pi\Gamma \bigl]}^{\omega}
\Bigr)=q\omega+\sum_{\val \alpha=\val \beta=2}c_{{\alpha\beta}}q^{(\alpha)}\omega^{(\beta)}
+S(\Lambda^{-1/2}     \mu^{1/2} \nu^{-1},\Gamma)
$$
so that it is enough to concentrate our attention
on the ``products" $q'' \omega''$.
We have 
$$\bigl(\delta_{0}(1+T^{-1}\eta)\bigr)''\ast\exp-2\pi\Gamma\in S(1,\Gamma) $$
and since
$q''\in S(\Lambda^{-1/2}\mu^{1/2},\nu^{-1}\Gamma),
$
we get a remainder in 
$
S(\Lambda^{-1/2}\mu^{1/2},\Gamma),
$
which is fine as long as
$\langle\delta_0\rangle\ge c \nu^{1/2}$.
However when
$\langle\delta_0\rangle\le c \nu^{1/2}$, we know that, for a good choice of the fixed positive constant
$c$,
the function $\delta_{0} $
satisfies the estimates of
$S(\nu^{1/2},\nu^{-1}\Gamma),$
since it is the $\Gamma$-distance function
to the set of (regular) zeroes of the function $q$
so that
$q'' \delta_{0}''\in S(\Lambda^{-1/2}\mu^{1/2}\nu^{-1/2},\nu^{-1}\Gamma)$
which is what we are looking for. However,
we are  left 
with
$$
q'' (\delta_{0}\eta\ast\exp-2\pi\Gamma)''T^{-1}.
$$
Since we have
$
(\delta_{0}\eta)''=\delta_{0}''\eta
+2\delta_{0}'\eta'
+\delta_{0}\eta''
$
and
$
\val{\delta_{0}''\eta
+2\delta_{0}'\eta'
}\le CT( \nu^{-1/2}+\Lambda^{-1/2}),
$
we have only to deal with the term
$$\multline
\delta_{0}\eta''\ast\exp-2\pi\Gamma=\int
\delta_{0}(Y)\eta''(Y)\exp-2\pi\Gamma(X-Y) dY
\\=
-\int
\underbrace{\delta_{0}'(Y)\eta'(Y)}_{\lesssim T \Lambda^{-1/2}}\exp-2\pi\Gamma(X-Y) dY
-\int
\underbrace{\delta_{0}(Y)\eta'(Y)}_{\lesssim T\Lambda^{-1/2} \langle \delta_{0}\rangle}4\pi(X-Y)\exp-2\pi\Gamma(X-Y) dY.
\endmultline
$$
For future reference we summarize part of the previous discussion by the following
result.
\proclaim{Lemma 3.3.2} With the notations above, we have
$$\gather
\Val{\Bigl(\delta_{0}(1+T^{-1}\eta)\ast \exp-2\pi\Gamma\Bigr)}\le C
\langle\delta_{0}\rangle,
\quad
\Val{\Bigl(\delta_{0}(1+T^{-1}\eta)\ast \exp-2\pi\Gamma\Bigr)'}\le C,
\\
\Val{\Bigl(\delta_{0}(1+T^{-1}\eta)\ast \exp-2\pi\Gamma\Bigr)''}\le C
\langle\delta_{0}\rangle\nu^{-1/2}.
\endgather
$$
\endproclaim
\demo{Proof}
Starting over the discussion, we have already seen that the
result is true whenever
$\langle\delta_{0}\rangle\gtrsim \nu^{1/2}$.
Moreover when
$\langle\delta_{0}\rangle\ll\nu^{1/2}$, 
we have seen that
$\val{\delta_{0}''}\lesssim  \nu^{-1/2}$ and $T^{-1}\val \eta\lesssim 1$;
moreover we have already checked
$\val{\eta'}\lesssim T\Lambda^{-1/2}$
and 
$
T^{-1}\val{\delta_{0}'\eta'
}\lesssim \Lambda^{-1/2} \lesssim \nu^{-1/2}
$
as well as
$
\val{\delta_{0}\eta''\ast\exp-2\pi\Gamma}\lesssim\Lambda^{-1/2} \langle\delta_{0}\rangle
\lesssim  \langle\delta_{0}\rangle \nu^{-1/2}
$.
\qed
\enddemo
Eventually, using the lemma A.1.3,
we get that the first integral above
is in 
$S(T\Lambda^{-1/2},\Gamma)$
whereas the second belongs to
$S(T\Lambda^{-1/2}\langle \delta_{0}\rangle,\Gamma)$. Finally,
it means that, up to terms in 
$S(\Lambda^{-1/2}\mu^{1/2}\nu^{-1/2}\langle \delta_{0}\rangle,\Gamma),$
the operator $P(t)$ has a Weyl symbol equal to 
the rhs of (3.3.3).\qed
\enddemo
We shall use a partition of unity $1=\sum_{k}\chi_{k}^2$ related to the metric
$ \nu(t,X)^{-1}\Gamma$ and a sequence $(\psi_{k})$
as in section 1.4.
We have, omitting the variable $t$,
with $p_{0}$ defined in (3.3.3),
$$\multline
p_{0}(X)=\sum_{k} \chi_{k}(X)^2q(X)\int \delta_{0}(Y)\bigl(1+T^{-1}\eta(Y)\bigr)2^n\exp-2\pi \Gamma(X-Y) dY
\\+ \sum_{k} \chi_{k}(X)^2q(X)\int \Theta(Y)2^n\exp-2\pi \Gamma(X-Y)dY.
\endmultline$$
Using the lemma A.1.6, we obtain,
assuming
$\delta_{0}=\delta_{0k}, \Theta=\Theta_{k}, q=q_{k}$ on $U_{k}$
$$\multline
p_{0}=\sum_{k}\chi_{k}^2q_{k}\bigl(\delta_{0k}(1+T^{-1}\eta)\ast 2^n \exp-2\pi\Gamma\bigr)
+
\sum_{k}\chi_{k}^2q_{k} \bigl(\Theta_{k}\ast 2^n \exp-2\pi\Gamma\bigr)
+ S(\Lambda^{-1/2}\mu^{1/2}\nu^{-\io},\Gamma).
\endmultline
\tag 3.3.4$$
\proclaim{Lemma 3.3.3}
With 
$\widetilde{\Theta}_{k}=\Theta_{k}\ast 2^n \exp-2\pi\Gamma,\quad
d_{k}=
\delta_{0k}(1+T^{-1}\eta)\ast 2^n \exp-2\pi\Gamma
$
and
$q_{k},\chi_{k}$ defined above,
we have
$$
\sum_{k}\chi_{k} \sharp q_{k} d_{k}\sharp \chi_{k}+ \sum_{k}\chi_{k} \sharp q_{k} \widetilde{\Theta}_{k}\sharp \chi_{k}
=p_{0}
+ S(\Lambda^{-1/2}\mu^{1/2}\nu^{-1/2}
\langle \delta_{0}\rangle
,\Gamma).
\tag 3.3.5$$
\endproclaim
\demo{Proof}
We already know that
$\val{d_{k}}\lesssim \langle \delta_{0}\rangle,\quad
\val{d_{k}'}\lesssim 1,\quad
\val{d''_{k}}\lesssim \langle \delta_{0}\rangle\nu^{-1/2},
$
so that
$$
\val{(q_{k} d_{k})''=q''_{k}d_{k}+2q'_{k}d'_{k}+q_{k}d''_{k}}
\lesssim
\Lambda^{-1/2}\mu^{1/2}\bigl(\langle \delta_{0}\rangle+\nu^{1/2}+\nu^{1/2} \langle \delta_{0}\rangle\bigr)
\lesssim
\Lambda^{-1/2}\mu^{1/2}\nu^{1/2} \langle \delta_{0}\rangle.
\tag 3.3.6$$
As a consequence, we get 
$$\align
\sum_{k}\chi_{k} \sharp q_{k} d_{k}\sharp \chi_{k}
&=\sum_{k}
\Bigl(\chi_{k}  q_{k} d_{k}+\frac{1}{4i\pi}
\poi{\chi_{k}}{q_{k} d_{k}}
+S(\nu^{-1}(
\Lambda^{-1/2}\mu^{1/2}\langle \delta_{0}\rangle \nu^{1/2}
),\Gamma)\Bigr)\sharp \chi_{k}
\\
&=\sum_{k}\Bigl(\chi_{k}  q_{k} d_{k}+\frac{1}{4i\pi}
\poi{\chi_{k}}{q_{k} d_{k}}\Bigr)\sharp \chi_{k}+
\sum_{k} S(
\Lambda^{-1/2}\mu^{1/2}\langle \delta_{0}\rangle \nu^{-1/2}
,\Gamma)\sharp \chi_{k}
\\
&=\sum_{k}\Bigl(\chi_{k}  q_{k} d_{k}+\frac{1}{4i\pi}
\poi{\chi_{k}}{q_{k} d_{k}}\Bigr)\chi_{k}+\frac{1}{4i\pi}\sum_{k}
\poi{\chi_{k}  q_{k} d_{k}+\frac{1}{4i\pi}
\poi{\chi_{k}}{q_{k} d_{k}}}{\chi_{k}}
\\
& \hbox to 200pt{} +
S(
\Lambda^{-1/2}\mu^{1/2}\langle \delta_{0}\rangle \nu^{-1/2}
,\Gamma)
\endalign
$$
since
$
\val{(\chi_{k}q_{k} d_{k})''\chi_{k}''}\lesssim
\Lambda^{-1/2}\mu^{1/2}(\langle \delta_{0}\rangle
+
\nu^{1/2}
+\nu\langle \delta_{0}\rangle\nu^{-1/2}) \nu^{-1}\lesssim
\langle \delta_{0}\rangle\nu^{-1/2}\Lambda^{-1/2}\mu^{1/2}.
$
Using now that $\chi_{k} \sharp q_{k} d_{k}\sharp \chi_{k}$ is real-valued,
we obtain 
$$
\sum_{k}\chi_{k} \sharp q_{k} d_{k}\sharp \chi_{k}
=\sum_{k}\chi_{k}^2  q_{k} d_{k}
-\frac{1}{16\pi^2}\sum_{k}
\poi{
\poi{\chi_{k}}{q_{k} d_{k}}}{\chi_{k}}
+
S(
\Lambda^{-1/2}\mu^{1/2}\langle \delta_{0}\rangle \nu^{-1/2}
,\Gamma).
\tag 3.3.7$$
We note now that, using (3.3.6), we have
$$
\poi{
\poi{\chi_{k}}{q_{k} d_{k}}}{\chi_{k}}=-H_{{\chi_{k}}}^2(q_{k} d_{k})
\in S(\Lambda^{-1/2}\mu^{1/2}\langle \delta_{0}\rangle\nu^{1/2}\nu^{-1},\Gamma).
\tag 3.3.8
$$
We examine now the term
$$
\chi_{k}\sharp q_{k}\widetilde{\Theta}_{k}\sharp\chi_{k}
=\bigl(
\chi_{k}q_{k}\widetilde{\Theta}_{k}
\bigr)\sharp \chi_{k}+\frac{1}{4i\pi}\poi{\chi_{k}}{q_{k}\widetilde{\Theta}_{k}}\sharp\chi_{k}
+
S(\nu^{-1}\Lambda^{-1/2}\mu^{1/2}
\nu\langle \delta_{0}\rangle\nu^{-1/2},\Gamma)\sharp\chi_{k}.
$$
We have
$$\align
\chi_{k}q_{k}\widetilde{\Theta}_{k}\sharp\chi_{k}
&\in \chi_{k}^2q_{k}\widetilde{\Theta}_{k}+S(\nu^{-1}
\Lambda^{-1/2}\mu^{1/2}
\nu\langle \delta_{0}\rangle\nu^{-1/2},\Gamma
),
\\
\frac{1}{4i\pi}\poi{\chi_{k}}{q_{k}\widetilde{\Theta}_{k}}&\in i\R
+
S(\nu^{-1/2}
\Lambda^{-1/2}\mu^{1/2}
\nu\langle \delta_{0}\rangle\nu^{-1/2},\Gamma).
\endalign
$$
Since $\chi_{k}\sharp q_{k}\widetilde{\Theta}_{k}\sharp\chi_{k}
$ is real-valued, we get
$$
\sum_{k}\chi_{k}\sharp q_{k}\widetilde{\Theta}_{k}\sharp\chi_{k}=\sum_{k}
\chi_{k}^2q_{k}\widetilde{\Theta}_{k}
+
S(\Lambda^{-1/2}\mu^{1/2}
\langle \delta_{0}\rangle\nu^{-1/2},\Gamma).\tag 3.3.9
$$
Collecting the information
(3.3.4), (3.3.7), (3.3.8) and (3.3.9) we obtain (3.3.5)
and the lemma.\qed
\enddemo
From this lemma and the lemma 3.3.1
we obtain that
$$
\re\bigl(Q(t)M(t)\bigr)
=\sum_{k}\chi_{k}^w\bigl(q_{k} d_{k}+q_{k} \widetilde{\Theta}_{k}\bigr)^w\chi_{k}^w
+ 
S(\Lambda^{-1/2}\mu^{1/2}
\langle \delta_{0}\rangle\nu^{-1/2},\Gamma)^w.
\tag 3.3.10$$
Moreover the same arguments as above in Lemma 3.3.1 give also that
$$
\re(q_{k }^wd_{k}^w+q_{k}^w \widetilde{\Theta}_{k}^w)
=(q_{k} d_{k}+q_{k}\widetilde{\Theta}_{k})^w
+
S(\Lambda^{-1/2}\mu^{1/2}
\langle \delta_{0}\rangle\nu^{-1/2},\Gamma)^w.\tag 3.3.11
$$
\proclaim{Proposition 3.3.4}
Let $T>0$ be given and
$Q(t)= q(t)^w$  given by $(2.1.1$-$2)$. We define $M(t)$ according to $(3.1.1)$.
Then, with a partition of unity $1=\sum_{k}\chi_{k}^2$
related to the metric
$\nu(t,X)^{-1}\Gamma$
we have
$$\gather
\re{(Q(t)M(t))}=
\sum_{k}\chi_{k}^w\re\bigl(q_{k}^w d_{k}^w+q_{k}^w\widetilde{\Theta}_{k}^w\bigr)\chi_{k}^w
+ 
S(\Lambda^{-1/2}\mu^{1/2}
\langle \delta_{0}\rangle\nu^{-1/2},\Gamma)^w
\tag 3.3.12
\\\text{and  }\quad
\re{(Q(t)M(t))}+ 
S(\Lambda^{-1/2}\mu^{1/2}
\langle \delta_{0}\rangle\nu^{-1/2},\Gamma)^w\ge 0.
\tag 3.3.13
\endgather
$$
\endproclaim
\demo{Proof}
The equality (3.3.12) follows from (3.3.10-11).
According to Lemma 2.1.9,
we have to deal with four subsets of indices, $E_{\pm},  E_{0}, E_{00}$.
The classification in Definition 2.1.8 shows that
section 3.2.[1]
takes care
of the cases $E_{0}$
and shows that, from (3.2.7),
$$
\text{for $k\in E_{0}$,\quad}
\re\bigl(q_{k}^w d_{k}^w+q_{k}^w\widetilde{\Theta}_{k}^w)+S(\Lambda^{-1/2}\mu^{1/2}
\nu^{-1/2},\Gamma)^w\ge 0.\tag 3.3.14
$$
Furthermore, the estimate (3.2.16) in section 3.2.[2]  shows that
$$
\text{for $k\in E_{\pm}$,\quad}
\re\bigl(q_{k}^w d_{k}^w+q_{k}^w\widetilde{\Theta}_{k}^w)+S(\Lambda^{-1/2}\mu^{1/2}
\nu^{-1/2}\langle \delta_{0}\rangle,\Gamma)^w\ge 0.\tag 3.3.15
$$
Moreover if $k\in E_{00}$, the weight
$\nu$ is bounded above and 
$$
q_{k}^w d_{k}^w+q_{k}^w\widetilde{\Theta}_{k}^w\in
S(\Lambda^{-1/2}\mu^{1/2}
\nu^{-1/2},\Gamma) ^w.
\tag 3.3.16
$$
The equality (3.3.12) and (3.3.14-15-16)
give (3.3.13).
\qed
\enddemo
\subhead
3.4. The multiplier method
\endsubhead
\proclaim{Theorem 3.4.1}
Let $T>0$ be given and
$Q(t)= q(t)^w$  given by $(2.1.1$-$2)$. We define $M(t)$ according to $(3.1.1)$.
There exist $T_{0}>0$ and $c_{0}>0$ depending only on a finite number of $\gamma_{k}$
in $(2.1.1)$
such that, for $0<T\le T_{0}$,
with $D(t,X)=\langle
\delta_{0}(t,X)\rangle 
$,
($D$ is Lipschitz continuous with Lipschitz constant 2, as $\delta_{0}$ in $(2.1.8)$
and thus a $\Gamma$-weight),
$$
\frac{d}{dt} M(t)+2\re{\bigl(Q(t)M(t)\bigr)}\ge
T^{-1}\w{(D^2)}\Lambda^{-1/2}c_{0}.
\tag 3.4.1$$
Moreover we have with $m$ defined in $(2.4.7)$, $\widetilde m(t,\cdot)=m(t,\cdot)\ast 2^n\exp-2\pi\Gamma$, \quad
$$\gather
M(t)= \w{m(t,X)}=
\widetilde m(t,X)^w,\
\text{with $\widetilde m\in S_{1}(D, D^{-2}
\Gamma)
+S(1,
\Gamma)$}.\tag 3.4.2
\\
m(t,X)= a(t,X)+b(t,X), \
\val{a/D}+\val{a'_{X}}+\val{b} \text{ bounded},
\
\dot m\ge 0,
\tag 3.4.3
\\
a=\delta_{0}(1+T^{-1}\eta),\quad b=\widetilde{\Theta}.
\endgather
$$
\endproclaim
\demo{Proof}
From the estimate
(3.1.2), we get, with a positive fixed constant $C_{0}$, 
$$
\frac{d}{dt}M(t)\ge 
\frac{1}{2C_{0}T} \w{(\Lambda^{-1/2}\mu^{1/2}
\nu^{-1/2}
\langle\delta_{0} \rangle)
}+T^{-1}\w{{(\delta_{0}^2 )}}\Lambda^{-1/2},
$$
and from
(3.3.13)
and Lemma A.1.4 we know that,
with a fixed (nonnegative) constant $C_{1}$,
$$
2\re{\bigl(Q(t)M(t)\bigr)}+ C_{1}
\w{(\Lambda^{-1/2}\mu^{1/2}
\nu^{-1/2}
\langle\delta_{0} \rangle)
}\ge0.
$$
As a result we get, if $4C_{1}C_{0}T\le 1$ (we shall choose $T_{0}=\frac{1}{4C_{0}({C_{1}}+1)}$),
$$
\frac{d}{dt} M(t)+2\re{\bigl(Q(t)M(t)\bigr)} 
\ge 
\frac{1}{4C_{0}T} \w{(\Lambda^{-1/2}\mu^{1/2}
\nu^{-1/2}
\langle\delta_{0} \rangle)
}+T^{-1}\w{{(\delta_{0}^2 )}}\Lambda^{-1/2}.
$$
Using (2.1.17)( $\mu\ge\nu/2$), this gives
$$
\frac{d}{dt} M(t)+2\re{\bigl(Q(t)M(t)\bigr)} 
\ge
T^{-1}\Lambda^{-1/2}\w{(1+\delta_{0}^2)}\bigl(\frac{1}{2^{5/2}C_{0}+1}\bigr),
$$
which is the sought  result.
\qed
\enddemo
\head
4. From semi-classical to local estimates
\endhead
\subhead
 4.1. From semi-classical to inhomogeneous estimates
\endsubhead
Let us consider
a smooth real-valued function $f$
defined on $\R\times \R^n\times \R^n$,
satisfiying
(2.1.2) and such that, for all multi-indices $\alpha, \beta$,
$$
\sup_{t\in \R\atop
(x,\xi)\in\RZ}\val{(\p_{x}^\alpha\partial_{\xi}^\beta f)(t,x,\xi)}(1+\val \xi)^{-1+\val\beta}= C_{{\alpha\beta}}<\io.\tag 4.1.1
$$
Using a Littlewood-Paley
decomposition, we have
$$\gather
f(t,x,\xi)=\sum_{j\in \N} f(t,x,\xi)\varphi_{
j}(\xi)^2,
\text{
$\supp{\varphi_{0}}$ compact}, 
\\
\text{for $j\ge 1$, 
$\supp{\varphi_{j}}\subset \{\xi\in \R^n, 2^{j-1}\le \val \xi\le
2^{j+1}\}$,
$\sup_{j,\xi}\val{\partial_{\xi}^\alpha\varphi_{j}(\xi)}2^{j\val\alpha}<\io$.
}
\endgather
$$
We 
introduce also some  smooth nonnegative compactly supported functions $\psi_{j}(\xi)$,
satisfying the same 
estimates than $\varphi_{j}$
and supported in $2^{j-2}\le \val\xi\le 2^{j+2}$ for $j\ge 1$,
identically 1 on the support of $\varphi_{j}$.
For each $j\in \N$, we define the symbol
$$
q_{j}(t,x,\xi)=f(t,x,\xi)\psi_{j}(\xi)
\tag 4.1.2$$
and we remark that (2.1.2) is satisfied for $q_{j}$
and the following estimates hold:
$
\val{(\p_{x}^\alpha\partial_{\xi}^\beta q_{j})}\le C'_{{\alpha\beta}}\Lambda_{j}^{1-\val\beta},
$
with $\Lambda_{j}=2^j$.
Note that the semi-norms of $q_{j}$ can be estimated from above
independently of $j$.
We can reformulate this by saying that
$$
q_{j}\in S(\Lambda_{j},\Lambda_{j}^{-1}\Gamma_{j}),\quad\text{with $\Gamma_{j}(t,\tau)=\val{t}^2\Lambda_{j}+\val{\tau}^2\Lambda_{j}^{-1}$}
\quad(\text{note that $\Gamma_{j}=\Gamma_{j}^\sigma$}).
\tag 4.1.3$$
\proclaim{Lemma 4.1.1}
There
exists $T_{0}>0, c_{0}>0,$
depending only on a finite number
of
semi-norms of $f$
such that, 
for each $j\in \N$, we can find
$ D_{j}$
a $\Gamma_{j}$--uniformly Lipschitz continuous function 
with Lipschitz constant $2$,
valued in $[1, \sqrt{2\Lambda_{j}}]$,
$a_{j}, \ b_{j}$
real-valued
such that 
$$
\sup_{j\in \N, \val t\le T_{0}\atop X\in \RZ}
\left(\Val{\frac{a_{j}(t,X)}{D_{j}(t,X)}}
+\norm{\nabla_{X}a_{j}(t,X)}_{\Gamma_{j}}+\val{b_{j}(t,X)}\right)<\io.
\tag 4.1.4$$
Moreover we have with 
$m_{j}=a_{j}+b_{j}
$,\quad
$\widetilde m_{j}(t,\cdot)=m_{j}(t,\cdot)\ast 2^n\exp-2\pi\Gamma_{j}$, \quad
$Q_{j}(t)=q_{j}(t)^w,$
$$
M_{j}(t)= m_{j}(t,X)^{
\text{Wick($\Gamma_{j}$)}
}
=
\widetilde m_{j}(t,X)^w,\
\text{with $\widetilde m_{j}\in S_{1}(D_{j}, D_{j}^{-2}
\Gamma_{j})
+S(1,
\Gamma_{j})$},\tag 4.1.5
$$
(the {\rm Wick($\Gamma_{j}$)}
quantization is defined in definition A.1.7) the estimate
$$
\frac{d}{dt} M_{j}(t)+2\re{\bigl(Q_{j}(t)M_{j}(t)\bigr)}\ge
T^{-1}\bigl({D_{j}^2}\bigr)^{\text{Wick($\Gamma_{j}$)}}\Lambda_{j}^{-1/2}c_{0}.
\tag 4.1.6$$
\endproclaim
\demo{Proof} It is a straightforward consequence of Definition A.1.7
and of Theorem 3.4.1: let us check this.
Considering the linear symplectic mapping $L:$
$(t,\tau)\mapsto (\Lambda_{j}^{-1/2}t,\Lambda_{j}^{1/2}\tau)$,
we see that the symbols
$q_{j}\circ L$
belong uniformly to 
$S(\Lambda_{j},\Lambda_{j}^{-1}\Gamma_{0})$.
Applying the theorem 3.4.1 to $q_{j}\circ L$, 
we find $D$
a $\Gamma_{0}$--uniformly Lipschitz continuous function $\ge 1$,
$a, b$
real-valued
such that 
$$
\sup_{j\in \N, \val t\le T_{0}\atop X\in \RZ}
\left(\Val{\frac{a(t,X)}{D(t,X)}}
+\norm{\nabla_{X}a(t,X)}_{\Gamma_{0}}+\val{b(t,X)}\right)<\io,
\tag 4.1.7$$
and so  that, with 
$m=a+b
$,\quad
$\widetilde m(t,\cdot)=m(t,\cdot)\ast 2^n\exp-2\pi\Gamma_{0}$, \quad
$Q(t)=(q_{j}(t)\circ L)^w,$
$$
M(t)= m(t,X)^{
\text{Wick}
}
=
\widetilde m(t,X)^w,\
\text{with $\widetilde m\in S_{1}(D, D^{-2}
\Gamma_{0})
+S(1,
\Gamma_{0})$},$$
$$
\frac{d}{dt} M(t)+2\re{\bigl(Q(t)M(t)\bigr)}\ge
T^{-1}\bigl({D^2}\bigr)^{\text{Wick($\Gamma_{0}$)}}\Lambda_{j}^{-1/2}c_{0}.
\tag 4.1.8$$
Now we define the real-valued functions
$a_{j} =a\circ L^{-1}, b_{j}=b\circ L^{-1}, D_{j}=D\circ L^{-1}$
and we have, since
$\Gamma_{0}(S)=\Gamma_{j}(LS),$
$$\multline
\Val{\frac{a_{j}(t,X)}{D_{j}(t,X)}}
+\norm{\nabla_{X}a_{j}(t,X)}_{\Gamma_{j}}+\val{b_{j}(t,X)}
\\
=
\Val{\frac{a(t,L^{-1}X)}{D(t,L^{-1}X)}}
+\sup_{T\in\RZ}\frac{\val{a'_{j}(t,X)\cdot T}}{\Gamma_{j}(T)^{1/2}}+\val{b(t,L^{-1}X)}
\\
\hskip85pt=
\Val{\frac{a(t,L^{-1}X)}{D(t,L^{-1}X)}}
+\sup_{T\in\RZ}\frac{\val{a'(t,X)\cdot L^{-1}T}}{\Gamma_{j}(T)^{1/2}}+\val{b(t,L^{-1}X)}
\\
=
\Val{\frac{a(t,L^{-1}X)}{D(t,L^{-1}X)}}
+\norm{a'(t,X)}_{\Gamma_{0}}+\val{b(t,L^{-1}X)},
\endmultline
$$
so that (4.1.7) implies (4.1.4).
Considering now $m_{j} =a_{j}+b_{j}$
and for a metaplectic $U$ in the fiber of the symplectic $L$
(see definition A.1.7), we have
$$
M_{j}(t)=m_{j}(t,X)^{\text{Wick($\Gamma_{j}$)}}=
U(m_{j}\circ L)^{\text{Wick($\Gamma_{0}$)}}U^*.
\tag 4.1.9$$
Thus we obtain
$$\align
\frac{d}{dt}& M_{j}(t)+2\re{\bigl(Q_{j}(t)M_{j}(t)\bigr)}
\\
\text{\sevenrm  from (4.1.9.)\quad}&=
U\frac{d}{dt} 
(m_{j}\circ L)^{\text{Wick($\Gamma_{0}$)}}U^*\hskip-1pt+\hskip-1pt
2\re{\Bigl(UU^*q_{j}(t)^wU(m_{j}\circ L)^{\text{Wick($\Gamma_{0}$)}}U^*
\Bigr)}
\\
\text{\sevenrm \qquad   }&=
U\Bigl[\frac{d}{dt} 
(m_{j}\circ L)^{\text{Wick($\Gamma_{0}$)}}+2\re{\Bigl(U^*q_{j}(t)^wU(m_{j}\circ L)^{\text{Wick($\Gamma_{0}$)}}
\Bigr)}\Bigr] U^*
\\\text{\sevenrm using $\scriptstyle 
\left.\matrix \scriptstyle 
m=m_{j}\circ L\\ 
\scriptstyle 
(q\circ L)^w=U^* q^w U
\endmatrix\right]
$\quad   }
&=
U\Bigl[\frac{d}{dt} 
(m)^{\text{Wick($\Gamma_{0}$)}}+2\re{\Bigl((q_{j}\circ L)^w(m)^{\text{Wick($\Gamma_{0}$)}}
\Bigr)}\Bigr] U^*
\\
\text{\sevenrm from (4.1.8)\quad }&\ge 
U\Bigl[T^{-1}\bigl({D^2}\bigr)^{\text{Wick($\Gamma_{0}$)}}\Lambda_{j}^{-1/2}c_{0}\Bigr] U^*
\\
\text{\sevenrm from (A.1.16)\quad }&=
T^{-1}UU^*\bigl({D^2}\circ L^{-1}\bigr)^{\text{Wick($\Gamma_{j}$)}}UU^*\Lambda_{j}^{-1/2}c_{0}
\\
&=
T^{-1}\bigl({D_{j}^2}\bigr)^{\text{Wick($\Gamma_{j}$)}}
\Lambda_{j}^{-1/2}c_{0},\endalign$$
which is (4.1.6), completing the proof of the lemma.
\qed
\enddemo
We define now, with $\varphi_{j}$ given after  (4.1.1), $M_{j}$ in (4.1.5)
$$
\Cal M(t)=\sum_{j\in\N}\varphi_{j}^w \Lambda_{j}^{-1/2}M_{j}(t)\varphi_{j}^w.
\tag 4.1.10$$
\proclaim{Lemma 4.1.2}
With $M_{j}$ defined in $(4.1.5)$
and $\varphi_{j}, \psi_{j}$
as above, 
$$
\sum_{j}\varphi_{j}^w M_{j}(t)\bigl((1-\psi_{j})f(t)\bigr)^w\varphi_{j}^w \in
S(\langle\xi\rangle^{-\io},\langle\xi\rangle\val{dx}^2+\langle\xi\rangle^{-1}\val{d\xi}^2)^w,
\tag 4.1.11$$
$$
\sum_{j}\varphi_{j}^w M_{j}(t)\varphi_{j}^w \bigl((1-\psi_{j})f(t) \bigr)^w
\in
S(\langle\xi\rangle^{-\io},\langle\xi\rangle\val{dx}^2+\langle\xi\rangle^{-1}\val{d\xi}^2)^w. 
\tag 4.1.12$$\endproclaim
\demo{Proof}
Since
$\psi_{j}\equiv1$
on the support of $\varphi_{j}$,
we get that, uniformly with respect to $j$, 
$$\text{
$\bigl((1-\psi_{j})f(t)\bigr)^w\varphi_{j}^w $
$\in$
$S(\langle\xi\rangle^{-\io},\val{dx}^2+\langle\xi\rangle^{-2}\val{d\xi}^2)^w$.}
\tag 4.1.13$$
Since
$
\widetilde{m_{j}}\in 
S(\Lambda_{j}^{1/2},\Lambda_{j}\val{dx}^2+\Lambda_{j}^{-1}\val{d\xi}^2),
$
we get that
$
\psi_{j }\widetilde{m_{j}}\in  S(\langle\xi\rangle^{1/2},
\langle\xi\rangle\val{dx}^2+\langle\xi\rangle^{-1}\val{d\xi}^2),
$ and consequently
$
\varphi_{j}\sharp\psi_{j }\widetilde{m_{j}}\in S(\langle\xi\rangle^{1/2},
\langle\xi\rangle\val{dx}^2+\langle\xi\rangle^{-1}\val{d\xi}^2)
$
so that 
$$
\varphi_{j}\sharp\psi_{j}\widetilde{m_{j}}\sharp(1-\psi_{j})f(t)\sharp \varphi_{j}
\in S(\langle\xi\rangle^{-\io},\langle\xi\rangle\val{dx}^2+\langle\xi\rangle^{-1}\val{d\xi}^2)
\subset 
S(\langle\xi\rangle^{-\io},\val{dx}^2+\val{d\xi}^2).
\tag 4.1.14$$
Moreover we have 
$\varphi_{j}\sharp(1-\psi_{j})\widetilde{m_{j}}\in S(\Lambda_{j}^{-\io},\Lambda_{j}\val{dx}^2+\Lambda_{j}^{-1}\val{d\xi}^2)
\subset 
S(\Lambda_{j}^{-\io},\val{dx}^2+\val{d\xi}^2)\
$
so that (4.1.13) implies
$$
\varphi_{j}\sharp(1-\psi_{j})\widetilde{m_{j}}\sharp(1-\psi_{j})f(t)\sharp \varphi_{j}
\in S(\langle\xi\rangle^{-\io},\val{dx}^2+\val{d\xi}^2)
\subset
S(\langle\xi\rangle^{-\io},\langle\xi\rangle\val{dx}^2+\langle\xi\rangle^{-1}\val{d\xi}^2).
\tag 4.1.15$$
As a consequence, from (4.1.14) and (4.1.15) we get, uniformly in $j$,  that 
$$
\varphi_{j}\sharp\widetilde{m_{j}}\sharp(1-\psi_{j})f(t)\sharp \varphi_{j}\in 
S(\langle\xi\rangle^{-\io},\langle\xi\rangle\val{dx}^2+\langle\xi\rangle^{-1}\val{d\xi}^2).
\tag 4.1.16$$
Since $\varphi_{j}, \psi_{j}$ depend only on the variable $\xi$, the support condition implies
$
\varphi_{j}^w\psi_{j}^w=\varphi_{j}^w
$
and we obtain that from (4.1.16)
$$
\sum_{j}\varphi_{j}\sharp\widetilde{m_{j}}\sharp(1-\psi_{j})f(t)\sharp \varphi_{j}
=
\sum_{j}\psi_{j}\sharp\varphi_{j}\sharp\widetilde{m_{j}}\sharp(1-\psi_{j})f(t)\sharp \varphi_{j}
\sharp \psi_{j}
\in
S(\langle\xi\rangle^{-\io},\langle\xi\rangle\val{dx}^2+\langle\xi\rangle^{-1}\val{d\xi}^2),
$$
completing the proof of (4.1.11). The proof of (4.1.12) follows almost in the same way: we get as in (4.1.16) that
$$
\varphi_{j}\sharp\widetilde{m_{j}}\sharp \varphi_{j}\sharp(1-\psi_{j})f(t)\in 
S(\langle\xi\rangle^{-\io},\langle\xi\rangle\val{dx}^2+\langle\xi\rangle^{-1}\val{d\xi}^2).
$$
Now with
$\Phi_{j}=\varphi_{j}\sharp(1-\psi_{j})f(t)$,
we have
$\Phi_{j}\in
S(\langle\xi\rangle^{-\io},\val{dx}^2+\langle\xi\rangle^{-2}\val{d\xi}^2)
$
and
from the formula (A.5.5) we have also
$\val{(\p_{x}^\alpha\p_{\xi}^\beta\Phi_{j})(x,\xi)}
\le C_{\alpha\beta N}2^{jn}(1+\val{\xi-\supp \varphi_{j}})^{-N}(1+\val \xi),
$
so that
$$
\val{(\p_{x}^\alpha\p_{\xi}^\beta\Phi_{j})(x,\xi)}
\le
\cases C_{\alpha\beta N}2^{jn} 2^{-j(N-1)}&\text{if  $\val{\xi}\ge 2^{j+2},$}
\\
C_{\alpha\beta N}2^{jn} 2^{-jN}&\text{if  $2^{j-2}<\val{\xi}<2^{j+2},$}
\\
C_{\alpha\beta N}2^{jn} 2^{-j(N-1)}&\text{if  $\val{\xi}\le 2^{j-2},$}
\endcases
$$
implying  that  $\sum_{j}\varphi_{j}\sharp \widetilde{m_{j}}\sharp \Phi_{j}$
belongs to
$
S(\langle\xi\rangle^{-\io},\val{dx}^2+\langle\xi\rangle^{-2}\val{d\xi}^2).
$\qed
\enddemo
\proclaim{Lemma 4.1.3}
With $F(t)=f(t,x,\xi)^w$,
$\Cal M$ defined in $(4.1.10$), $M_{j}$ in $(4.1.5)$
$$\multline
\frac{d}{dt}\Cal M(t)+2\re(\Cal M(t)F(t))
\\=
\sum_{j}\Lambda_{j}^{-1/2}
\varphi_{j}^w \biggl(
\dot M_{j}(t)+2\re\Bigl(M_{j}(t) \bigl(\psi_{j}f(t) \bigr)^w\Bigr)
\biggr)
\varphi_{j}^w
+\sum_{j}
2\re\bigl(\varphi_{j}^w M_{j}(t)[\varphi_{j}^w, (\psi_{j} f(t))^w]\Lambda_{j}^{-1/2}
\bigr)\\
+
S(\langle\xi\rangle^{-\io},\langle\xi\rangle\val{dx}^2+\langle\xi\rangle^{-1}\val{d\xi}^2)^w.
\endmultline
$$
\endproclaim
\demo{Proof}
We have
$$\multline
\frac{d}{dt}\Cal M(t)+2\re(\Cal M(t)F(t))=\sum_{j}
\varphi_{j}^w \dot M_{j}(t)\Lambda_{j}^{-1/2}
\varphi_{j}^w+
2\re\bigl(\varphi_{j}^w M_{j}(t)\Lambda_{j}^{-1/2}
\varphi_{j}^w F(t) \bigr)\\=
\sum_{j}
\varphi_{j}^w \dot M_{j}(t)\Lambda_{j}^{-1/2}
\varphi_{j}^w+
2\re\bigl(\varphi_{j}^w \Lambda_{j}^{-1/2}
M_{j}(t) F(t) \varphi_{j}^w\bigr)
+
2\re\bigl(\varphi_{j}^w \Lambda_{j}^{-1/2}
M_{j}(t)[\varphi_{j}^w, F(t) ]\bigr).\endmultline
\tag 4.1.17
$$
On the other hand, we have
$$
2\re\bigl(\varphi_{j}^w \Lambda_{j}^{-1/2}
M_{j}(t)F(t)\varphi_{j}^w  \bigr)
=
2\re\Bigl(\varphi_{j}^w \Lambda_{j}^{-1/2}
M_{j}(t)\bigl(\psi_{j}f(t)\bigr)^w\varphi_{j}^w  \Bigr)
+
2\re\Bigl(\varphi_{j}^w \Lambda_{j}^{-1/2}
M_{j}(t)\bigl((1-\psi_{j})f(t)\bigr)^w\varphi_{j}^w  \Bigr)
$$
and since we have also
$$\multline
2\re\bigl(\varphi_{j}^w \Lambda_{j}^{-1/2}
M_{j}(t)[\varphi_{j}^w,F(t)]  \bigr)
\\=
2\re\bigl(\varphi_{j}^w \Lambda_{j}^{-1/2}
M_{j}(t)[\varphi_{j}^w,\bigl(\psi_{j}f(t)\bigr)^w]  \bigr)
+2\re\bigl(\varphi_{j}^w \Lambda_{j}^{-1/2}
M_{j}(t)[\varphi_{j}^w,\bigl((1-\psi_{j})f(t)\bigr)^w]  \bigr),
\endmultline$$
we get the
result of the lemma
from
Lemma 4.1.2 and (4.1.17). \qed
\enddemo
\proclaim{Lemma 4.1.4}
With the above notations, we have
$$\sum_{j}
\re\bigl(\varphi_{j}^w \Lambda_{j}^{-1/2}
M_{j}(t)[\varphi_{j}^w, (\psi_{j} f(t))^w]\bigr)\in
S(\langle\xi\rangle^{-1},
\langle\xi\rangle^{}\val{dx}^2+\langle\xi\rangle^{-1}\val{d\xi}^2)^w.
\tag 4.1.18$$
\endproclaim
\demo{Proof}
The Weyl symbol of the bracket
$[\varphi_{j}^w, (\psi_{j} f(t))^w]$
is
$
\frac{1}{2i\pi}
\poi{\varphi_{j}}{\psi_{j} f(t)} +r_{j}, \quad
r_{j}\in S(\Lambda_{j}^{-1}
,\Lambda_{j}^{-1}\Gamma_{j}
)
$
where $(r_{j})$
is a confined sequence in
$S(\langle\xi\rangle^{-1},\val{dx}^2+\langle\xi\rangle^{-2}\val{d\xi}^2)$.
As a consequence, we have
$$
\sum_{j}
\varphi_{j}^w \Lambda_{j}^{-1/2}
M_{j}(t)r_{j}^w\in S(\langle\xi\rangle^{-1},
\langle\xi\rangle^{}\val{dx}^2+\langle\xi\rangle^{-1}\val{d\xi}^2)^w.
$$
With $\Psi_{j}=-\frac{1}{2\pi}
\poi{\varphi_{j}}{\psi_{j} f(t)} $(real-valued $\in S(1,\Lambda_{j}^{-1}\Gamma_{j}
)
)$,
we are left with
$
\sum_{j}\Lambda_{j}^{-1/2}\re(\varphi_{j}\sharp\widetilde{m}_{j}(t)\sharp i \Psi_{j})
$
which belongs to
$
S(\langle\xi\rangle^{-1},
\langle\xi\rangle^{}\val{dx}^2+\langle\xi\rangle^{-1}\val{d\xi}^2).
$
\qed
\enddemo
\definition{Definition 4.1.5}
The symplectic metric $\Upsilon$ on $\RZ$ is defined as
$$
\Upsilon_{\xi}=\langle\xi\rangle^{}\val{dx}^2+\langle\xi\rangle^{-1}\val{d\xi}^2.
\tag 4.1.19$$
With $D_{j}$ given in lemma 4.1.1, we define
$$
d(t,x,\xi)=\sum_{j}\varphi_{j}(\xi)^2D_{j}(t,x,\xi).
\tag 4.1.20$$
 \enddefinition
 \proclaim{Lemma 4.1.6}
 The function $d(t,\cdot)$ is uniformly Lipschitz continuous for the metric
 $\Upsilon$ in the strongest sense, namely, there exists a positive fixed constant $C$ such that
 $$
 C^{-1}\val{d(t,x,\xi)-d(t,y,\eta)}
 \le \min\bigl(\langle\xi\rangle^{1/2},\langle\eta\rangle^{1/2}\bigr)\val{x-y}
 +\frac{\val{\xi-\eta}}{\max\bigl(\langle\xi\rangle^{1/2},\langle\eta\rangle^{1/2}\bigr)}.
 \tag 4.1.21$$
Moreover it  satisfies
$
d(t,x,\xi)\in[1,2\langle\xi\rangle^{1/2}].
$
It is thus a weight for that metric $\Upsilon$.
\endproclaim
\demo{Proof} 
 Since the $\varphi_{j}$ are nonnegative with $\sum_{j}\varphi_{j}^2=1$,
 we get from Lemma 4.1.1
 that
 $$1=\sum_{j}\varphi_{j}^2\le \sum_{j} \varphi_{j}^2D_{j}=
 d\le \sum_{j} \varphi_{j} (\xi)^2\Lambda_{j}^{1/2} 2^{1/2}
 \le\sum_{j} \varphi_{j} (\xi)^2\langle\xi\rangle^{1/2} 2=\langle\xi\rangle^{1/2} 2.
 $$
 Also, we have 
 $$
 d(t,x,\xi)-d(t,y,\eta)=
 \sum_{j}\varphi_{j}(\xi)^2\bigl( D_{j}(t,x,\xi)-D_{j}(t,y,\eta)\bigr)
 +\sum_{j}D_{j}(t,y,\eta)\bigl(\varphi_{j}(\xi)^2-\varphi_{j}(\eta)^2\bigr),
 $$
 so that, with $X=(x,\xi), Y=(y,\eta)$, $\Gamma_{j}$ given in (4.1.3),
 $$\align
 \val{ d(t,x,\xi)-d(t,y,\eta)}&\le \sum_{j}\varphi_{j}(\xi)^2 2\Gamma_{j}(X-Y)^{1/2}
 +\sum_{j,\atop\varphi_{j}(\xi)\not=0 \text{ or } \varphi_{j}(\eta)\not=0}
 2^{1/2} 2^{j/2}\val{\xi-\eta}
 2^{-j}C
 \\
 &\lesssim
 \sum_{j}\varphi_{j}(\xi)^24\bigl(\langle\xi\rangle^{1/2}\val{x-y}+
 \langle\xi\rangle^{-1/2}\val{\xi-\eta}
 \bigr)
 +
  \val{\xi-\eta}\sum_{j,\atop\varphi_{j}(\xi)\not=0 \text{ or } \varphi_{j}(\eta)\not=0} 2^{-j/2}
    \\
  &\lesssim
   \langle\xi\rangle^{1/2}\val{x-y}+
 \langle\xi\rangle^{-1/2}\val{\xi-\eta}
 +\val{\xi-\eta} \bigl(\langle\xi\rangle^{-1/2}+\langle\eta\rangle^{-1/2}\bigr).
 \endalign$$
 We get thus, if $ \langle\xi\rangle\sim  \langle\eta\rangle$,
 $$
 \val{ d(t,x,\xi)-d(t,y,\eta)}\lesssim
 \langle\xi\rangle^{1/2}\val{x-y}+
 \langle\xi\rangle^{-1/2}\val{\xi-\eta} .
 \tag 4.1.22$$
 If $2^{j_{0}}\sim \langle\xi\rangle\ll \langle\eta\rangle\sim 2^{k_{0}}$,  we have
 $$\multline
 \val{ d(t,x,\xi)-d(t,y,\eta)}\le
 \sum_{j, \varphi_{j}(\xi)\not=0} \varphi_{j}(\xi)^22^{(j+1)/2}
 +
  \sum_{j, \varphi_{j}(\eta)\not=0} \varphi_{j}(\eta)^22^{(j+1)/2}
 \\ \lesssim
  2^{j_{0}/2}+2^{k_{0}/2}
  \sim
  2^{k_{0}/2}\sim \val{\eta-\xi}2^{-k_{0}/2}\sim \langle\eta\rangle^{-1/2}\val{\eta-\xi}.
 \endmultline
 \tag 4.1.23$$
 Eventually, (4.1.23) and (4.1.22) give (4.1.21), completing the proof of the lemma.\qed
\enddemo
 Note also that $ \langle\xi\rangle$ is a
 $\Upsilon$-weight and is even such that
 $$
  \val{\langle\xi\rangle^{1/2}- \langle\eta\rangle^{1/2}}\le \frac{\val{\xi-\eta}}
  {\langle\xi\rangle^{1/2}+\langle\eta\rangle^{1/2}}.
 \tag 4.1.24$$
\proclaim{Lemma 4.1.7}
With $F(t)=f(t,x,\xi)^w$,
$\Cal M$ defined in $(4.1.10)$, $M_{j}$ in $(4.1.5)$,
the positive constant $c_{0}$
defined in lemma 4.1.1,
$$
\frac{d}{dt}\Cal M(t)+2\re(\Cal M(t)F(t))\ge 
c_{0}T^{-1} 
\sum_{j}
\varphi_{j}^w 
\bigl(\Lambda_{j}^{-1}{D_{j}^2}\bigr)^{\text{Wick($\Gamma_{j}$)}}
\varphi_{j}^w
+
S(\langle\xi\rangle^{-1},\Upsilon)^w.\tag 4.1.25
$$
The operator $\Cal M(t)$ has a Weyl symbol
in the class $
S_{1}(\langle\xi\rangle^{-1/2}d,
d^{-2}\Upsilon).
$
Moreover the selfadjoint operator $\Cal M(t)$
satisfies, with a fixed constant $C$, 
$$
\Cal M(t)\Cal M(t)\le
C^2
\sum_{j}
\varphi_{j}^w 
\bigl(\Lambda_{j}^{-1}{D_{j}^2}\bigr)^{\text{Wick($\Gamma_{j}$)}}
\varphi_{j}^w. 
\tag 4.1.26$$
\endproclaim
\demo{Proof}
The estimate
(4.1.25)
is a consequence
of the lemmas 4.1.3, 4.1.4 and 4.1.1.
From (4.1.10), we get that
$$
\Cal M(t)\in \sum_{j} \varphi_{j}^w S_{1}(D_{j}\Lambda_{j}^{-1/2},D_{j}^{-2}\Gamma_{j})^w 
\varphi_{j}^w
\subset
S_{1}(d\langle\xi\rangle^{-1/2},d^{-2}\Upsilon)^w.
$$
From the lemma 4.1.1 and the finite overlap of the $\varphi_{j}$, we get
$$\multline
\norm{\Cal M(t) u}^2
\lesssim \sum_{j}
\Lambda_{j}^{-1}\norm{\varphi_{j}^wM_{j}(t)\varphi_{j}^wu}^2
= \sum_{j}
\Lambda_{j}^{-1}\poscal{\varphi_{j}^wM_{j}\varphi_{j}^wu}{\varphi_{j}^wM_{j}\varphi_{j}^wu}
\\=
\sum_{j}
\Lambda_{j}^{-1}\poscal{\varphi_{j}^wu}{\underbrace{M_{j}(\varphi_{j}^2)^w
M_{j}}_{\in S(D_{j}^2,\Gamma_{j})^w
}\varphi_{j}^wu}
\underbrace{\lesssim }_{\text{\sevenrm from lemma A.1.4}}
\sum_{j}
\poscal{\varphi_{j}^wu}{
\bigl(\Lambda_{j}^{-1}
D_{j}^2\bigr)^{\text{Wick($\Gamma_{j}$)}}
\varphi_{j}^wu},
\endmultline
$$
which is (4.1.26).
\qed
\enddemo
\proclaim{Lemma 4.1.8}
Let $a$ be a symbol in $S(\langle\xi\rangle^{-1},\Upsilon)$.
Then, with  
constants
$C_{1}, C_{2}$ depending on a finite number of semi-norms 
of $a$, 
we have
$$
\val{\poscal{a^w u}{ u}}
\le C_{1}\norm{u}^2_{{H^{-1/2}}}
\le C_{2}
\sum_{j}
\poscal
{\bigl(\Lambda_{j}^{-1}{D_{j}^2}\bigr)^{\text{Wick($\Gamma_{j}$)}}
\varphi_{j}^w u}{\varphi_{j}^w u}.
$$
\endproclaim
\demo{Proof}
We have, since $D_{j}\ge 1$ and the Wick quantizations are nonnegative
$$
\sum_{j}
\poscal
{\bigl(\Lambda_{j}^{-1}{D_{j}^2}\bigr)^{\text{Wick($\Gamma_{j}$)}}
\varphi_{j}^w u}{\varphi_{j}^w u}
\ge 
\sum_{j}
\poscal
{\bigl(\Lambda_{j}^{-1}{}\bigr)^{\text{Wick($\Gamma_{j}$)}}
\varphi_{j}^w u}{\varphi_{j}^w u}
=\poscal{(\sum_{j}\Lambda_{j}^{-1}\varphi_{j}^2)^w u}{u}
\sim \norm{u} ^2_{{H^{-1/2}}},
$$
where $H^{-1/2}$ is the standard Sobolev space of index $-1/2$.
Now, it is a classical result that
$$
\poscal{a^w u}{ u}=
\poscal{\
\underbrace{(\langle\xi\rangle^{1/2})^w
a^w(\langle\xi\rangle^{1/2})^w}_{\in S(1,\Upsilon)^w\subset\Cal L(L^2)}(\langle\xi\rangle^{-1/2})^w u}
{(\langle\xi\rangle^{-1/2})^w u}
$$
which implies that
$
\val{\poscal{a^w u}{ u}}\lesssim \norm{u}^2_{{H^{-1/2}}}.
$
\qed
\enddemo
\proclaim{Theorem 4.1.9}
Let $f(t,x,\xi)$ be 
a smooth real-valued function 
defined on $\R\times \R^n\times \R^n$,
satisfiying
$(2.1.2)$ and $(4.1.1)$.
Let $f_{0}(t,x,\xi)$ be 
a smooth complex-valued function 
defined on $\R\times \R^n\times \R^n$,
such that $ \langle\xi \rangle f_{0}(t,x,\xi)$ satisfies $(4.1.1)$.
Then there exists $T_{0}>0, c_{0}>0$ depending on a finite number of seminorms of $f, f_{0},$
such that,
for all $T\le T_{0}$ and all $u \in C^\io_{c}\bigl((-T,T);\Cal S(\R^n)\bigr)$
$$
\norm{D_{t}u+i f(t,x,\xi)^w u +f_{0}(t,x,\xi)^w u}_{L^2(\R^{n+1})} 
\ge c_{0}T^{-1} 
\left(\int \norm{u(t)}_{H^{-1/2}(\R^n)}^2 dt\right)^{1/2}
$$
\endproclaim
\demo{Proof}
{\it \underbar{(i) We assume first that} $f_{0}\equiv 0$. }
Using the lemmas 4.1.7-8, we get 
$$
2\re\poscal{D_{t}u+i f(t)^w u}{i\Cal M(t) u}
\ge 
(c_{0}T^{-1}-C_{2})
\sum_{j}\poscal
{\varphi_{j}^w 
\bigl(\Lambda_{j}^{-1}{D_{j}^2}\bigr)^{\text{Wick($\Gamma_{j}$)}}
\varphi_{j}^wu} {u}, 
\tag 4.1.27$$
and from the estimate
(4.1.26),  provided that
$$
c_{0}/(2C_{2})\ge T,
\tag 4.1.28$$
we get 
$$
2\norm{D_{t}u+i f(t)^w u}_{{L^2(\R^{n+1})}}
\Bigl[
\sum_{j}\poscal
{\varphi_{j}^w 
\bigl(\Lambda_{j}^{-1}{D_{j}^2}\bigr)^{\text{Wick($\Gamma_{j}$)}}
\varphi_{j}^wu} {u}
\Bigr]^{1/2}C
\ge \frac{c_{0}}{2T}
\sum_{j}\poscal
{\varphi_{j}^w 
\bigl(\Lambda_{j}^{-1}{D_{j}^2}\bigr)^{\text{Wick($\Gamma_{j}$)}}
\varphi_{j}^wu} {u}
$$
so that, with fixed positive constants $c_{1}, c_{2}$,
using again the lemma 4.1.8
$$
\norm{D_{t}u+i f(t)^w u}_{{L^2(\R^{n+1})}}
\ge \frac{c_{1}}{T}
\Bigl[
\sum_{j}\poscal
{\varphi_{j}^w 
\bigl(\Lambda_{j}^{-1}{D_{j}^2}\bigr)^{\text{Wick($\Gamma_{j}$)}}
\varphi_{j}^wu} {u}
\Bigr]^{1/2}
\ge
\frac{c_{2}}{T}
\left(\int \norm{u(t)}_{H^{-1/2}(\R^n)}^2 dt\right)^{1/2},
$$
which is our result.
Let us check now the case $f_{0}\not\equiv 0$.\par
{\it \underbar{(ii) 
Let us assume that} $\im(f_{0})\in S(\langle\xi\rangle^{-1},\langle\xi\rangle^{-1}\Upsilon)$. }
Going back to the computation in (4.1.27), with (4.1.28) fulfilled, we have
$$\multline
2\re\poscal{D_{t}u+i f(t)^w +f_{0}(t)^wu}{i\Cal M(t) u}
\ge 
\frac{c_{0}}{2T}
\sum_{j}\poscal
{\varphi_{j}^w 
\bigl(\Lambda_{j}^{-1}{D_{j}^2}\bigr)^{\text{Wick($\Gamma_{j}$)}}
\varphi_{j}^wu} {u}\\
+2\re\poscal{\re(f_{0}(t))^wu}{i\Cal M(t) u}
+2\re\poscal{\im(f_{0}(t))^wu}{\Cal M(t) u}.
\endmultline
$$
From the identity
$
2\re\poscal{\re(f_{0}(t))^wu}{i\Cal M(t) u}=\poscal{\bigl[\re(f_{0}(t))^w,i\Cal M(t)\bigr]u}{ u}
$
and
the fact that, from Theorem 18.5.5 in \cite{H6} we have
$$
\bigl[\re(f_{0}(t))^w,i\Cal M(t)\bigr]\in 
S(\langle\xi\rangle^{-1/2}d d^{-1}\langle\xi\rangle^{-1/2},\Upsilon)^w
=
S(\langle\xi\rangle^{-1},\Upsilon)^w
$$
we can use the lemma 4.1.8 to control this term by
$C\sum_{j}\poscal
{\varphi_{j}^w 
\bigl(\Lambda_{j}^{-1}{D_{j}^2}\bigr)^{\text{Wick($\Gamma_{j}$)}}
\varphi_{j}^wu} {u}$.
On the other hand, from our assumption on $\im f_{0}$, we get that 
$$
\Cal M(t)\im(f_{0}(t))^w\in
S(\langle\xi\rangle^{-1/2}d \langle\xi\rangle^{-1},\Upsilon)^w\subset
S(\langle\xi\rangle^{-1},\Upsilon)^w,
$$
which can be also controlled by
$C\sum_{j}\poscal
{\varphi_{j}^w 
\bigl(\Lambda_{j}^{-1}{D_{j}^2}\bigr)^{\text{Wick($\Gamma_{j}$)}}
\varphi_{j}^wu} {u}$.
Eventually, we obtain the result in that case too, for $T$ small enough.
\par
{\it \underbar{(iii) We are left with the general case} $\im(f_{0})\in S(1,\langle\xi\rangle^{-1}\Upsilon)$}; we note that,
with
$$
\omega_{0}(t,x,\xi)=\int_{0}^t\im f_{0}(s,x,\xi) ds,\quad\text{(which belongs to
$S(1,\langle\xi\rangle^{-1}\Upsilon)$),}
\tag 4.1.29$$
we have
$$\multline
D_{t}+i f(t)^w+(\re f_{0}(t))^w+i(\im f_{0}(t))^w=(e^{\omega_{0}(t)})^wD_{t}(e^{-\omega_{0}(t)})^w
+i f(t)^w+(\re f_{0}(t))^w
\\=
(e^{\omega_{0}(t)})^w\Bigl(
D_{t}
+i f(t)^w+(\re f_{0}(t))^w
\Bigr)
(e^{-\omega_{0}(t)})^w
+ \bigl(if(t)-e^{\omega_{0}(t)}\sharp if(t)\sharp e^{-\omega_{0}(t)}\bigr)^w+
S(\langle\xi\rangle^{-1},\langle\xi\rangle^{-1}\Upsilon)^w.
\endmultline$$
Noting that $e^{\pm\omega_{0}}$  belongs to 
$S(1,\langle\xi\rangle^{-1}\Upsilon)$,
we compute
$$\align
e^{\omega_{0}}\sharp if\sharp e^{-\omega_{0}}&=
\Bigl(e^{\omega_{0}}if+\frac{1}{4i\pi}\poi{e^{\omega_{0}}}{if}\Bigr)\sharp e^{-\omega_{0}}
+S(\langle\xi\rangle^{-1},\langle\xi\rangle^{-1}\Upsilon)
\\
&=if+\frac{1}{4i\pi}
\poi{e^{\omega_{0}}if}
{e^{-\omega_{0}}}+\frac{1}{4i\pi}
\poi{e^{\omega_{0}}}
{if}e^{-\omega_{0}}
+S(\langle\xi\rangle^{-1},\langle\xi\rangle^{-1}\Upsilon)
\\
&=if+\frac{1}{2\pi}\poi{{\omega_{0}}}
{f}+S(\langle\xi\rangle^{-1},\langle\xi\rangle^{-1}\Upsilon).
\endalign
$$
We obtain
$$\multline
L=D_{t}+i f(t)^w+f_{0}(t)^w=
(e^{\omega_{0}(t)})^w\Bigl(
D_{t}
+i f(t)^w+(\re f_{0}(t)
+\frac{1}{2\pi}\poi{f}{{\omega_{0}}}
)^w
\Bigr)
(e^{-\omega_{0}(t)})^w
\\+ 
S(\langle\xi\rangle^{-1},\langle\xi\rangle^{-1}\Upsilon)^w, 
\endmultline
\tag 4.1.30$$
and analogously
$$
L_{0}=
D_{t}
+i f(t)^w+(\re f_{0}(t)
+\frac{1}{2\pi}\poi{f}{{\omega_{0}}}
)^w+ 
S(\langle\xi\rangle^{-1},\langle\xi\rangle^{-1}\Upsilon)^w
=
(e^{-\omega_{0}(t)})^w L
(e^{\omega_{0}(t)})^w.
\tag 4.1.31$$
Using now the fact
that the symbol
$\re f_{0}(t)
+\frac{1}{2\pi}\poi{f}{{\omega_{0}}}$ is real-valued in 
$S(1,\langle\xi\rangle^{-1}\Upsilon)$,
we can use (ii) to prove the estimate in the theorem
for the operator
$$L_{0}=D_{t}
+i f(t)^w+(\re f_{0}(t)
+\frac{1}{2\pi}\poi{f}{{\omega_{0}}}
)^w+S(\langle\xi\rangle^{-1},\langle\xi\rangle^{-1}\Upsilon)^w.
$$
We note also that
$
e^{\omega_{0}}\sharp e^{-\omega_{0}}=1+t^2
S(\langle\xi\rangle^{-2},\langle\xi\rangle^{-1}\Upsilon)
$
so that, for $\val{t}$ small enough,
$$%
\spreadmatrixlines {5pt}
\left.
\matrix
\text{
the operators
$(e^{\pm\omega_{0}})^w$ are invertible in $L^2(\R^n)$}
\\
\text{
and their inverses are  pseudodifferential operators in $S(1,\langle\xi\rangle^{-1}\Upsilon)^w.
$}
\endmatrix\right\}
\tag 4.1.32
$$
From the previous identity and (ii), we get for $u\in \mooc((-T,T), \Cal S(\R^n))$
$$
\int \nuorm{(e^{-\omega_{0}(t)})^wL(e^{\omega_{0}(t)})^w u(t)}_{{L^2(\R^n)}}^2 dt 
\ge \frac{c_{0}^2}{T^2}
\int \nuorm{u(t)}_{{H^{-1/2}(\R^n)}}^2 dt.
$$
Applying this to 
$$
u(t)= \Bigl((e^{\omega_{0}(t)})^w\Bigr)^{-1}v(t),
\tag 4.1.33$$
we obtain
$$
\int \nuorm{(e^{-\omega_{0}(t)})^wLv(t)}_{{L^2(\R^n)}}^2 dt 
\ge \frac{c_{0}^2}{T^2}
\int \nuorm{\Bigl((e^{\omega_{0}(t)})^w\Bigr)^{-1}v(t)}_{{H^{-1/2}(\R^n)}}^2 dt.
\tag 4.1.34$$
We have
$$
 \nuorm{\Bigl((e^{\omega_{0}(t)})^w\Bigr)^{-1}v(t)}_{{H^{-1/2}(\R^n)}}^2=
 \nuorm{(\langle\xi\rangle^{-1/2})^w\Bigl((e^{\omega_{0}(t)})^w\Bigr)^{-1}
 (\langle\xi\rangle^{1/2})^w
 (\langle\xi\rangle^{-1/2})^w
 v(t)}_{L^2(\R^n)}^2.
$$
Now the operator
$(\langle\xi\rangle^{-1/2})^w\Bigl((e^{\omega_{0}(t)})^w\Bigr)^{-1}
 (\langle\xi\rangle^{1/2})^w$
 is invertible
 with inverse
 $$\Omega(t)=(\langle\xi\rangle^{-1/2})^w(e^{\omega_{0}(t)})^w
 (\langle\xi\rangle^{1/2})^w
 \tag 4.1.35$$
 which is a bounded operator on $L^2(\R^n)$ so that
 $$
 \norm{v}_{L^2}=\norm{\Omega^{}\Omega^{-1} v}_{{L^2}}\le \norm{\Omega}_{\Cal L(L^2)}
 \norm{\Omega^{-1} v}_{L^2}.
 \tag 4.1.36$$
 As a result, from the inequality (4.1.34), we get
 $$\multline
\int \nuorm{(e^{-\omega_{0}(t)})^wLv(t)}_{{L^2(\R^n)}}^2 dt 
\ge \frac{c_{0}^2}{T^2}
\int \nuorm{\Omega(t)^{-1}(\langle\xi\rangle^{-1/2})^wv(t)}_{{L^2(\R^n)}}^2 dt
\\
\ge 
\frac{c_{0}^2}{T^2}\int\nuorm{(\langle\xi\rangle^{-1/2})^wv(t)}_{{L^2(\R^n)}}^2 \frac{1}{\norm{\Omega(t)}^2}dt
\ge 
\frac{c_{1}^2}{T^2}\int\nuorm{v(t)}_{{H^{-1/2}(\R^n)}}^2 dt,
\endmultline$$
which is the result.
The proof of Theorem 4.1.9 is complete.
\qed
\enddemo
\remark{Comment 4.1.10}{
Although Theorem 4.1.9 is providing
a solvability
result
with loss
of 3/2 derivatives for the evolution equation
$$
\partial_{t}+f(t,x,\xi)^w+ f_{0}(t,x,\xi)^w,
$$
where $f, f_{0}$ are satisfying the assumptions of this theorem, the statement
does not seem  quite sufficient
to handle operators with homogeneous symbols
for two reasons. The first one is that
the reduction of homogeneous symbols in the cotangent bundle
of a manifold
 will lead to a model operator like the one above, but only at the cost of some microlocalization
 in the cotangent bundle.
 We need thus to get a microlocal version of  our estimates.
 The second reason is that the 
 function
$f(t,x,\xi)$ is not a classical symbol in
 the phase space $\R_{t}\times \R^n_{x}\times \R_{\tau}\times \R^n_{\xi}$
 and we have to pay attention
 to the discrepancy between
 homogeneous localization in the phase space $\R^{2n+2}$
 and localization in $\RZ$ with parameter $t$.
 That difficulty should be taken seriously,
 since the loss of derivatives is strictly larger than 1;  in fact,
 commuting a cutoff function with the operator
 will produce an error of order 0, larger than what is controlled by the estimate.
 In the next section, we prove a localized version
 of the theorem 4.1.9, which will be suitable for future use
 in the
 homogeneous
 framework.
 }
 \endremark
 \subhead
4.2.  From  semi-classical to localized inhomogeneous estimates
\endsubhead
We begin with a modified version of Lemma 4.1.7, involving a microlocalization in $\RZ$.
\proclaim{Lemma 4.2.1}
Let $f(t,x,\xi)$ be real-valued satisfying $(2.1.2)$ and $(4.1.1)$; we shall note
$F(t)=f(t,x,\xi)^w$. Let 
$\Cal M$ be defined in $(4.1.10)$.
We define  $c_{1}=c_{0}/ C^2$, where $c_{0}$ is given by lemma 4.1.1
and $C$ appears in $(4.1.26)$.
Let $\psi(x,\xi)$ be a real-valued
symbol in $S(1,\langle\xi\rangle^{-1}\Upsilon)$. We have
$$
\frac{d}{dt}\bigl(\psi^w\Cal M(t)\psi^w\bigr)+2\re
\bigl(\psi^w  \Cal M(t)\psi^wF(t)\bigr)\hskip-3pt 
\ge c_{1} T^{-1}
\psi^w\Cal M(t)\Cal M(t)
\psi^w+ S(\langle \xi\rangle^{-1},\Upsilon)^w.\tag 4.2.1
$$
\endproclaim
\demo{Proof}
We compute, using (4.1.25) on the fourth line below,
$$\align
\frac{d}{dt}&\bigl(\psi^w\Cal M(t)\psi^w\bigr)+2\re
\bigl(\psi^w  \Cal M(t)\psi^wF(t)\bigr)
\\&=
\psi^w\dot{\Cal M}(t)\psi^w+
\psi^w  \Cal M(t)\psi^wF(t)+F(t)\psi^w  \Cal M(t)\psi^w
\\
&=
\psi^w\Bigl(\dot{\Cal M}(t)+ 2\re
\Cal M(t) F(t)
\Bigr)\psi^w+
\psi^w  \Cal M(t)\bigl[\psi^w,F(t)\bigr]+\bigl[F(t),\psi^w \bigr] \Cal M(t)\psi^w
\\
&\ge c_{1} T^{-1}
\psi^w\Cal M(t)\Cal M(t)
\psi^w
\\&\hskip25pt +
\psi^w  \Bigl[\Cal M(t), \bigl[\psi^w,F(t)\bigr]\Bigr]+
\psi^w\bigl[\psi^w,F(t)\bigr]
 \Cal M(t)-
 \bigl[\psi^w, F(t) \bigr] \Cal M(t)\psi^w
 \\
&= c_{1} T^{-1}
\psi^w\Cal M(t)\Cal M(t)
\psi^w
\\&\hskip25pt +
\psi^w  \Bigl[\Cal M(t), \bigl[\psi^w,F(t)\bigr]\Bigr]+
\Bigl[\psi^w,\bigl[\psi^w,F(t)\bigr]\Bigr]
 \Cal M(t)
+\bigl[\psi^w,F(t)\bigr]\bigl[\psi^w, \Cal M(t)\bigr].
\endalign
$$
Next we analyze each term on the last line. We have
\roster
\item"{$\bullet$}"$\psi^w  \Bigl[\Cal M(t), \bigl[\psi^w,F(t)\bigr]\Bigr]
\in S(d\langle \xi\rangle^{-1/2} 1d^{-1}\langle \xi\rangle^{-1/2} ,\Upsilon)^w=
S(\langle \xi\rangle^{-1},\Upsilon)^w$
since
$$\psi^w,\ [\psi^w,F(t)\bigr]\in S(1,\langle \xi\rangle^{-1} \Upsilon)^w,\quad
\Cal M(t)\in S_{1}(d\langle \xi\rangle^{-1/2},d^{-2} \Upsilon)^w,$$
\item"{$\bullet$}"
$\Bigl[\psi^w,\bigl[\psi^w,F(t)\bigr]\Bigr]
 \Cal M(t)
 \in S(d\langle \xi\rangle^{-3/2} ,\Upsilon)^w\subset
S(\langle \xi\rangle^{-1},\Upsilon)^w$
since $d\le 2 \langle \xi\rangle^{1/2}$ and
$$\Bigl[\psi^w,\bigl[\psi^w,F(t)\bigr]\Bigr]
\in S(\langle \xi\rangle^{-1} ,\langle \xi\rangle^{-1} \Upsilon)^w,\quad
\Cal M(t)\in S_{1}(d\langle \xi\rangle^{-1/2},d^{-2} \Upsilon)^w,$$
\item"{$\bullet$}"
$\bigl[\psi^w,F(t)\bigr]\bigl[\psi^w, \Cal M(t)\bigr],
 \in S( d\langle \xi\rangle^{-1/2} \langle \xi\rangle^{-1/2}d^{-1},\Upsilon)^w=
S(\langle \xi\rangle^{-1},\Upsilon)^w$
since 
$$\bigl[\psi^w,F(t)\bigr]
\in S(1,\langle \xi\rangle^{-1} \Upsilon)^w,\quad
\Cal M(t)\in S_{1}(d\langle \xi\rangle^{-1/2},d^{-2} \Upsilon)^w.$$
\endroster
We have proven in particular that
$$
\frac{d}{dt}\bigl(\psi^w\Cal M(t)\psi^w\bigr)+2\re
\bigl(\psi^w  \Cal M(t)\psi^wF(t)\bigr)
=
\psi^w\Bigl(\dot{\Cal M}(t)+ 2\re
\Cal M(t) F(t)
\Bigr)\psi^w+S(\langle \xi\rangle^{-1},\Upsilon)^w.
\tag 4.2.2$$
Also, we have
$
\frac{d}{dt}\bigl(\psi^w\Cal M(t)\psi^w\bigr)+2\re
\bigl(\psi^w  \Cal M(t)\psi^wF(t)\bigr)
\ge c_{1} T^{-1}
\psi^w\Cal M(t)\Cal M(t)
\psi^w+ S(\langle \xi\rangle^{-1},\Upsilon)^w,
$
which is (4.2.1).\qed
\enddemo
\proclaim{Theorem 4.2.2}
Let $f(t,x,\xi)$ be 
a smooth real-valued function 
defined on $\R\times \R^n\times \R^n$,
satisfiying
$(2.1.2)$ and $(4.1.1)$.
Let $f_{0}(t,x,\xi)$ be 
a smooth complex-valued function 
defined on $\R\times \R^n\times \R^n$,
such that $ \langle\xi \rangle f_{0}(t,x,\xi)$ satisfies $(4.1.1)$.
We define
$$
L=D_{t}+i f(t,x,\xi)^w +f_{0}(t,x,\xi)^w.
$$
Let   $\psi(x,\xi)\in S(1,\langle \xi\rangle^{-1}\Upsilon)$ be a real-valued symbol.
Then there exists $T_{0}>0, c_{0}>0, C\ge 0,$ depending on a finite number of seminorms of $f, f_{0},\psi, $
such that,
for all $T\le T_{0}$, 
all $u \in C^\io_{c}\bigl((-T,T);\Cal S(\R^n)\bigr)$,
with $\omega_{0}$
given by $(4.1.29)$,
$$\multline
T\norm{\psi^w (e^{-\omega_{0}})^w Lu}_{L^2(\R^{n+1})} +C T^{1/2}
\left(\int \norm{u(t)}_{H^{-1/2}(\R^n)}^2 dt\right)^{1/2}
+C
\left(\int \norm{u(t)}_{H^{-3/2}(\R^n)}^2 dt\right)^{1/2}
\\
\ge c_{0}
\left(\int \norm{\psi^wu(t)}_{H^{-1/2}(\R^n)}^2 dt\right)^{1/2}.
\endmultline
\tag 4.2.3$$
\endproclaim
\demo{Proof}
We compute, noting  $F(t)= f(t,x,\xi)^w$,
$$\multline
2\re\poscal{Lu}{i\psi^w \Cal M(t)\psi^w u}
=
\left\langle{\Bigl(\psi^w\dot{\Cal M}(t)\psi^w+2\re
\bigl(\psi^w  \Cal M(t)\psi^wF(t)\bigr)\Bigr)u
},{u}\right\rangle
\\+
\left\langle{
\Bigl[\bigl(\re f_{0}(t)\bigr)^w ,    i\psi^w  \Cal M(t)\psi^w\Bigr]
u
},{u}\right\rangle
+2\re 
\left\langle{
\psi^w  \Cal M(t)\psi^w
\im f_{0}(t)^w   u
},{u}\right\rangle.
\endmultline$$
{\it \underbar{(i) 
Let us assume that} $\im(f_{0})\in S(\langle\xi\rangle^{-1},\langle\xi\rangle^{-1}\Upsilon)$. }
Then we get that
$$
\psi^w  \Cal M(t)\psi^w
\im f_{0}(t)^w  \in S(d\langle\xi\rangle^{-1/2}\langle\xi\rangle^{-1},\Upsilon)^w
\subset S(\langle\xi\rangle^{-1},\Upsilon)^w
$$
and since
$
\Bigl[\bigl(\re f_{0}(t)\bigr)^w ,    i\psi^w  \Cal M(t)\psi^w\Bigr]
 \in S(d\langle\xi\rangle^{-1/2}\langle\xi\rangle^{-1/2}d^{-1},\Upsilon)^w
=S(\langle\xi\rangle^{-1},\Upsilon)^w,
$
the inequality
(4.1.25) , the identity (4.2.2) and  lemmas 4.1.8 -- 4.2.1
show that
$$\multline
2\re\poscal{Lu}{i\psi^w \Cal M(t)\psi^w u}
=
\left\langle{\Bigl(\psi^w\dot{\Cal M}(t)\psi^w+2\re
\bigl(\psi^w  \Cal M(t)\psi^wF(t)\bigr)\Bigr)u
},{u}\right\rangle
\\
\ge
\frac{c_{1}}{2}T^{-1}\int\norm{\Cal M(t) \psi^w u(t)}^2_{{L^2(\R^n)}} dt 
+\frac{c_{0}}{2}T^{-1}\int\norm{\psi^w u(t)}^2_{{H^{-1/2}(\R^n)}} dt 
\\
-C\int\norm{ u(t)}^2_{{H^{-1/2}(\R^n)}} dt.
\endmultline
$$
As a consequence, we have
$$\multline
2T\int \norm{\psi^w Lu(t)}_{L^2(\R^n)}\norm{\Cal M(t)\psi^w u(t)}_{L^2(\R^n)}dt
+CT\int\norm{ u(t)}^2_{{H^{-1/2}(\R^n)}} dt
\\
\ge
\frac{c_{1}}{2}\int\norm{\Cal M(t) \psi^w u(t)}^2_{{L^2(\R^n)}} dt 
+\frac{c_{0}}{2}\int\norm{\psi^w u(t)}^2_{{H^{-1/2}(\R^n)}} dt ,
\endmultline
$$
so that, with $\alpha>0$,
$$\multline
T\int\Bigl( T\alpha^{-1} \norm{\psi^w Lu(t)}_{L^2(\R^n)}^2+\alpha T^{-1}\norm{\Cal M(t)\psi^w u(t)}_{L^2(\R^n)}^2\Bigr)dt
+CT\int\norm{ u(t)}^2_{{H^{-1/2}(\R^n)}} dt
\\
\ge
\frac{c_{1}}{2}\int\norm{\Cal M(t) \psi^w u(t)}^2_{{L^2(\R^n)}} dt 
+\frac{c_{0}}{2}\int\norm{\psi^w u(t)}^2_{{H^{-1/2}(\R^n)}} dt .
\endmultline
$$
Choosing $\alpha\le c_{1}/2$ yields the  
result
$$
T^2\alpha^{-1} \int\norm{\psi^w Lu(t)}_{L^2(\R^n)}^2dt
+CT\int\norm{ u(t)}^2_{{H^{-1/2}(\R^n)}} dt
\ge
\frac{c_{0}}{2}\int\norm{\psi^w u(t)}^2_{{H^{-1/2}(\R^n)}} dt ,
$$
which is a better estimate than the sought one.\par\no
{\it \underbar{(ii) 
Let us deal now with the general case} $\im(f_{0})\in S(1,\langle\xi\rangle^{-1}\Upsilon)$. }
Using the definitions (4.1.29), (4.1.31) and the property (4.1.30),
we can use (i) above to get the estimate for $L_{0}$, so that
with a fixed $c_{2}>0$
$$
T\norm{\psi^w L_{0}u}_{L^2(\R^{n+1})} + T^{1/2}
\left(\int \norm{u(t)}_{H^{-1/2}(\R^n)}^2 dt\right)^{1/2}
\ge c_{2}
\left(\int \norm{\psi^wu(t)}_{H^{-1/2}(\R^n)}^2 dt\right)^{1/2},
\tag 4.2.4$$
so that
$$
T\norm{\psi^w
(e^{-\omega_{0} })^w
L(e^{\omega_{0} })^wu}_{L^2(\R^{n+1})} + T^{1/2}
\left(\int \norm{u(t)}_{H^{-1/2}(\R^n)}^2 dt\right)^{1/2}
\ge c_{2}
\left(\int \norm{\psi^wu(t)}_{H^{-1/2}(\R^n)}^2 dt\right)^{1/2}.
\tag 4.2.5$$
Applying this to 
$u(t)$ given by (4.1.33), we obtain
$$\multline
T\norm{\psi^w
(e^{-\omega_{0} })^w
Lv}_{L^2(\R^{n+1})} + T^{1/2}
\left(\int \nuorm{\Bigl((e^{\omega_{0} })^w\Bigr)^{-1}v(t)}_{H^{-1/2}(\R^n)}^2 dt\right)^{1/2}
\\
\ge c_{2}
\left(\int \nuorm{\psi^w\Bigl((e^{\omega_{0} })^w\Bigr)^{-1}v(t)}_{H^{-1/2}(\R^n)}^2 dt\right)^{1/2}.
\endmultline
\tag 4.2.6$$
Using that
$\Bigl((e^{\omega_{0} })^w\Bigr)^{-1}$ is a
pseudodifferential operator with symbol in 
$S(1,\langle \xi\rangle^{-1}\Upsilon)$,
we obtain, using the notation (4.1.35),
$$\multline
T\norm{\psi^w
(e^{-\omega_{0} })^w
Lv}_{L^2(\R^{n+1})} +C T^{1/2}
\left(\int \norm{v(t)}_{H^{-1/2}(\R^n)}^2 dt\right)^{1/2}
\\
\ge c_{2}
\left(\int \norm{
\Omega(t)^{-1}
(\langle\xi\rangle^{-1/2})^w
\psi^wv(t)}_{L^2(\R^n)}^2 dt\right)^{1/2}
-C_{1}
\left(\int \norm{v(t)}_{H^{-3/2}(\R^n)}^2 dt\right)^{1/2},
\endmultline
\tag 4.2.7$$
so that, using (4.1.36),
$$\align
T&\norm{\psi^w
(e^{-\omega_{0} })^w
Lv}_{L^2(\R^{n+1})}
 +C T^{1/2}
\left(\int \norm{v(t)}_{H^{-1/2}(\R^n)}^2 dt\right)^{1/2}
 +
C_{1}
\left(\int \norm{v(t)}_{H^{-3/2}(\R^n)}^2 dt\right)^{1/2}
\\
&\ge c_{2}
\left(\int \nuorm{
(\langle\xi\rangle^{-1/2})^w
\psi^wv(t)}_{L^2(\R^n)}^2 \frac{1}{	\norm{\Omega(t)}^2}dt\right)^{1/2}
\\
&\ge c_{3}
\left(\int \nuorm{
(\langle\xi\rangle^{-1/2})^w
\psi^wv(t)}_{L^2(\R^n)}^2dt\right)^{1/2}
=c_{3}
\left(\int \nuorm{
\psi^wv(t)}_{H^{-1/2}(\R^n)}^2dt\right)^{1/2},\tag 4.2.8
\endalign
$$
which is the result.
The proof of the theorem is complete.
\qed
\enddemo
\subhead
4.3.  From inhomogeneous localization to homogeneous localization
\endsubhead
In this section, we are given a positive integer $n$, and we define $N=n+1$.
The running point of $T^\ast(\R^N)$ will be denoted by $(y,\eta).$
We are also given a
point
$(y_0;\eta_0)\in \R^{N}\times \Bbb S^{N-1}$
such that
$$
Y_0=(y_0;\eta_0)=(t_0,x_0; \tau_0,\xi_0)\in \R\times\R^n\times \R\times\R^n,\quad\text{with $\tau_0=0,\ \xi_0\in \Bbb S^{n-1}, t_{0}=0$.}
\tag 4.3.1$$
We consider
$F(t,x,\xi)=f(t,x,\xi)-if_0(t,x,\xi)$,
with $f,f_0$ satisfying the assumptions of Theorem 4.2.2.
Let $\psi_0 (\xi)$
be a function supported in a conic neighborhood
of $\xi_0$ and
$\chi_0(\tau,\xi)$ be  an homogeneous localization 
near $\tau=0$ as in the appendix A.7 with some positive $r_{0}$.
We consider also a classical first-order pseudodifferential operator $R$ in $\R^N$
such that
$
Y_{0}\notin WF R.
$
We consider the first-order operator
$$
\Cal L= D_t+ i\bigl(F(t,x,\xi)\psi_0(\xi)\chi_0(\tau,\xi)\bigr)^w+R.
\tag 4.3.2$$
We have
$$
\Cal L= D_t+ i\bigl(F(t,x,\xi)\psi_0(\xi)\bigr)^w
+{\underbrace{i
\Bigl(F(t,x,\xi)\psi_0(\xi)\bigl(\chi_0(\tau,\xi)-1\bigr)\Bigr)^w}_{=F_1(t,x,\tau,\xi)^w}}+
R.\tag 4.3.3$$
Let $\psi_1 (\xi)$
be a function supported in a conic neighborhood
of $\xi_0$ and
$\chi_1(\tau,\xi)$ be  an homogeneous localization 
near $\tau=0$ as in the appendix A.7 with some positive $r_{1}<r_{0}$
and such that
$$\gather
\supp \chi_{1}\subset\{\chi_{0}=1\}, \quad\supp( \psi_{1} \chi_{1})\subset
\{\psi_{0}\chi_{0}=1\},\tag 4.3.4
\\
[-T_{1},T_{1}]\times K_{1}\times\supp \psi_{1}\chi_{1}\subset (WF R)^c,
\tag 4.3.5
\endgather$$
where $T_{1}>0$ and $K_{1}$ is a compact neighborhood of $x_{0}$.
Let $\psi(x,\xi)$
be a symbol satisfying the assumptions of Theorem 4.2.2
and let
 $\rho_{1}\in \mooc(\R)$,
such that
$$
\supp \psi\subset K_{1}\times\{\psi_{1}=1\},\quad \supp\rho_{1}\subset[-T_{1},T_{1}].
\tag 4.3.6$$
We can apply the theorem 4.2.2 to the operator
$L=D_t+ i\bigl(F(t,x,\xi)\psi_0(\xi)\bigr)^w$.
We have,
with
$u\in\Cal S(\R^N)$, 
$$\multline
T_{1}\norm{\psi^w (e^{-\omega_{0}})^w (\Cal L-F_{1}-R)\rho_{1}\chi_{1}^wu}_{L^2(\R^{n+1})} +C T_{1}^{1/2}
\left(\int \norm{\rho_{1}\chi_{1}^wu(t)}_{H^{-1/2}(\R^n)}^2 dt\right)^{1/2}
\\+C
\left(\int \norm{\rho_{1}\chi_{1}^wu(t)}_{H^{-3/2}(\R^n)}^2 dt\right)^{1/2}
\ge c_{0}
\left(\int \norm{\psi^w\rho_{1}\chi_{1}^wu(t)}_{H^{-1/2}(\R^n)}^2 dt\right)^{1/2}.
\endmultline
$$
We get then
$$\aligned
T_{1}&\norm{\psi^w (e^{-\omega_{0}})^w \rho_{1} \chi_1^w \Cal L  u
+\psi^w (e^{-\omega_{0}})^w \rho_{1}  [\Cal L, \chi_1^w ]u
}_{L^2(\R^{n+1})}
\\+T_{1}&\norm{\psi^w (e^{-\omega_{0}})^w [\Cal L, \rho_{1}] \chi_1^w u}_{L^2(\R^{n+1})}
\\+T_{1}&\norm{\psi^w (e^{-\omega_{0}})^w F_1^w \rho_{1}\chi_1^w u}_{L^2(\R^{n+1})} 
+T_{1}\norm{\psi^w (e^{-\omega_{0}})^w R \rho_{1} \chi_1^w u}_{L^2(\R^{n+1})} 
\\+C &T_{1}^{1/2}
\nuorm{(\langle\xi\rangle^{-1/2})^w\rho_{1} \chi_1^w u}_{L^2(\R^{n+1})}
+C
\nuorm{(\langle\xi\rangle^{-3/2})^w\rho_{1} \chi_1^w u}_{L^2(\R^{n+1})}
\\
&\ge c_{0}
\left(\int \norm{\psi^w\rho_{1} \chi_1^w u}_{H^{-1/2}(\R^n)}^2 dt\right)^{1/2}.
\endaligned
\tag 4.3.7$$
We assume now that $u\in \Cal S(\R^N)$, $\supp u\subset\{(t,x),\val t\le T_{1}/2\}$
and also that
$\rho_{1}$ is 1 on $[-3T_{1}/4,3T_{1}/4]$.
We introduce two admissible
\footnote{The properties of definition 1.3.1 are classical for $G$ and easily
checked for $g$.  One can check also that
$(1+\val \xi+\val \tau)^s$ are $G$-weights and $(1+\val \xi)^s$
are $g$-weights.}
  metrics on $\R^{2N}$,
$$
G=\val{dt}^2+\val{dx}^2+\frac{\val{d\xi}^2+\val{d\tau}^2}{1+\val\xi^2+\tau^2}
\le
g=\val{dt}^2+\val{dx}^2+\frac{\val{d\xi}^2}{1+\val\xi^2}
+\frac{\val{d\tau}^2}{1+\val\xi^2+\tau^2}.
\tag 4.3.8$$
\item{(1)} The operator $[\Cal L,\chi_1^w]$ has a symbol in $S(1, G)$ which is essentially supported
in the  region where
$\val\tau\sim \val \xi$.
\item{(2)} The quantity
$[\Cal L,\rho_{1}]\chi_1^w u=[\Cal L,\rho_{1}]\chi_1^w \rho_2u$
if $\rho_2(t)$ is 1 on $[-T_{1}/2,T_{1}/2]$ and supported in $[-3T_{1}/4,3T_{1}/4]$
and thus the operator
$[\Cal L,\rho_{1}]\chi_1^w \rho_2$
has a symbol in $S((1+\val\xi+\val\tau)^{-\io}, G)$.
\item{(3)}
The operator $F_1^w \rho_{1}\chi_1^w$ is the composition
of the symbol $F_{1}\in S(\langle\xi\rangle,g)$
with the  symbol in $\rho_{1} \sharp \chi_{1} \in S(1,G)$ and thus is a priori in 
$S(\langle\xi\rangle,g)$; however, looking at the expansion,
and using (4.3.4), 
we see that it has a symbol in  $S((1+\val\xi+\val\tau)^{-\io}, G)$: it is not completely
 obvious though and we refer the reader to the lemma A.8.1 for a complete argument.
 \item{(4)}
The operator $\psi^w (e^{-\omega_{0}})^w R \rho_{1} \chi_1^w $
is also the composition of an operator in $S(1, g)^w$
with an operator in $S(\langle \xi,\tau\rangle, G)^w$ ; however,
using (4.3.4-5-6)
and the appendix A.8, we see 
that 
$\psi^w (e^{-\omega_{0}})^w R\rho_{1}\chi_{1}^w$
has a symbol in 
$S((1+\val\xi+\val\tau)^{-\io}, G)$. 
 \item{(5)} The operator
 $(\langle\xi\rangle^s)^w \rho_{1} \chi_1^w$
 is also the sum of an operator
 in $S(\langle\tau,\xi\rangle^s,G)$
 plus a symbol  in $S((1+\val\xi+\val\tau)^{-\io}, G)$.
\vs
We write now, with
$R_{1}$ of order $-\io$ (weight $\langle \xi,\tau\rangle^{-\io}$) for $G$,
$E_{0}$ of order 0 (weight 1) for $G$, supported in 
$\{(t,x,\tau,\xi), \val t\le T_{1}, x\in K_{1}, (\tau,\xi)\in \supp \nabla \chi_{1}, (x,\xi)\in \supp \psi\}$,
$$\multline
T_{1}\norm{\psi^w (e^{-\omega_{0}})^w \rho_{1} \chi_1^w \Cal L  u
+E_{0}u
}_{L^2(\R^{n+1})}
+T_{1}\norm{R_{1}u}_{L^2(\R^{n+1})}
+ CT_{1}^{1/2}
\norm{ u}_{H^{-1/2}(\R^{n+1})}
+C
\norm{ u}_{H^{-3/2}(\R^{n+1})}
\\
\ge c_{0}
\left(\int \norm{\psi^w\rho_{1} \chi_1^w u}_{H^{-1/2}(\R^n)}^2 dt\right)^{1/2}.
\endmultline
\tag 4.3.9$$
\proclaim{Theorem 4.3.1}
Let $\Cal L$ be the pseudodifferential operator given by $(4.3.2)$
and $Y_0=(y_0,\eta_0)$ be given by $(4.3.1)$.
We assume  that
$
\{Y_{0}\}\subset \Delta_{0}\subset (WF R)^c,
$
where $\Delta_{0}$ is a compact-conic neighborhood
of $Y_{0}$.
Then,
there exists two pseudodifferential operators $\Phi_0, \Psi_{0}$
of order $0$
(weight $1$) for $G$,
both essentially supported
in $\Delta_{0}$
with
$\Phi_0$ is elliptic at $Y_0$, 
and there exists $r>0$
such that, for all
$u\in \Cal S(\R^N), \supp u\subset \{(t,x),\val{t}\le r\},$
$$
r\norm{ \Psi_{0}\Cal L  u}_{L^2(\R^{N})}
+ r^{1/2}
\norm{ u}_{H^{-1/2}(\R^{N})}
+
\norm{ u}_{H^{-3/2}(\R^{N})}
\ge 
\norm{\Phi_0 u}_{H^{-1/2}(\R^N)}.
\tag 4.3.10$$
\endproclaim
\demo{Proof}
It is a direct consequence of (4.3.9)
since,
using the ellipticity of
$\Cal L$ in the support of the symbol of $E_{0}$, we get 
$E_0= \Cal K \Cal L+R_{2}$,
where $\Cal K$ is a pseudodifferential operator of order 0
such that
$
WF \Cal K\subset \Delta_{0}
$
 and 
$R_{2}$
is a pseudodifferential operator of order $-\io$ for $G$.\qed
\enddemo
\subhead
4.4.  Proof of the solvability result stated in Theorem 1.2.2
\endsubhead
Let $P$ be a first-order pseudodifferential operator with principal symbol $p$
satisfying the assumptions of Theorem 1.2.2 and let
$(y_0,\eta_0)$ be a point in the cosphere bundle.
If $p(y_0,\eta_0)\not=0$, then there exists a pseudodifferential operator $\Phi_0$
of order 0, elliptic at 
$(y_0,\eta_0)$
such that
$$
\norm{P^* u}_{0}+\norm{u}_{-1}\ge \norm{\Phi_0 u}_{1}.
\tag 4.4.1$$
In fact, the ellipticity assumption implies that
there
exist a pseudodifferential operator $K$
of order $-1$ and a pseudodifferential operator $R$
of order $0$
such that
$$
\Id= KP^*+R,\quad (y_0,\eta_0)\notin WF R.
$$
As  consequence, 
for $\Phi_0$ of order $0$ essentially supported
close enough to  $(y_0,\eta_0)$, we get
$
\Phi_0=\Phi_0 KP^* +\Phi_0 R
$
with 
$\Phi_0 R$ of order $-\io$, which gives (4.4.1).
\par
Let us assume now that $p(y_0,\eta_0)=0$.
We know from the assumption (1.2.1) that
$\p_\eta p(y_0,\eta_0)\not=0$ and we may suppose that
$
(\p_\eta \re p)(y_0,\eta_0)\not=0.
$
Using the Malgrange-Weierstrass theorem, we can find a conic
neighborhood
of
$(y_0,\eta_0)$ in which 
$$
p(y,\eta)= \bigl(\sigma + a(s,z,\zeta)+ib(s,z,\zeta)\bigr)e_0(y,\eta)
$$
where
$a,b$ are real-valued positively homogeneous of degree 1,
$e_0$ is homogeneous of degree 0, elliptic near $(y_0,\eta_0)$,
$(s,z;\sigma,\zeta)\in\R\times \R^n\times\R\times \R^n$
a choice of symplectic coordinates in $T^*(\R^N)$ ($N=n+1$),
with 
$y_0=(0,0), \eta_0=(0,\dots,0,1)$.
Noting that the Poisson bracket
$$
\poi{\sigma+a}{s}=1
$$
we see that there exists an homogeneous canonical transformation
$\Xi^{-1}$,
from a (conic)
neighborhood
of $(y_0,\eta_0)$ to a
conic neighborhood
of $(0;0,\dots 0,1) $ in $\R^N\times \R^N$
such that
$$
p\circ \Xi=\bigl(\tau +i q(t,x,\xi)\bigr) (e\circ\Xi).
$$
Note in particular that, setting $\tau=\sigma+a, t=s,$
{\it (which preserves the coordinate $s$)}
 yields 
$$-\p_\tau q=\poi{t}{q}=\poi{s}{b}\circ \chi=0.
$$
We see now that there exists some
elliptic Fourier integral operators $A,B$ and $E$
a  pseudodifferential 
operator of order 0, elliptic at $(y_{0},\eta_{0})$
such that
$$\align
AEP^*B&= D_t+i(f(t,x,\xi)\chi_0(\tau,\xi))^w +R,\\
\quad BA&=\Id+S, (y_0,\eta_0)\in \Gamma_{0}(\text{\sevenrm conic neighborhood of
$\scriptstyle (y_0,\eta_0)$})\subset (WF S)^c,
\endalign$$
where
$f$ satisfies (2.1.2),
$R$ is a pseudodifferential operator
of order 0,
and $\chi_0$ 
is a nonnegative
homogeneous localization near $\tau =0$.
Using the fact that the coordinate $s$ is preserved by the canonical transformation,
we can  assume
that $A, B$ are local operators in the $t$ variable, i.e.
are such that
$$
u\in\mooc, \supp u\subset \{(t, x)\in \R\times \R^n,\val t\le r\}\Longrightarrow
\supp Bu\subset \{(s, z)\in \R\times \R^n,\val s\le r\}.
$$
Using the fact that
the operator
$P$ is polyhomogeneous,
one can
iterate the use of the Malgrange-Weierstrass theorem
to reduce our case to $AEP^* B=\Cal L$ of the type given in (4.3.2).
We can apply the theorem 4.3.1, giving the existence
of a pseudodifferential operator $\Psi_0$
of order 0, elliptic at  $\Xi^{-1}(y_0,\eta_0)$,
essentially supported in 
$\Xi^{-1}(\Gamma_{0})$
such that for all
$u\in\mooc(\R^N), \supp u\subset\{\val t\le r\}$,
$$
r\norm{\Psi_{0} AEP^* B  u}_{0}
+ r^{1/2}
\norm{ u}_{-1/2}
+
\norm{ u}_{-3/2}
\ge 
\norm{\Phi_0 u}_{-1/2}.
$$
We may assume that $A$ and $B$ are properly supported and apply the previous inequality to 
$u=Av$,
whose support in the $s$ variable is unchanged.
We get
$$
r\norm{\Psi_{0} AEP^* B  Av}_{0}
+ r^{1/2}
\norm{ Av}_{-1/2}
+
\norm{ Av}_{-3/2}
\ge 
\norm{\Phi_0 Av}_{-1/2},
$$
so that
$$
r\norm{\Psi_{0} AEP^*v}_{0}
+ Cr^{1/2}
\norm{ v}_{-1/2}
+
C_{1}\norm{ v}_{-3/2}
\ge\norm{\Phi_0 Av}_{-1/2}\ge C_{2}^{-1}
\norm{B\Phi_0 Av}_{-1/2},
$$
which gives, for all $v\in \mooc(\R^N), \supp v\subset\{y\in\R^N,\val{y-y_{0}}\le r\}$,
$$
r\norm{ P^* v}_{0}
+ r^{1/2}
\norm{ v}_{-1/2}
+
\norm{ v}_{-3/2}
\ge 
\norm{\Phi v}_{-1/2},
\tag 4.4.2$$
where $\Phi=c B\Phi_0A$ is a pseudodifferential operator
of order 0, elliptic near $(y_0,\eta_0).$
By compactness of the cosphere bundle, one gets, using (4.4.2) or (4.4.1),
$$
\norm{v}_{-1/2}
\le C
\sum_{1\le \kappa\le l}
\norm{\Phi_{0\kappa}v}_{-1/2}+C\norm{v}_{-1}
\le C_{1}r \norm{P^*v}_0+C_{1}r^{1/2}\norm{v}_{-1/2}+C_{1}\norm{v}_{-1},
\tag 4.4.3$$
which entails,  by shrinking $r$, the existence of $r_{0}>0, C_{0}>0$,
such that for
$v\in\mooc(\R^N)$, $\supp v\subset\{y\in\R^N,
\val{y-y_{0}}\le r_{0}\}=B_{r_{0}},$
$$
\norm{v}_{-1/2}
\le C_0 \norm{P^*v}_0.
\tag 4.4.4$$
Let $s$ be a real number and $P$ be an operator of order $m$,
satisfying the assumptions of Theorem 1.2.2.
Let $E_{\sigma}$ be a properly supported operator with symbol $\langle \xi\rangle^\sigma$.
Then the operator 
$E_{1-m-s} P E_{s}$ is of first order, satisfies \cps\ and from the previous discussion,
there exists $C_{0}>0, r_{0}>0$
such that
$$
\norm{v}_{-1/2}
\le C_0 \norm{E_{s}P^*E_{1-m-s} v}_0,\quad v\in\mooc(\R^N),\ \supp v\subset B_{r_{0}}.
$$
We get, with $ \chi_{{r}}$ supported in $B_{r}$ and 
$\chi_{r}=1$ on $B_{r/2}$, with $\supp u\subset B_{{r_{0}/4}},$
$$\multline
\norm{\chi_{{r_{0}}}E_{m+s-1}\chi_{{r_{0}/2}}u}_{-1/2}
\le C_0 \norm{E_{s}P^*E_{1-m-s} \chi_{{r_{0}}}E_{m+s-1}\chi_{{r_{0}/2}}u}_0 
\\
\le 
C_0 \nuorm{E_{s}P^*E_{1-m-s}
\underbrace{[ \chi_{{r_{0}}},E_{m+s-1}]\chi_{{r_{0}/2}}}_{S^{-\io}}u}_0
+
C_0 \nuorm{E_{s}P^*\underbrace{E_{1-m-s} E_{m+s-1}}_{{=\Id +S^{-\io}}}\underbrace{\chi_{{r_{0}/2}}u}_{=u}}_0
\\
\le
C_{0}\norm{P^* u}_{s}+\norm{ Ru}_{0},
\endmultline$$
where $R$ is of order $-\io$.
Since we have 
$$
\chi_{{r_{0}}}E_{m+s-1}\chi_{{r_{0}/2}}u=
\underbrace{
[\chi_{{r_{0}}},E_{m+s-1}]\chi_{{r_{0}/2}}}_{S^{-\io}}u+
E_{m+s-1}\underbrace{\chi_{{r_{0}}} \chi_{{r_{0}/2}}u}_{=u},
$$
we get
$
\norm{u}_{s+m-\frac{3}{2}}\le C_{0}\norm{P^*u}_{s} +C_{1}\norm{u}_{s+m-2}
$
and,
shrinking the support of $u$, we obtain the estimate
$$
\norm{u}_{s+m-\frac{3}{2}}\le C_{2}\norm{P^*u}_{s}, 
\tag 4.4.5$$
for $u\in \mooc$ with support in a neighborhood of $y_{0}$.
This implies  the local solvability of $P$, with the loss of derivatives
claimed by the theorem 1.2.2, whose proof is now 
complete.
\head
A. Appendix
\endhead
\subhead
A.1. Wick quantization
\endsubhead
We recall here some facts
on the so-called
Wick quantization, as used in \cite{L4-5-6}.
\definition
{Definition A.1.1}
{
Let 
$Y =(y,\eta)$ be a point 
in $\R^{2n}$. The operator $\Sigma_Y$ is defined as  
$\bigl[2^n e^{-2\pi\val{ \cdot -Y}^2}\bigr]^w$.
This is a rank-one orthogonal projection:
$\Sigma_Y u = (Wu)(Y)\tau_Y\varphi$  with 
$(Wu)(Y)=\poscal{u}{\tau_Y \varphi}_{\l2}$, where
$\varphi(x) = 2^{n/4} e^{-\pi\val x^2}$ and 
$(\tau_{y,\eta}\varphi)(x) = \varphi(x-y) 
e^{2i\pi\langle x-\frac{y}{2},\eta\rangle}.$
Let $a$ be in $L^\infty(\R^{2n})$. 
The Wick quantization of $a$ is defined as
$$
\w{a} =
 \int_{\R^{2n}} a(Y) \Sigma_Y dY.
\tag A.1.1$$
}
\enddefinition
The following proposition is classical and easy 
(see e.g. section 5 in [L5]).
\proclaim{Proposition A.1.2}
\item{1.}
Let $a$ be in $L^\infty(\R^{2n})$. Then $\w a = W^* a^\mu W$
and $\w 1 = {\text{Id}}_{\l2}$
where $W$ is the isometric mapping from $\l2$ to $\L2$
given above, and $a^\mu$ the operator of multiplication by $a$ in $\L2$.
The operator $\pi_H = W W^*$ is the orthogonal projection
on a closed proper subspace $H$ of $\L2$. 
Moreover, we have  
$$\norm{\w a}_{{\Cal L}(\l2)} \le  
\norm a _{L^\infty(\R^{2n})},
\tag A.1.2$$
$$
a(X)\ge 0\  \text{for all}\ X \ \text{implies}\ \w a \ge 0.
\tag A.1.3$$
\item {2.}
Let $m$ be  a real number,and
$p\in S(\Lambda^m, \Lambda^{-1}\Gamma)$. 
Then
$
\w p = p^w + r(p)^w,
$
with $r(p) \in S(\Lambda^{m-1}, \Lambda^{-1}\Gamma)$ 
so that the mapping $p\mapsto r(p)$ is continuous.
More precisely, one has
$$
r(p)(X)=
\int_0^1 \!\!\!\int_{\RZ} 
(1-\theta) p''(X + \theta Y) Y^2 
e^{-2\pi\Gamma(Y)} 2^n dY d\theta.
$$
Note that $r(p)=0$ if $p$ is affine.
\item{3.}
For $a\in L^\io(\RZ)$,
the Weyl symbol
of
$\w{a}$
is
$$
a\ast 2^n\exp-2\pi\Gamma\quad
\text{which belongs to $S(1,\Gamma)$ 
with $k^{th}$-seminorm 
$c(k)\norm{a}_{L^\io}$.}
\tag A.1.4$$
\item{4.}
Let $\R\ni t\mapsto a(t,X)\in \R$
such that, for $t\le s$,
$a(t,X)\le a(s,X)$. Then,
for $u\in C^1_c\bigl(\R_t,L^2(\R^n)\bigr)$,
assuming $a(t,\cdot)\in L^\io(\RZ)$, 
$$
\int_\R\re\poscal{D_t u(t)}{i\w{a(t)}
 u(t)}_{L^2(\R^n)}dt \ge 0.
\tag A.1.5$$
\item{5.} With the operator $\Sigma_Y$ 
given in definition A.1.1,
we have the estimate
$$
\norm{\Sigma_Y \Sigma_Z
}_{{\Cal L}(\l2)} 
\le  2^n e^{-\frac{\pi}{2}\Gamma(Y-Z)}.\hskip6pt 
\tag A.1.6$$
\item{6.} More precisely, the Weyl symbol of
$\Sigma_Y \Sigma_Z $ 
is, as a function of the variable $X\in \RZ$,
setting $\Gamma(T)=\val T^2$
$$
 e^{-\frac{\pi}{2}\val{Y-Z}^2} e^{-2i\pi[X-Y,X-Z]}
2^n e^{-2\pi\val{ X -\frac{Y+Z}{2}}^2}.
\tag A.1.7$$
\endproclaim
Since for the Weyl quantization, one has
$
\norm{a^{w}}_{{\Cal L}(\l2)} 
\le  2^n \norm a _{L^1(\RZ)},
$
we get the result (A.1.7) from (A.1.6).
Note that (A.1.5)
is simply a way of writing
that
$
\frac{d}{dt}\bigl(\w{a(t)}\bigr)\ge 0,
$
which is a consequence of (A.1.3)
and of
the non-decreasing assumption
made on $t\mapsto a(t,X)$.
\proclaim{Lemma A.1.3}
Let $M$ be a $\Gamma$-weight, i.e. a positive function such that
$
M(X) M(Y)^{-1}\le C(1+\Gamma(X-Y))^N
$
(see definition 1.3.1). Then if a measurable
function $a$ defined on $\RZ$
satisfies
for all $X$, $\val{a(X)}\le C_{1}M(X)$,
the symbol $a\ast \exp-2\pi\Gamma$
belongs to $S(M,\Gamma)$
with semi-norms depending only on $C_{1}$.
More generally, for a polynomial $p$ the symbol
$A$ defined by
$$
A(X)=\int a(Y) p(X-Y)\exp-2\pi\Gamma(X-Y) dY
$$
belongs to  to $S(M,\Gamma)$.
\endproclaim
\demo{Proof}
We check first
$$
(a\ast 2^n\exp-2\pi\Gamma)^{(k)}(X)=\int a(Y)P_{k}(X-Y)2^n\exp-2\pi\Gamma(X-Y) dY
\tag A.1.8$$
with a polynomial
$P_{k}$, which gives
$$\multline
M(X)^{-1}\val{(a\ast 2^n\exp-2\pi\Gamma)^{(k)}(X)}\le C_{1}\int\frac{M(Y)}{M(X)}
\val{P_{k}(X-Y)}2^n\exp-2\pi\Gamma(X-Y) dY
\\\le
C_{1}\int C\bigl(1+\Gamma(X-Y)\bigr)^{N}
\val{P_{k}(X-Y)}2^n\exp-2\pi\Gamma(X-Y) dY
=C_{1}C\gamma(k,N,n).
\endmultline$$
Let us examine
$A^{(k)}$: it is a sum of terms of type (A.1.8) and thus the above argument works.
\qed
\enddemo
\proclaim{Lemma A.1.4}
Let $g$
be an admissible metric on $\RZ$
(see definition 1.3.1)
such that, with $\Gamma$ a given symplectic norm,
there exists $C_0> 0$,
$n_0\ge 0$ such that
$$
\forall X,Y,T,\
g_X(T)\le C_0 \Gamma(T),\quad\frac{g_X(T)}
{g_Y(T)}\le C_0\bigl(1+\Gamma(X-Y)\bigr)^{n_0}.
\tag A.1.9$$
Let $m$ be a weight for
$g$
(definition 1.3.1)
such that
$$
\frac{m(Y)}{m(Z)}
\le C_0 \bigl(1+\Gamma(Z-Y)\bigr)^{n_0}.
\tag A.1.10$$
Then,
if 
$A\in \op
{S(m,g)}$,
there exists
a semi-norm
$\gamma$
 of the 
symbol of $A$ such that
$$
\val{\poscal{Av}{v}}\le \gamma
\poscal{\w{m} v}{v}=\gamma\int_{\RZ}m(Y)
\norm{\Sigma_Y v}_{L^2}^2dY.
\tag A.1.11$$
\endproclaim
\demo{Proof}
Theorem 6.9
in \cite{BC}
shows
that the space
 $\Cal H(m^{1/2},g)$ 
is equal to 
$\Cal H(m^{1/2},\Gamma)$
provided that $m^{1/2}$ is regular.
In fact we may assume that
$m$ is regular
since it is anyhow always equivalent
to a regular weight.
Using definition 7.1 in
\cite{BC}, we check that
$g$
is dominated by a
strongly temperate metric, 
namely the constant metric $\Gamma$.
Moreover the corollary 6.7 and theorem 7.8 in \cite{BC}
imply
$$
\val{\poscal{Av}{v}}
\le \nuorm{Av}_{\Cal H(m^{-1/2},g)}
\nuorm{v}_{\Cal H(m^{1/2},g)}
\le \gamma
\nuorm{v}_{\Cal H(m^{1/2},g)}^2
= \gamma
\nuorm{v}_{\Cal H(m^{1/2},\Gamma)}^2
=
\gamma\int_{\RZ}
m(Y)\norm{\theta_Y^w u}_{L^2}^2 dY,
$$
where $(\theta_Y)$
is a partition of unity related 
to the metric $\Gamma$.
We have, using the results of this section,
(A.1.10)
 and 
(A.1.6),
with 
$\langle T\rangle^2=1+\Gamma(T)$,
for all $N_1,N_2$,
$$\align
\int
m(Y)\norm{\theta_Y^w u}^2dY
&=
\iiint_{\RZ}m(Y)
\poscal{\theta_Y^w\Sigma_{Z_1}
\Sigma_{Z_1} u}{\theta_Y^w\Sigma_{Z_2}
\Sigma_{Z_2} u} 
dYdZ_1dZ_2
\\
&\le
\iiint m(Z_1)^{1/2}m(Z_2)^{1/2}
\norm{\Sigma_{Z_1}u}
\norm{\Sigma_{Z_2}u}
\langle Y-Z_1\rangle^{-N_1}
\langle Z_2-Z_1\rangle^{-N_2}
dYdZ_1dZ_2 C_{N_1,N_2}
\\
&\le
\iint m(Z_1)^{1/2}m(Z_2)^{1/2}
\norm{\Sigma_{Z_1}u}
\norm{\Sigma_{Z_2}u}
\langle Z_2-Z_1\rangle^{-N_2}
dZ_1dZ_2 C_{N_1,N_2}
\\
&\le
\int m(Z)
\norm{\Sigma_{Z}u}^2
dZ,
\endalign$$
which completes the proof of the lemma.
\qed
\enddemo
\proclaim{Lemma A.1.5}
Let $m_{1}, m_{2}$ be two $\Gamma$-weights (see definition 1.3.1)
and $a_{1},a_{2}$ be two locally Lipschitz continuous functions such that
$
\val{a_{1}(X)}\le m_{1}(X),\quad
\val{a'_{2}(X)}\le m_{2}(X).
$
Then the operator
$$
\w{a_{1}}\w{a_{2}}\in\w{(a_{1}a_{2})}+\op{S(m_{1}m_{2},\Gamma)}.
\tag A.1.12$$
\endproclaim
\demo{Proof}
We use the definition A.1.1 and Taylor's
formula
to write
$$
\w{a_{1}}\w{a_{2}}=\iint a_{1}(Y)\Bigl(a_{2}(Y)+
\int_{0}^{1}a'_{2}\bigl(Y+\theta(Z-Y)\bigr)d\theta(Z-Y)\Bigr)\Sigma_{Y}\Sigma_{Z}dY dZ
=\w{(a_{1}a_{2})}+R^w,
$$
with
$$R(X)=\iiint_{0}^1
a_{1}(Y)a'_{2}\bigl(Y+\theta(Z-Y)\bigr)(Z-Y)
e^{-\frac{\pi}{2}\val{Y-Z}^2} e^{-2i\pi[X-Y,X-Z]}
2^n e^{-2\pi\val{ X -\frac{Y+Z}{2}}^2}dY dZd\theta.
\tag A.1.13$$
We have, using (5) in definition 1.3.1,
$$\align
\val{R(X)}
&\le
\iiint_{0}^1 m_{1}(Y)m_{2}(Y)\frac{m_{2}(Y+\theta(Z-Y))}{m_{2}(Y)}\val{Z-Y}
e^{-\frac{\pi}{2}\val{Y-Z}^2} 
2^n e^{-2\pi\val{ X -\frac{Y+Z}{2}}^2}dY dZd\theta
\\
&\le C m_{1}(X)m_{2}(X)\iiint_{0}^1 
(1+\val{Y-X}^2)^N(1+\val{Y-Z}^2)^{N+1/2}
e^{-\frac{\pi}{2}\val{Y-Z}^2} 
e^{-2\pi\val{ \frac{Y+Z}{2}-X}^2}dY dZd\theta
\\&=
C m_{1}(X)m_{2}(X)\iint(1+\val{T/2+S}^2)^N(1+\val{T}^2)^{N+1/2}
e^{-\frac{\pi}{2}\val{T}^2} 
e^{-2\pi\val{ S}^2}dT dS
\\&=C'm_{1}(X)m_{2}(X).
\endalign$$
Moreover taking derivatives of $R$ in its defining formula (A.1.13)
above leads to the same estimate for $R^{(k)}(X)$. The proof of the lemma is complete.
\qed
\enddemo
\proclaim{Lemma A.1.6 }
Let
$(\chi_{k})$ be a partition of unity and $(\psi_{k})$ be a sequence
  as in lemma 1.4.1
for an admissible metric of type
$\lambda^{-1}(X)\Gamma$,
where $\lambda$ is a $\Gamma$-weight and 
$\Gamma=\Gamma^\sigma.$ 
Let $\omega$ be a locally bounded function
such that
$
\val{\omega(X)}\le M(X)
$
where $M$ is a $\Gamma$-weight.
Assume that, for each $k$, there
exist a bounded function $\omega_{k}$
such that
$\omega(X)=\omega_{k}(X)$ for all $X\in \supp{\chi_{k}}$
and such that for all $X\in \RZ$,
$\val{\omega_{k}(X)}\le M(X) \lambda(X)^{N_{0}}.$
Then with
$
\widetilde{\omega}(X)=\int \omega(Y)2^n\exp-2\pi\Gamma(X-Y) dY,$
we have
$$
\chi_{k}(X)\widetilde{\omega}(X)=\chi_{k}(X)\widetilde{\omega_{k}}(X)+r_{k}(X),\
\sum_{k}r_{k}\in S(\lambda^{-\io},\Gamma).
\tag A.1.14$$
\endproclaim
\demo{Proof}
We already know from the lemma A.1.3 that
$X\mapsto \widetilde{\omega}(X)=\int \omega(Y)2^n\exp-2\pi\Gamma(X-Y) dY
$
belongs to $S(M,\Gamma)$.
We check now
$$\align\chi_{k}(X)\widetilde{\omega}(X)&=
\chi_{k}(X)\int \omega(Y)2^n\exp-2\pi \Gamma(X-Y)dY
\\&=
\chi_{k}(X)\int \psi_{k}(Y)\omega(Y)2^n\exp-2\pi \Gamma(X-Y)dY
\\&\hskip55pt+
\chi_{k}(X)\int_{ Y, \psi_{k}(Y)\not=1} (1-\psi_{k}(Y))\omega(Y)2^n\exp-2\pi \Gamma(X-Y)dY
\\&=\chi_{k}(X)\int \psi_{k}(Y)\omega_{k}(Y)2^n\exp-2\pi \Gamma(X-Y)dY+r_{k}(X).\tag A.1.15
\endalign$$
We have 
$
\Gamma(U_{k}-(U_{k}^*)^c)=\inf_{\Gamma(T)<1\le \Gamma(S)}\Gamma(X_{k}+r_{0}\lambda(X_{k})^{1/2}T
-X_{k}-\lambda(X_{k})^{1/2}2r_{0}S)
$
and thus
$
\Gamma(U_{k}-(U_{k}^*)^c)\ge \lambda(X_{k})r_{0}^2.
$
Since
$\psi_{k}$
is equal to 1 on $U_{k}^*$
(notations of section 1.4)
we obtain from (A.1.15)
$$
\val{r_{k}^{(j)}(X)}_{\Gamma}\le C_{j}\psi_{k}(X)\exp-\pi\Gamma(U_{k}-(U_{k}^*)^c)
\le C_{{j,N,r_{0}}}\psi_{k}(X)\lambda(X)^{-N}
$$
and thus
$
\sum_{k} r_{k}\in S(\lambda^{-\io}, \Gamma).
$
We obtain
$$
\chi_{k}\widetilde{\omega}=\chi_{k}\bigl(\psi_{k}\omega_{k}
\ast 2^n\exp-2\pi\Gamma\bigr)+r_{k}
=\chi_{k}\bigl(\omega_{k}
\ast 2^n\exp-2\pi\Gamma\bigr)+
\chi_{k}\bigl(\omega_{k}(\psi_{k}-1)
\ast 2^n\exp-2\pi\Gamma\bigr)+
r_{k},
$$
and applying again the same reasoning to the penultimate term above,
we get for $Y\in {(U_{k}^*)}^c$ and $X\in U_{k}$, that
$
\Gamma(X-Y)\ge \lambda(X_{k})r_{0}^2 
$
the following estimate
for the integrand
$$\align
\exp-\pi\Gamma(X-Y)&\exp-\pi\lambda(X_{k})r_{0}^2\times M(Y)\lambda(Y)^{N_{0}}
\\&\le
CM(X)\lambda(X)^{N_{0}}(1+\Gamma(X-Y))^{N_{0}}
\exp-\pi\Gamma(X-Y)\exp-\pi\lambda(X_{k})r_{0}^2
\\
&\le
C'M(X)\lambda(X_{k})^{N_{0}}
(1+\Gamma(X-X_{k}))^{N_{0}}
\exp-\frac{\pi}{2}\Gamma(X-Y)\exp-\pi\lambda(X_{k})r_{0}^2
\\&\le
C''M(X)\lambda(X_{k})^{3N_{0}}
\exp-\frac{\pi}{2}\Gamma(X-Y)\exp-\pi\lambda(X_{k})r_{0}^2
\\&\le
C'''M(X)\lambda(X_{k})^{3N_{0}}
\exp-\frac{\pi}{2}\Gamma(X-Y)\exp-\pi\lambda(X_{k})r_{0}^2
\endalign
$$
which yields the result.
\qed
\enddemo
\definition{Definition A.1.7}
Let $\Gamma$ be a symplectic quadratic form  on $\R^n\times \R^n$, i.e. a positive definite quadratic form
such that $\Gamma=\Gamma^\sigma$(see definition 1.3.2(2)).
There exists a unique linear symplectic mapping $A$ such that for all $X=(x,\xi)$,
$\Gamma(AX)=\sum_{1\le j\le n}x_{j}^2+\xi_{j}^2.
$
Let $U$ be a metaplectic transformation in the fiber of $A$. Then for $a\in L^{\io}(\RZ)$,
we define
$$
a^{\text{Wick($\Gamma$)}}=\int a(Y) 2^n
\bigl(\exp-2\pi\Gamma(\cdot-Y)\bigr)^wdY
=U\w{(a\circ A)} U^\ast.
\tag A.1.16$$
\enddefinition
\remark{Remark A.1.8}
Note that since $U$ is uniquely determined up to a factor of modulus one, that definition is consistent.
We remark also that, defining for $X\in\RZ$,
$
\Phi(X)= 2^n\exp-2\pi\Gamma(X),
$
we have
$
\Phi(AX-AY)=2^n\exp-2\pi\val{X-Y}^2,
$
which is the Weyl symbol of $\Sigma_{Y}$
(definition A.1.1).
From the Segal formula, we have, with a metaplectic $U$ in the fiber of
$A$
$$
\Phi(X -Z)^w=U \Phi(AX -Z)^wU^\ast
$$
and thus we can justify the equality in formula (A.1.16) since
$$\multline\int a(Y) 2^n\bigl(
\exp-2\pi\Gamma(X-Y)\bigr)^w dY
=
\int a(AY) \Phi(X-AY)^w dY
\\=
\int a(AY) U \Phi(AX -AY)^wU^\ast
=\int a(AY) U \Sigma_{Y}U^\ast dY=
U\w{(a\circ A)} U^\ast.
\endmultline$$ 
\endremark
\remark{Remark A.1.9}
We can also notice that the definition above is consistent with the fact that
Wick and Weyl quantization coincide
for linear forms: if $a$ is a linear form, we have
$$
a^{\text{Wick($\Gamma$)}}
=U\w{(a\circ A)} U^\ast=U{(a\circ A)^w} U^\ast=UU^\ast a^wUU^\ast=a^w.
\tag A.1.17$$
Also, it is easy with the formula (A.1.16) to check that the results of section A.1
on the Wick quantization can be extended,
{\it mutatis mutandis},
to the Wick($\Gamma$) quantization.
\endremark
\subhead
A.2. Properties of some metrics
\endsubhead\par\no
{\bf Proof of the remark following definition 1.3.1.}
Using a partition of unity related to the slowly varying $g$,
as in \cite{BL}, we define
$
M_*(X)=\int_{\RZ} M(Y)\varphi_Y(X)\val{g_Y}^{1/2}dY.
$
It is a simple matter left to the reader to check that
$M_*$ belongs to $S(M,g)$
and satisfies (1.3.3).
\proclaim{Lemma A.2.1} Let $\Gamma$ be a positive definite quadratic form
on $\RZ$ such that
$\Gamma =\Gamma^\sigma$
and let $g_X=\lambda(X)^{-1}\Gamma$ be a metric conformal to $\Gamma$ such that
$g$ is slowly varying and $\inf_X\lambda(X)\ge 1$.
Then the metric $g$ satisfies\quad
$
g_X(T)\le Cg_Y(T)\bigl(1+\Gamma(X-Y)\bigr),
$
i.e. 
$$
\frac{\lambda(Y)}{\lambda(X)}\le C\bigl(1+\Gamma(X-Y)\bigr),
\tag A.2.1$$
implying that $g$ is admissible.
\endproclaim
\demo{Proof}
Since $g$ is slowly varying, we may assume, with a positive $r_0$,
$g_Y(Y-X)\ge r_0^2$, which means
$
{\Gamma(Y-X)}\ge r_0^2{\lambda(Y)}
$
and using
$ {\lambda(X)}\ge 1$ 
we get
${\lambda(Y)}/{\lambda(X)}\le r_0^{-2} \Gamma(Y-X).$
\qed
\enddemo
\proclaim{Lemma A.2.2}
 Let $\Gamma$ be a positive definite quadratic form
on $\RZ$ such that
$\Gamma =\Gamma^\sigma$
and let $g_X=\lambda(X)^{-1}\Gamma$ be a metric conformal to $\Gamma$.
 Assume that
 $
 \lambda(X)=d(X)^2+\lambda_1(X)
 $
 with a  function
 $d$ uniformly Lipschitz continuous (with respect to $\Gamma$)
 and
 $\lambda_1^{-1}\Gamma$ slowly varying with
 $\lambda_1\ge 1$.
 Then the metric $g$ is slowly varying.
\endproclaim
\demo{Proof}
Let us assume that
$\val{X-Y}^2\le r^2 \bigl(d(X)^2+\lambda_1(X)\bigr)
$.
If $d(X)^2\le \lambda_1(X)$,
using the fact that
$\lambda_1^{-1}\Gamma$ is slowly varying,
we can choose $r$ small enough so that
$
\lambda_1(X)\le C_1\lambda_1(Y)
$
and thus
$$
\lambda(X)\le 2C_1\lambda_1(Y)\le 2C_1\lambda(Y).
$$
If $d(X)^2>\lambda_1(X)$,
we have, with $L$ standing for the Lipschitz constant of $d$,
$$
2^{-1/2}\lambda(X)^{1/2}< d(X)\le 
d(Y)+L\val{X-Y}\le 
\lambda(Y)^{1/2}+
L r\lambda(X)^{1/2}$$
so that, for $r\le \frac{1}{2^{3/2}L+1}$
we get
$
\lambda(X)\le 8\lambda(Y).
$
\remark{Remark A.2.3}
It is a simple exercise left to the reader
to show that (1) in Definition 1.3.1
is satisfied whenever there exists $r_0>0, C_0>0$ such that
for all $X,Y,T\in\RZ$,
 $g_X(Y-X)\le r_0^2$ implies
$g_Y(T)\le C_0 g_X(T)$.
\endremark
Taking this remark into account,  we complete 
the proof of the lemma.\qed
\enddemo
\subhead A.3. Proof of Lemma 2.1.5 on the proper class
\endsubhead
All norms in this proof are taken with respect to the constant
quadratic form $\Gamma$, so we omit the index everywhere
and denote $\norm{\cdot}_\Gamma$ by $\val{\cdot}$.
Since for all $j\in \N$,
$
\val{f^{(j)}(X)}\le \gamma_j\Lambda^{m-\frac{j}{2}},
$
we get
$
1\le \lambda(X)\le 1+\Lambda
\max_{0\le j<2m\atop j\in\N }\gamma_j^{\frac{2}{2m-j}}=1+\gamma \Lambda\le(1+\gamma)\Lambda
$
and (2.1.12).
For $0\le j<2m$, we have from the definition of $\lambda$, the estimate
$
\val{f^{(j)}(X)}\le\lambda(X)^{m-\frac{j}{2}},
$
and for $j\ge 2m$, we can use
$$
\val{f^{(j)}(X)}\le\gamma_j\Lambda^{m-\frac{j}{2}}=
\gamma_j\Lambda^{-\frac{(j-2m)}{2}}\le
\gamma_j\lambda^{-\frac{(j-2m)}{2}}(1+\gamma)^{\frac{(j-2m)}{2}},
$$
so that $f\in S(\lambda^m,\lambda^{-1}\Gamma)$
with a $j$-th semi-norm less than $1$ for $j<2m$ and less than
$\gamma_j(1+\gamma)^{\frac{(j-2m)}{2}}$
 for $j\ge 2m$ .\par
 Let us now prove that
 $\lambda^{-1}\Gamma$ is slowly varying. Let us assume that
 $\val{X-Y}^2\le r^2 \lambda(X)$. 
 Using Taylor's formula, we get for  the smallest
 integer 
 $N\ge 2m$ ($N=-[-2m]$) and $0\le j<2m$,
 $$
 \val{f^{(j)}(X)}\le\sum_{l, j+l<2m}\val{f^{(j+l)}(Y)}\frac{r^l}{l!}
 \lambda (X)^{l/2}
 +\gamma_N
\Lambda^{m-\frac{N}{2}}\frac{r^{N- j}}{(N -j)!}
 \lambda (X)^{(N- j)/2},
 $$
 so that 
 $
 \val{f^{(j)}(X)}\le\sum_{l, j+l<2m}\lambda (Y)^{\frac{2m-j-l}{2}}
  \lambda (X)^{\frac{l}{2}}
 \frac{r^l}{l!}
  +\gamma_N
\Lambda^{\frac{2m-N}{2}}
 \lambda (X)^{\frac{N- j}{2}}\frac{r^{N- j}}{(N -j)!},
 $
 and 
$$
\multline
 \val{f^{(j)}(X)}\le
  \sum_{l, j+l<2m}
 ( \lambda (Y)^{\frac{2m-j}{2}})^{\frac{2m-j-l}{2m-j}}
  (\lambda (X)^{\frac{2m-j}{2}})^{\frac{l}{2m-j}}
 \frac{r^l}{l!}
  +\gamma_N
\Lambda^{\overbrace{\scriptstyle\frac{2m-N}{2}}^{\le 0}}
 \lambda (X)^{\frac{N- j}{2}}\frac{r^{N- j}}{(N -j)!}\\
 \le \sum_{l, j+l<2m}\!\!\!
  {\frac{2m-j-l}{2m-j}}\lambda (Y)^{\frac{2m-j}{2}}
  \frac{r^l}{l!}
+
  {\frac{l}{2m-j}}
  \lambda (X)^{\frac{2m-j}{2}}
 \frac{r^l}{l!}
    \\
    \qquad+\gamma_N(1+\gamma)^{\frac{N-2m}{2}}
 \lambda (X)^{\frac{2m- j}{2}}\frac{r^{N- j}}{(N -j)!}
 \endmultline
 $$
 implying
 $$
 \multline
 \val{f^{(j)}(X)}
 \\
 \le
\lambda (Y)^{\frac{2m-j}{2}}
\overbrace{
 \sum_{l, j+l<2m}
  {\frac{2m-j-l}{2m-j}}
  \frac{r^l}{l!}}^{\text{$=p(r)$ a polynomial in $r$}}
  + \lambda (X)^{\frac{2m-j}{2}}
  \underbrace{
  \Bigl(\sum_{1\le l, j+l<2m}
   {\frac{l}{2m-j}}
  \frac{r^l}{l!}  +\gamma_N(1+\gamma)^{\frac{N-2m}{2}}
 \frac{r^{N- j}}{(N -j)!}\Bigr)}_{\text{$=\epsilon(r)$ goes to zero with $r$.}
}. 
\endmultline
$$
  Assuming then that $j$ was chosen so that
  $\lambda(X)=1+ \val{f^{(j)}(X)}^{\frac{2}{2m-j}}$,
  we get
$$\lambda(X)\le 1+\Bigl(
\lambda(Y)^{\frac{2m-j}{2}}p(r)
+  \lambda(X)^{\frac{2m-j}{2}}
\epsilon(r)\Bigr)^{\frac{2}{2m-j}}, $$
  so that there exist
  $r_0>0, C_0\ge 1$,
  depending only on the $N$ first semi-norms of 
  $f$, such that  for
  $r\le r_0,$ 
  we have $$ 
  \val{X-Y}^2\le r^2\lambda(X)\Longrightarrow
  \lambda(X)\le C_0
\lambda(Y),$$
and thus
$
 r\le r_0,
  \val{X-Y}^2\le r^2C_0^{-1}\lambda(X)
\Longrightarrow
 C_0^{-1}\lambda(X)\le \lambda(Y)\le C_0
\lambda(X),
$
which is the property (1) in Definition A.1.1.
The property (2) in that definition is obviously satisfied since $\lambda(X)\ge 1$.
Moreover, we get a stronger property than (3) from the Lemma A.2.1 above in this appendix.
\qed
\subhead
A.4. Some a priori estimates
and loss of derivatives
\endsubhead
In this section,
we prove that, at least when a factorization occurs,
it is possible to limit the loss of derivatives
to 3/2
(the loss is always counted with respect to the elliptic case).
Let us study the model-case
$$
L=D_t+iA_0 B_1, \quad A_0\in \op{S^0}, B_1\in \op{S^1}
$$
with real-valued Weyl symbols such that  $A_0\ge c_0\Lambda^{-1}, \dot{B_1}\ge 0$.
We compute, using the notation
$$
\norm{u}=\left(\int \val{u(t)}^2dt\right)^{1/2},\quad
\val{v}=\norm{v}_{\Bbb H},\quad \Bbb H=L^2(\R^{n}), \quad
\val{u}_\io=\sup_{t\in \R}\val{u(t)},
$$
$$
2\re\poscal{Lu}{iB_1u}=\poscal{\dot{B_1}(t) u(t)}{u(t)}+2\re
\poscal{A_0 B_1 u}{B_1 u}\ge 2c_0\Lambda^{-1}\norm{B_1 u}^2.
$$
As a consequence, for $\supp{u}\subset[-T, T]$,
$$\multline
2\re\poscal{Lu}{iB_1u}+2\re\poscal{Lu}{iH(t-T_0)u}
\\
\ge
c_0\Lambda^{-1}\norm{B_1 u}^2+\val{u}_\io^2
+\nuorm{A_{0}^{1/2}B_{1}u}^2
+2\re\poscal{A_{0}^{1/2}B_{1}u}{i H_{T_{0}}A_{0}^{1/2}u}
\\
\ge
c_0\Lambda^{-1}\norm{B_1 u}^2
+\val{u}_\io^2(1-\sup_{\val t\le T}\norm{A_0(t)}T)
\\ \text{\sevenrm (for T small enough) }\quad
\ge
c_0\Lambda^{-1}\norm{B_1 u}^2
+\frac{1}{2}\val{u}_\io^2 ,
\endmultline$$
so that
$
c_0^{-1}\Lambda{}\norm{Lu}^2+
c_0{}\Lambda^{-1}\norm{B_1u}^2+2\norm{Lu}\norm{u}
\ge
c_0\Lambda^{-1}\norm{B_1 u}^2
+\frac{1}{2}\val{u}_\io^2 
$
and thus
$$
(c_0^{-1}\Lambda{}+1)\norm{Lu}^2+
T\val{u}_\io^2\ge \frac{1}{2}\val{u}_\io^2 
$$
entailing for $T\le 1/4$,\quad
$
(c_0^{-1}\Lambda{}+1)\norm{Lu}^2\ge\frac{1}{4}\val{u}_\io^2 ,$
which gives
$
\norm{Lu}\gg \Lambda^{-1/2}\norm u
$,
an estimate with loss of 3/2 derivatives.\par
The next question is obviously: how do we manage to get the estimate
$A_0\ge \Lambda^{-1}$? Assuming
$A_0\ge -C\Lambda^{-1}$, 
we can always consider instead
$A_0+(C+1)\Lambda^{-1}\ge \Lambda^{-1}$; now this modifies the operator $L$ and although
our estimate
is too weak to absorb a zeroth order
perturbation,
it is enough to check that the energy method is stable by zeroth order perturbation.
We consider then
$$
D_t+iA_0 B_1+S+iR,\quad
A_0\ge \Lambda^{-1},\quad S,R\in \op{S^0}.
$$
Inspecting the method above, we see that
$S$ will not produce any trouble, since we shall commute it with $B_1$, 
producing an operator of order 0.
The term produced by $R$ are more delicate to handle: we shall have to deal with
$$
2\poscal{Ru}{B_1u}+2\poscal{Ru}{H_{T_0}u}.
$$
The second term is $L^2$ bounded and can be absorbed.
There is no simple way to absorb the first term, which is of size
$\norm{B_1u}\norm u$
which is too large with respect to the
terms that we dominate.
However we can consider the $L^2$-bounded invertible operator
$U(t)$ (which is in $\op{S^0}$ and self-adjoint) such that
$U(0)=\Id$
and 
$\dot U(t)=-U(t)R(t)$
so that
$$\multline
L =D_t+iR+iA_0 B_1+S= U(t)^{-1}D_t U(t)+iA_0 B_1+S
\\=
U(t)^{-1}\bigl(D_t+iA_0 B_1 +S\bigr)U(t)
-U(t)^{-1}\bigl[iA_0B_1+S,U(t)\bigr]
\\
=
U(t)^{-1}\Bigl(D_t+iA_0 B_1 +S+\bigl[U(t),iA_0B_1+S\bigr]U(t)^{-1}\Bigr)U(t).
\endmultline$$
Now the term
$\bigl[U(t),iA_0B_1\bigr]U(t)^{-1}$
has a real-valued principal symbol in $S^0$ and amounts to a modification
of S, up to unimportant terms of order $-1$.
The term
$\bigl[U(t),S\bigr]U(t)^{-1}$
is of order $-1$ and can be absorbed.
We have proven the following lemma.
\proclaim{Lemma A.4.1}
Let $\Lambda\ge 1$ be  given.
We consider the metric $G=\val{dx}^2+\Lambda^{-2}\val{d\xi}^2$ on $\R^n\times \R^n$.
Let $a_{0}(t,x,\xi)$ be in $S(1,G)$ such that
$a_{0}(t,x,\xi)\ge 0$.
Let $b _{1}(t,x,\xi)$ be real-valued and in $S(\Lambda,G)$ such that
$$\bigl(b_{1}(t,x,\xi)-b(s,x,\xi)\bigr)(t-s)\ge 0.$$
Let $r(t,x,\xi)$ be a complex-valued symbol in $S(1,G)$.
Assuming that $a_{0}, b_{1}, r_{0}$ are continuous functions, there exists a constant $C>0$ depending only on the semi-norms of the symbols
$a_{0}, b_{1}, r_{0}$, such that, for all 
$u\in C^1_{c}([-T,T],L^2(\R^n))$ with $CT\le 1$, 
$$
C\norm{Lu}_{{L^2(\R^{n+1})}}
\ge \Lambda^{-1/2}T^{-1}\norm{u}_{L^2(\R^{n+1})}.
$$
\endproclaim
\subhead
A.5. Some lemmas on symbolic calculus
\endsubhead
Let $g$ be an admissible metric on
$\RZ$ and $m$ be a $g$-weight (see definition 1.3.1).
Then, 
at each
point $X\in\RZ$,
we can define a metric $g_X^\sharp$ by
taking
the geometric mean of $g_X,g_X^\sigma$
so that in particular
$$
g_X\le g_X^\sharp=(g_X^\sharp)^\sigma\le 
g_X^\sigma.
\tag A.5.1$$
We define
$$
h(X)=
\sup_{g^\sharp_X(T)=1}
{g_X(T)}
\tag A.5.2$$
and we note that whenever
$g^\sigma=\lambda^2 g$
we get from the definition 1.3.1 that
$g^\sharp=\lambda_{g} g$ and 
$\lambda_{g}=1/h$.
\definition{Definition A.5.1}
Let $l$ be a nonnegative integer.
We define the set
$S_{l}(m,g)$ as the set of smooth functions $a$ defined on $\RZ$
such that $a$ satisfies the estimates of $S(m,g)$ for derivatives of order $\le l$,
and 
the estimates
of
$S(m,g^\sharp)$
for derivatives of order $\ge l+1$,
which means
$$
\val{a^{(k)}(X)T^k}\le C_km(X)\times
\cases 
g_X(T)^{k/2}&\text{if $k\le l,$}
\\
g_X^\sharp(T)^{k/2}h(X)^{\frac{l+1}{2}}&\text{if $k\ge l+1$, with 
$h(X)=\dis
\sup_{g^\sharp_X(T)=1}
{g_X(T)}$.}
\endcases
$$
Note that since $h\le 1$ and
$
g\le h g^\sharp
$,
we get
$S(m,g)\subset S_{l}(m,g)$.
If $g=\lambda(X)^{-1}\Gamma_{X}$, where $\lambda(X)$ is positive (scalar)
and 
$\Gamma_{X}=\Gamma_{X}^\sigma$,
then
$g_{X}^\sharp=\Gamma_{X}$
and $a$ belongs to $S_{l}(m, \lambda^{-1}\Gamma)$
means
$$
\val{a^{(k)}(X)}_{{\Gamma_{X}}}\le C_{k} m(X)
\times
\cases 
\lambda(X)^{-k/2}&\text{if $k\le l,$}
\\
\lambda(X)^{-l/2}&\text{if $k\ge l+1$}.
\endcases
$$
Moreover, if $g\equiv g^\sharp$,
then for all $l$, 
$S(m,g)=S_{l}(m,g)$.
\enddefinition
\proclaim{Lemma A.5.2}
Let $\Gamma$ be a positive definite quadratic form on $\RZ$ such that
 $\Gamma=\Gamma^\sigma$
 and $\lambda$ be a $\Gamma$-weight. 
Let $b$ be a symbol in $S_{1}(\lambda^m,\lambda^{-1}\Gamma)$, where $m$ is a real number.
Then
$
b\sharp b-b^2\in S(\lambda^{2m-1},\Gamma)
$
\endproclaim
\demo{Proof}
We have
$
(b\sharp b)(X)=\exp i\pi[D_{{X_{1}}},D_{X_{2}}]\bigl(b(X_{1})\otimes
b(X_{2})\bigr)_{\vert X_{1}=X_{2}=X}
$
so that using Taylor's formula with integral remainder for $s\mapsto e^s$
yields
$$
(b\sharp b)(X)= b(X)^2+\int_{0}^1\exp i\pi\theta [D_{{X_{1}}},D_{X_{2}}]d\theta
i\pi[D_{{X_{1}}},D_{X_{2}}]
b(X_{1})\otimes
b(X_{2})_{\vert X_{1}=X_{2}=X}.
$$
Since $b'\in S(\lambda^{m-1/2},\Gamma)$
and 
$$\multline
\exp i\pi\theta [D_{{X_{1}}},D_{X_{2}}](a_{1}(X_{1})\otimes a_{2}(X_{2}) )\\=
\exp i\pi [D_{{X_{1}}\theta^{-1/2}},D_{X_{2}\theta^{-1/2}}]
(a_{1}(\theta^{-1/2}X_{1}\theta^{1/2})\otimes a_{2}(\theta^{-1/2}X_{2}\theta^{1/2}) )_{\vert X_{1}=X_{2}=X}
\\
=
\exp i\pi [D_{{Y_{1}}},D_{Y_{2}}](a_{1}(\theta^{1/2}Y_{1})\otimes a_{2}(\theta^{1/2}Y_{2}) )
_{\vert Y_{1}=Y_{2}=\theta^{-1/2}X}
\\
=\bigl((a_{1}\circ \theta^{1/2})\sharp (a_{2}\circ \theta^{1/2})\bigr)(\theta^{-1/2}X),
\endmultline$$
we get that, if $a_{j}\in S(\lambda^{m_{j}},\Gamma)$,
we have
$
a_{j}\circ \theta^{1/2}\in S(\lambda^{m_{j}},\theta \Gamma)
$
so that 
the symbolic calculus for the metric
$\theta \Gamma$
(observe that it is admissible for $\theta$ bounded)
gives
$$(a_{1}\circ \theta^{1/2})\sharp (a_{2}\circ \theta^{1/2})\in 
S(\lambda^{m_{1}+m_{2}},\theta \Gamma)
$$
which implies
$
\bigl((a_{1}\circ \theta^{1/2})\sharp (a_{2}\circ \theta^{1/2})\bigr)\circ \theta^{-1/2}\in 
S(\lambda^{m_{1}+m_{2}},\Gamma)$.
Applying this to the integral above gives the result of the lemma.
\qed
\enddemo
 \proclaim{Lemma A.5.3}
 Let $\Gamma$ be a positive definite quadratic form on $\RZ$ such that
 $\Gamma=\Gamma^\sigma$
 and $\lambda$ be a $\Gamma$-weight.
 Let $l\in \N, \mu\in \R$
 and $a$ be a locally bounded function defined on $\RZ$ such that
 $$\forall j\in\{0,\dots, l\},\quad \val {a^{(j)}(X)}\le C\lambda(X)^{\mu-\frac{j}{2}}. $$
 Then the function 
 $a\ast \exp -2\pi\Gamma$ belongs to 
 $S_{l}(\lambda^\mu,\lambda^{-1}\Gamma).$
\endproclaim
\demo{Proof}
We use the formula
$
(a\ast \exp -2\pi\Gamma)(X)=\int a(X-Y)\exp-2\pi \Gamma(Y) dY
$
to obtain the estimate for the derivatives of order $\le l$: we get for $k\le l$
$$\multline
\val{(a\ast \exp -2\pi\Gamma)^{(k)}(X)}\le C\lambda(X)^{\mu-\frac{k}{2}}
\int
\frac{\lambda(X-Y)^{\mu-\frac{k}{2}}}{\lambda(X)^{\mu-\frac{k}{2}}
}
\exp-2\pi \Gamma(Y) dY
\\\le
C\lambda(X)^{\mu-\frac{k}{2}}
\int
(1+\Gamma(Y))^{N\val{\mu-\frac{k}{2}}}
\exp-2\pi \Gamma(Y) dY
=C'\lambda(X)^{\mu-\frac{k}{2}},
\endmultline$$
and for $k>l$
we have
$
(a\ast \exp -2\pi\Gamma)^{(k)}=(a^{(l)}\ast (\exp -2\pi\Gamma)^{(k-l)})
$
yielding immediately the result.
\qed
\enddemo
Let us recall the composition formula in the Weyl quantization,
with the symplectic form
$[,]$
given in (1.3.2).
We have $a^wb^w=(a\sharp b)^w$ and, for $X\in \RZ$, 
$$\multline
(a\sharp b)(X)=2^{2n } \iint _{{\RZ\times \RZ}}a(Y) b(Z)
\exp-4i\pi [X-Y,X-Z] dY dZ
\\=
2^{2n } \iint _{{\RZ\times \RZ}}a(Y+X) b(Z+X)
\exp-4i\pi [Y,Z] dY dZ.
\endmultline
\tag A.5.3$$
We note also that
$$
(a\sharp b)'=a'\sharp b+a\sharp b'.
\tag A.5.4$$
Moreover, if $a$ is a function only of $\xi$,
we have
$$\align
&(a\sharp b)(x,\xi)=
2^{2n } \int _{\R^{4n}}a(\eta) b(z,\zeta)
e^{-4i\pi (\xi-\eta)(x-z)}
e^{4i\pi (x-y)(\xi-\zeta)}dy d\eta dz d\zeta
\\
&=
2^{n } \int _{\R^{2n}}a(\eta) b(z,\xi)
e^{-4i\pi (\xi-\eta)(x-z)}
d\eta dz 
\\
&=
2^{n } \int _{\R^{2n}}
((1+D_{\eta}^2/4)^N a)(\eta) b(z,\xi)
(1+\val{x-z}^2)^{-N}
e^{-4i\pi (\xi-\eta)(x-z)}
d\eta dz 
\\
&=
2^{n } \int _{\R^{2n}}
((1+D_{\eta}^2/4)^N a)(\eta) (1+\val{\xi-\eta}^2)^{-N}
(1+D_{z}^2/4)^{N}
\Bigl(
b(z,\xi)
(1+\val{x-z}^2)^{-N}
\Bigr)
\\
&\hskip255pt
e^{-4i\pi (\xi-\eta)(x-z)}
d\eta dz
\endalign$$
so that with $N\ge E(n/2)+1$
$$
\val{(a\sharp b)(x,\xi)}\le 
\max_{j\le 2N}\nuorm{a^{(j)}}_{L^\io}
\max_{j\le 2N}\nuorm{b^{(j)}}_{L^\io}(1+ \val{\xi-\supp a})^{-N/2}
c(n,N).
\tag A.5.5
$$
\subhead
A.6. The Beals-Fefferman reduction
\endsubhead
\proclaim{Lemma A.6.1}
Let 
$F:\R\rightarrow\R$
be a $C^2$
function such that
$$
16\val{F(0)}
< F'(0)^2,
\quad
\nuorm{F''}_{L^\io(\R)}
\le 1.
\tag A.6.1$$
We set
$
\rho=\val{F'(0)}/4.
$
Then there exists
$t_0\in[-\rho/2,\rho/2]$
and $e\in C^1(\R)$ such that
$$
\text{for $\val t\le \rho$},
\
F(t)=(t-t_0) e(t),
\quad
8\rho\ge e(t)\ge \rho,
\quad
\nuorm{e'}_{L^\io(\R)}\le 1.
\tag A.6.2$$
\endproclaim
\demo{Proof}
Assume first that
$F(0)=0$
and $F'(0)=4\rho$. Then,
for
$\val{t}\le 2\rho,$
$$
F(t) =te(t),\quad
6\rho\ge e(t)\ge4\rho-2\rho=2\rho,
\nuorm{e'}_{L^\io(\R)}\le 1.
$$
Now if $F(0)>0$
and $F'(0)=4\rho$,
$
F(-\frac{\rho}{2})\le\rho^2
-\frac{\rho}{2}4\rho+\frac{\rho^2}{4}
<0,
$
so that,
for some $t_0\in ]-\rho/2,0[$
we have
$
F(t_0)=0.
$
Using what was done above,
we have
for
$\val s\le\val{F'(t_0)}/2$,
$$
F(s+t_0)=(s+t_0)e_0(s),\quad
3\val{F'(t_0)}/2
\ge
e_0(s)\ge
\val{F'(t_0)}/2,\quad
\nuorm{e_0'}_{L^\io(\R)}\le 1.
$$
But
since
$$
\frac{\val{F'(t_0)}}{2}\ge
\frac{1}{2}(4\rho-\frac{\rho}{2})=
\frac{7\rho}{4}
\quad
\text{and}\quad
\frac{7\rho}{4}-\frac{\rho}{2}=
\frac{5\rho}{4}\ge \rho
$$
we have
on
$
[t_0-\frac{7\rho}{4},
t_0+\frac{7\rho}{4}]$
which contains
$
[-\rho,\rho],
$
$$
F(t)=(t-t_0)e(t),
\quad
\val{t_0}\le \rho/2,
\quad
8\rho\ge \frac{27\rho}{4}\ge 
e(t)\ge 7\rho/4\ge \rho,
\nuorm{e'}_{L^\io(\R)}\le 1.\qed
$$
\enddemo
\proclaim{Lemma A.6.2}
Let 
$F:\R^d\rightarrow\R$
be a $C^2$
function such that
$$
2^6\val{F(0)}
< \nuorm{\nabla F(0)}^2,
\quad
\nuorm{F''}_{L^\io(\R^d)}
\le 1.
\tag A.6.3$$
We set
$
\rho=\nuorm{\nabla F(0)}2^{-5}.
$
There exists
two $C^1$ functions
$\alpha:\R^{d-1}\rightarrow[-5\rho,5\rho]
\quad\text{and}\quad e:\R^d\rightarrow[7\rho,70\rho],
$
a set of orthonormal coordinates
$(x_1,x')\in\R\times \R^{d-1}$ such that
for $\max\bigl(\val{x_1},\val {x'}\bigr)\le \rho$,
$$
F(x)=\bigl(x_1+\alpha(x')\bigr) e(x),
\quad
\nuorm{e'}_{L^\io(\R^d)}\le 1,\quad
\nuorm{\alpha'}_{L^\io(\R^{d-1})}\le 1.
\tag A.6.4$$
\endproclaim
\demo{Proof}
We can choose the coordinates
so that
$
\nabla F(0)=\frac{\p  F}{\p x_1}(0)\overset
\rightarrow\to {e_1}
$.
Then for $\val{x'}\le \rho$, we have
$\val{F(0,x')}\le 2^{-6+10}\rho^2
+\rho2^5\rho+\frac{1}{2}\rho^2
=\rho^2(2^5+2^4+2^{-1})
$
and
$$
\Val{\frac{\p  F}{\p x_1}(0,x')}\ge
\Val{\frac{\p  F}{\p x_1}(0,0)}-\rho=
(2^5-1)\rho
$$
so that
$$
\frac{
16\val{F(0,x')}
}
{\Val{\frac{\p  F}{\p x_1}(0,x')}^2}
\le\frac{16\times 48.5}{31^2}<1.
$$
Applying the lemma A.5.1, we get 
for all $\val{x'}\le \rho$
the existence
of $\alpha(x')$
such that,
when
$\val{x_1}\le 31\rho/4$
$$
F(x_1,x')=\bigl (x_1+\alpha(x')\bigr) e(x),
\quad
\val{\alpha(x')}\le\frac{33\rho}{8}<5\rho,\quad
70\rho\ge 8\times 33 \rho/4\ge \val{e(x)}\ge 31 \rho/4>7\rho.
$$
The implicit function theorem
guarantees
the $C^1$
regularity of the function $\alpha$
and the
Taylor-Lagrange
formula with integral remainder provides
the
regularity of $e$.\qed
\enddemo
\remark{Remark A.6.3}
If the function $F$ in the lemma A.6.2
is $\moo$, since the function $\alpha$ is obtained by the implicit function theorem, and $e$ by Taylor's formula with integral remainder, both function $\alpha,e$ are $\moo$.
Moreover, the identity $F(-\alpha(x'),x')=0$ implies that
$$
\val{\alpha^{(k)}(x')}\le C_{k}\rho^{1-k},\quad
\val{e^{(k)}(x')}\le C_{k}\rho^{-k}
$$
where $C_{k}$ are semi-norms of the function $F$ in $\max(\val{x_{1}},\val{x'})\le \rho$.
In particular, if we apply this result to the function (2.1.21)
$$
F(T)= \Lambda^{1/2}q\bigl(t,Y+\nu(t,Y)^{1/2} T\bigr)\mu(t,Y)^{-1/2}\nu(t,Y)^{-1}
$$
we get that
$
 \val {F^{(k)}}$
 is bounded above by $\gamma_{k}(q)$
 and $1/2\le \rho\le \gamma_{1}(q)$
 as defined in (2.1.1).
 We get then from the lemma A.6.2
 $$
  \Lambda^{1/2}q\bigl(t,Y+\nu(t,Y)^{1/2} T\bigr)\mu(t,Y)^{-1/2}\nu(t,Y)^{-1}= e_{0}(T)
  (T_{1}+\alpha_{0}(T'))
 $$
 so that $e_{0},\alpha_{0}$ are smooth with fixed bounds and thus
 $$
  \Lambda^{1/2}q\bigl(t,X\bigr)\mu(t,Y)^{-1/2}  
  =e_{0}\bigl((X-Y)\nu(t,Y)^{-1/2}\bigr)
  \nu(t,Y)^{1/2}\Bigl(X_{1}-Y_{1}+\alpha_{0}\bigl((X'-Y')\nu(t,Y)^{-1/2}\bigr)\nu(t,Y)^{1/2}
  \Bigr)
 $$
  which corresponds exactly to (2.1.22-23-24).
\endremark
\subhead
A.7. On tensor products of homogeneous functions
\endsubhead
Let $n\ge 1$ be an integer and $N=n+1$.
Let $(y_0;\eta_0)\in \R^{N}\times \Bbb S^{N-1}$
such that 
$$
(y_0;\eta_0)=(t_0,x_0; \tau_0,\xi_0)\in \R\times\R^n\times \R\times\R^n,\quad\text{with $\tau_0=0,\ \xi_0\in \Bbb S^{n-1}$.}
$$
Let $r\in]0,1/4]$ be given.
There
exists
a function $\chi_0\in\moo(\R^N; [0,1])$ 
such that for $\lambda\ge 1$ and $\eta \in \R^N$ with $\val \eta\ge 1$,
we have
$\chi_0(\lambda \eta)=\chi_0(\eta)$
(``homogeneity of degree zero outside the unit ball")
and
$$
\chi_0(\tau,\xi)=\cases
1&\text{if $\tau^2+\val\xi^2\ge 1$ and $\val{\tau} \le r\val\xi$,        }
\\
0&\text{if $\tau^2+\val\xi^2\le 1/4$ or $\val{\tau} \ge 2r\val\xi$.}
\endcases
$$
There exists a function $\psi_0\in\moo(\R^n; [0,1])$
such that
for $\lambda\ge 1$ and $\xi \in \R^n$ with $\val \xi\ge 1$,
we have
$\psi_0(\lambda \xi)=\psi_0(\xi)$
and, $$
\psi_0(\xi)=\cases
1&\text{if $\val\xi\ge 1$ and $\val{\frac{\xi}{\val \xi}-\xi_0}\le r$,        }
\\
0&\text{if $\val\xi\le 1/2$
or $\val{\frac{\xi}{\val \xi}-\xi_0}\ge 2r$
.}
\endcases
$$
We define the function $\Phi_0$ by
$$
\Phi_0(\tau, \xi)=\chi_0(\tau,\xi)\psi_0(\xi).
\tag A.7.1$$
\proclaim{Lemma A.7.1}
The function $\Phi_{0}$ is such that
for $\lambda\ge 1$ and $\eta \in \R^N$ with $\val \eta\ge 2$,
we have
$\Phi_0(\lambda \eta)=\Phi_0(\eta)$. Moreover,
with
$\eta_{0} =(0,\xi_{0})$, we have 
$$
\text{$\Phi_{0}(\eta)=1$ for $\val \eta\ge 2$ and }
\Val{\frac{\eta}{\val \eta}-\eta_{0}}\le r/2,
\qquad
\text{$\Phi_{0}(\eta)=0$ for $\val \eta\ge 2$ and }
\Val{\frac{\eta}{\val \eta}-\eta_{0}}\ge 4r.
$$
\endproclaim
\demo{Proof}
The function
$
\Phi_0$
is such that
for $\lambda\ge 1$ and $\eta \in \R^N$ with $\val \eta\ge 2$,
we have
$\Phi_0(\lambda \eta)=\Phi_0(\eta)$:
in fact, if
$\tau^2+\val\xi^2\ge 4$ and $\val \tau\le 2r \val\xi$,
we get
$\val\xi^2\ge  4(1+4r^2)^{-1}\ge 1$,
so that
$\psi_0(\lambda \xi)=\psi_0(\xi)$
and since we have also in that case
$\chi_0(\lambda \eta)=\chi_0(\eta)$,
we get the sought property.
 Now if 
 $\tau^2+\val\xi^2\ge 4$ and $\val \tau> 2 r \val\xi$,
 we see that
 $\chi_0(\lambda\tau,\lambda\xi)=\chi_0(\tau,\xi)=0$
 so that, 
 $\Phi_0(\lambda \eta)=0=\Phi_0(\eta).$
 Moreover, if $\tau^2+\val\xi^2\ge 4$
 and
$$
\frac{\tau^2}{\tau^2+\val\xi^2}+\Val{\frac{\xi}{(\tau^2+\val\xi^2)^{1/2}}-\xi_0}^2\le r^2/4,
$$
we get that
$
\val\tau\le r\val \xi(4-r^2)^{-1/2}\le r\val \xi
$
and thus $\chi_0(\tau,\xi)=1$;
also this implies $\val \xi\ge 2(1+r^2)^{-1/2}\ge 1$,
so that
$\psi_0(\xi)=\psi_0(\xi/\val\xi)$.
We have then
$$
\Val{\frac{\xi}{\val\xi}-\xi_0}\le \frac{r}{ 2}+\Val{
\frac{\xi}{\val\xi}
-
\frac{\xi}{(\tau^2+\val\xi^2)^{1/2}}
}
\le
\frac{r}{2}+\val\xi\val\tau^2\val\xi^{-3}
\le \frac{r}{2}+\frac{r^2}{4-r^2}\le r,
$$
which implies
$\psi_0(\xi)=\psi_0(\xi/\val\xi)=1$,
so that
$\Phi_0$ is equal to 1 on a conic neighborhood
of
$(0,\xi_0)$ in $\R^N$ minus a ball.
Similarly, if $\tau^2+\val \xi^2\ge 4$
and
$$
\frac{\tau^2}{\tau^2+\val\xi^2}+\Val{\frac{\xi}{(\tau^2+\val\xi^2)^{1/2}}-\xi_0}^2\ge 16r^2,
$$
either $\val\tau\ge2r\val\xi$ and $\chi_0(\tau,\xi)=0$,
entailing
$\Phi_0(\tau,\xi)=0$
or
$\val\tau\le2r\val\xi$ 
and then
$\Val{\frac{\xi}{(\tau^2+\val\xi^2)^{1/2}}-\xi_0}^2\ge 12r^2$
and
$\val{\xi}\ge 2(1+4r^2)^{-1/2}\ge 1$
so that
$\psi_0(\xi)=\psi_0(\xi/\val\xi)$.
In this case, we have
$$
\Val{\frac{\xi}{\val\xi}-\xi_0}\ge 2\sqrt 3r-\Val{
\frac{\xi}{\val\xi}
-
\frac{\xi}{(\tau^2+\val\xi^2)^{1/2}}
}\ge 2\sqrt 3r-\frac{\tau^2}{\val \xi^2}\ge
2\sqrt 3r -4r^2\ge 2r,
$$
implying $\psi_0(\xi)=0$ and thus $\Phi_0(\tau, \xi)=0$.
Eventually, we have proven that 
$\Phi_0$ is also supported in a conic neighborhood of
$(0,\xi_0)$ in $\R^N$.
\qed
\enddemo
\subhead
A.8. Composition of symbols
\endsubhead
\proclaim{Lemma A.8.1}
Let $G, g$ be the metrics on $\R^{2N}$ defined in $(4.3.8)$
and let $s_{1},s_{2}$ be two real numbers.
Let $a$ be a symbol in $S(\langle \xi\rangle^{s_{1}}, g)$ and $b$ be a symbol in 
$S(\langle \xi, \tau\rangle^{s_{2}}, G)$
such that
$\supp b\subset Z_{C}=\{(t,x,\tau,\xi)\in \R^{2N}, \val\tau\le 1+C\val \xi\}$.
Then the symbols
$a\sharp b, b\sharp a, a\circ b, b\circ a$
belong to 
$S(\langle \xi, \tau\rangle^{s_{1}+s_{2}}, G)$ and are essentially supported in 
$Z_{C}$, i.e. are the sum of a symbol of $
S(\langle \xi, \tau\rangle^{s_{1}+s_{2}}, G)$
supported in $Z_{C}$
and of a symbol in 
$S(\langle \xi, \tau\rangle^{-\io}, G)=\cap_{N}S(\langle \xi, \tau\rangle^{-N}, G)$.
\endproclaim
\demo{Proof}
 We have
 $$
 (a\circ b)(t,x,\tau,\xi)=
 \int e^{-2i\pi(s\sigma+y\eta)}
 a(t,x,\tau+\sigma,\xi+\eta) b(t+s, x+y, \tau,\xi) dsd\sigma dy d\eta,
 \tag A.8.1 $$
 so that,
 using the standard expansion of the symbols and the fact that $b$ is supported in $Z_{C}$,
 $$
 \multline 
 a\circ b=
 \sum_{\val \alpha<\nu}\frac{1}{\alpha !}
 \overbrace{D_{\tau,\xi}^\alpha a\ \p_{t,x}^\alpha b}^{\in S(\langle\tau,\xi\rangle^{1-\val \alpha},G)}
 \\
 +
 \int_0^1 \frac{(1-\theta)^{\nu-1}}{(\nu-1)!}e^{-2i\pi(s\sigma+y\eta)}
 D_{\tau,\xi}^\nu a(t,x,\tau+\theta\sigma,\xi+\theta\eta) 
 \p_{t,x}^\nu b(t+ s, x+y, \tau,\xi) dsd\sigma dy d\eta d\theta.
 \endmultline
 $$ 
 We define
 $$
  I_{\theta}(\tau ,\xi)=
\int e^{-2i\pi(s\sigma+y\eta)}
 D_{\tau,\xi}^\nu a(t,x,\tau+\theta\sigma,\xi+\theta\eta) 
 \p_{t,x}^\nu b(t+ s, x+y, \tau,\xi) dsd\sigma dy d\eta\tag A.8.2
 $$
 and integrating by parts, we obtain for all nonnegative even integers $m$ that
 $$\multline
 I_{\theta}(\tau ,\xi)=
 \int e^{-2i\pi(s\sigma+y\eta)}
 \langle \sigma \rangle^{-m}
 \langle D_s\rangle^m
 \langle s \rangle^{-m}
 \langle D_\sigma\rangle^m
  \langle y \rangle^{-m}
 \langle D_\eta\rangle^m
 \langle \eta \rangle^{-m}
 \langle D_y\rangle^m
  \\
  D_{\tau,\xi}^\nu a(t,x,\tau+\theta\sigma,\xi+\theta\eta) 
 \p_{t,x}^\nu b(t+s, x+ y, \tau,\xi)
   dsd\sigma dy d\eta,
 \endmultline$$
 and consequently
 $$
 \val{ I_{\theta}(\tau ,\xi)
}\lesssim \int
\langle \sigma \rangle^{-m}
 \langle s \rangle^{-m}
  \langle y \rangle^{-m}
 \langle \eta \rangle^{-m}
 (1+\val{\xi+\theta \eta})^{s_{1}-\nu}  dsd\sigma dy d\eta(1+\val{\xi}+\val \tau)^{s_{2}}\bold 1(\val \tau\lesssim \val \xi).
 $$
 In the integrand, when $\val \eta\le\val \xi/2$,
 we get, since $\theta\in [0,1]$, 
 $
 \val{\xi+\theta\eta}\ge \val \xi-\val \eta\ge \val \xi/2.
 $
 As a result, we get for this part of the integral the estimate
 $$(1+\val \xi)^{\val {s_{1}}-\nu}
 (1+\val{\xi}+\val \tau)^{s_{2}}
 \bold 1(\val \tau\lesssim \val \xi)
 \lesssim
  (1+\val{\xi}+\val \tau)^{-\nu/2}, \quad\text{for $\nu$ large enough.}
 $$
 When $\val \eta>\val \xi/2$,
 we use the term $\langle\eta \rangle^{-m}$
 and the estimate
 $$
 (1+\val \xi)^{-m/2} (1+\val{\xi}+\val \tau)^{s_{2}}
 \bold 1(\val \tau\lesssim \val \xi)
 \lesssim
  (1+\val{\xi}+\val \tau)^{-m/4}, \quad\text{for $m$ large enough.}
 $$
To check that the derivatives of $I_{\theta}$
 will satisfy the expected estimates,
 we differentiate the expression (A.8.2)
 and repeat the previous proof.
 We know now that, for all $\nu\ge N_{0}$
 $$
a\circ b=
 \sum_{\val \alpha<\nu}\frac{1}{\alpha !}
 D_{\tau,\xi}^\alpha a\ \p_{t,x}^\alpha b + r_{\nu}, \quad r_{\nu}\in
 S(\langle\tau,\xi\rangle^{-\nu/2},G).
$$
Using the standard Borel argument, we find 
$c\in S(\langle\tau,\xi\rangle^{s_{1}+s_{2}},G)$,
essentially supported in $Z_{C}$ such that, for all $\nu$
$$
c-\sum_{\val \alpha<\nu}\frac{1}{\alpha !}
 D_{\tau,\xi}^\alpha a\ \p_{t,x}^\alpha b\in  S(\langle\tau,\xi\rangle^{s_{1}+s_{2}-\nu},G),
$$
entailing that, for all $\nu\ge N_{0}$,
$$
a\circ b-c=-c+
\sum_{\val \alpha<\nu}\frac{1}{\alpha !}
 D_{\tau,\xi}^\alpha a\ \p_{t,x}^\alpha b+r_{\nu}\in S(\langle\tau,\xi\rangle^{\max(-\nu/2,s_{1}+s_{2}-\nu)},G),
$$
implying that
$a\circ b-c\in
S(\langle\tau,\xi\rangle^{-\io},G)$,
which gives the result of the lemma for $a\circ b$.
To get the result for $b\circ a$ is somewhat easier by looking at (A.8.2), to obtain the estimate
$$\multline
 \val{ I_{\theta}(\tau ,\xi)
}\lesssim \int
\langle \sigma \rangle^{-m}
 \langle s \rangle^{-m}
  \langle y \rangle^{-m}
 \langle \eta \rangle^{-m}
 (1+\val{\xi+\theta \eta}+\val{\tau+\theta\sigma})^{s_{2}-\nu}  
 \bold 1({\val{\tau+\theta \eta}\lesssim
 \val{\xi+\theta \sigma}
  })
  \\dsd\sigma dy d\eta(1+\val{\xi})^{s_{1}}.
  \endmultline$$
 When $\val \tau\lesssim \val \xi$ the discussion is the same as for $a\circ b$.
 When $\val \tau\gg  \val \xi$, we split the integral in two parts:
  the region where
  $\val \sigma\le\val \tau/2$, in which we get negative powers of $(1+\val \tau)$
  from the term with the exponent $s_{2}-\nu$, and the region
  where
  $\val \sigma>\val \tau/2$ 
  in which we use the term
  $\langle \sigma \rangle^{-m}$.
  The last part of the discussion is the same.
  To obtain the result for $a\sharp b$
  (which will give also $b\sharp a$ since $\overline{a\sharp b}=\bar b\sharp \bar a$),
  we use the group $J^t=\exp 2i\pi tD_{x}D_{\xi}$ and the formula
  $
  a\sharp b=J^{-1/2}\bigl(J^{1/2}a\circ J^{1/2} b\bigr).
$
  Using the assumptions of the lemma, we see that $J^{1/2}a$ satisfies the same hypothesis as $a$
  and 
  $J^{1/2}b$ is essentially supported in $Z_{C}$. The proofs above give thus that
  $J^{1/2}a\circ J^{1/2} b$ satisfies the conclusion of the lemma,
  which is ``stable" by the action of $J^{-1/2}$.
  The proof of the lemma A.8.1 is complete.
 \qed
\enddemo

\enddocument
\bye